\documentclass[leqno,11pt]{article}  
\usepackage{amsmath}
\usepackage{amssymb}

\usepackage[german,english]{babel}
\usepackage{amsthm}
\usepackage{xypic}  

\author{\textsc{Elmar Grosse-Kl\"onne}}
\date{}

\theoremstyle{plain} 
\newtheorem{satz}{Theorem}[section]  
\newtheorem{kor}[satz]{Corollary}  
\newtheorem{lem}[satz]{Lemma}  
\newtheorem{pro}[satz]{Proposition}  
 
\theoremstyle{remark}

\theoremstyle{definition}

\begin{document}

%

\title{ On the universal module of $p$-adic spherical Hecke algebras}

\begin{center}{\bf On the universal module of $p$-adic spherical Hecke algebras}\\Elmar Gro\ss{}e-Kl\"onne
\end{center}

\begin{abstract} Let $\widetilde{G}$ be a split connected reductive group with
  connected center $Z$ over a local
non-Archimedean field $F$ of residue characteristic $p$, let $\widetilde{K}$ be a hyperspecial
maximal compact open subgroup in $\widetilde{G}$. Let $R$ be a commutative ring, let $V$ be a finitely
generated $R$-free $R[\widetilde{K}]$-module. For an $R$-algebra $B$ and a character $\chi:{\mathfrak
    H}_V(\widetilde{G},\widetilde{K})\to B$ of the
spherical Hecke algebra ${\mathfrak
    H}_V(\widetilde{G},\widetilde{K})={\rm End}_{R[\widetilde{G}]}{\rm
  ind}_{\widetilde{K}}^{\widetilde{G}}(V)$ we consider the specialization $$M_{\chi}(V)={\rm
  ind}_{\widetilde{K}}^{\widetilde{G}}V\otimes_{{\mathfrak
    H}_V(\widetilde{G},\widetilde{K}),\chi}B$$ of the universal ${\mathfrak
    H}_V(\widetilde{G},\widetilde{K})$-module ${\rm
  ind}_{\widetilde{K}}^{\widetilde{G}}(V)$. For large classes of $R$ (including
${\mathcal O}_F$ and $\overline{\mathbb F}_p$), $V$, $B$ and $\chi$, arguing geometrically on
the Bruhat Tits building we give a sufficient criterion for $M_{\chi}(V)$ to
be $B$-free and to admit a $\widetilde{G}$-equivariant
resolution by a Koszul complex built from finitely many copies
of ${\rm
  ind}_{\widetilde{K}Z}^{\widetilde{G}}(V)$. This criterion is the exactness of certain
fairly small and explicit ${\mathfrak N}$-equivariant $R$-module complexes, where
${\mathfrak N}$ is the group of ${\mathcal O}_F$-valued points of the unipotent
radical of a Borel subgroup in $\widetilde{G}$. We verify it if $F={\mathbb
  Q}_p$ and if $V$ is an irreducible
$\overline{\mathbb F}_p[\widetilde{K}]$-representation with highest weight in the (closed) bottom
$p$-alcove, or a lift of it to ${\mathcal O}_F$. We use this to construct $p$-adic
integral structures in certain locally algebraic representations of $\widetilde{G}$.

\end{abstract}


{\it AMS 2010} Mathematics subject classification: 22E50, 20G25, 20E42

\tableofcontents

\section{Introduction}
 Let $F$ be a local non-Archimedean field with finite residue field $k_F$ of
characteristic $p>0$. Let $\widetilde{G}=\widetilde{\bf
  G}(F)$ be the group of $F$-rational points of a split connected reductive algebraic group $\widetilde{\bf
  G}$ over ${\mathcal O}_F$. Let $\widetilde{K}=\widetilde{\bf
  G}({\mathcal O}_F)$. \\

{\bf (1) The results.} 

(1a) {\it Lattices in smooth/locally algebraic $\widetilde{G}$-representations.} The general problem of finding $\widetilde{G}$-invariant norms on smooth $\widetilde{G}$-representations over $p$-adic fields is of vital interest for the so
called $p$-adic local Langlands program. Certain quite obvious and handy {\it necessary} conditions for the existence of $\widetilde{G}$-invariant norms on smooth (or, more generally, locally
algebraic) $\widetilde{G}$-representations over $p$-adic fields are given
e.g. in \cite{st}, \cite{eme}. However, it is widely unknown if these
conditions are also sufficient. The first non-trivial (smooth) case in which
these conditions were shown to be sufficient was that of unramified smooth
principal series representations of ${\rm GL}_2({\mathbb Q}_p)$, see
\cite{berbre}. Later on, more general results for tamely ramified smooth
principal series representations of ${\rm GL}_2(F)$ were obtained in
\cite{vig} using the Bruhat-Tits tree of ${\rm GL}_2(F)$. However, even for
${\rm GL}_2({\mathbb Q}_p)$ the general problem of constructing ${\rm
  GL}_2({\mathbb Q}_p)$-invariant norms in smooth principal series
representations (when they are predicted to exist) remains open. Moreover, for
groups $\widetilde{G}$ other than ${\rm GL}_2(F)$ we are not aware of a single
non-trivial example where $\widetilde{G}$-invariant norms on smooth
$\widetilde{G}$-representations have been constructed. In this paper we
construct $\widetilde{G}$-invariant norms on unramified smooth principal
series representations of $\widetilde{G}$ whenever the aforementioned
necessary conditions are fulfilled --- under the additional hypothesis that
$F={\mathbb Q}_p$, that $\widetilde{G}$ has connected center, and that the
Coxeter number of ${\widetilde{G}}$ is $p$-small (more precisely, the Coxeter
numbers of the irreducible components of $\widetilde{G}$ are at most
$p+1$). See Theorem \ref{normexis} for the precise (and in fact more general)
statement. As explained in Section 9 one may use Proposition \ref{bettercond}, Corollary \ref{rank2}, Theorem
\ref{modest} to construct $\widetilde{G}$-invariant norms in other, more
general locally algebraic representations whose algebraic parts 'belong to the bottom alcove'. 

(1b) {\it Freeness of the universal module of the spherical Hecke algebra.} A
central object for the representation theory of $\widetilde{G}$ is the
spherical Hecke algebra, formed with respect to $\widetilde{K}$. For example,
the smooth representation theory of ${\rm GL}_2(F)$ on $\overline{\mathbb
  F}_p$-vector spaces can at present be approached only through restriction of
${\rm GL}_2(F)$-representations to $\widetilde{K}$, and hence through the
spherical Hecke algebra. This theory had been initiated by Barthel and
Livn\'{e} who proved that the universal module over this algebra is
free. Except for the result of Bella${{\ddot{\i}}}$che and Otwinowska (on the
spherical algebra of ${\rm PGL}_3(F)$ with trivial coefficients, see
\cite{bo}, also for further motivations and applications) we are not aware of
other groups $\widetilde{G}$ for which the universal modules over the
spherical Hecke algebra have been shown to be free. In this paper we prove
that the the universal module over the spherical Hecke algebra with trivial
coefficients is free --- under the additional hypothesis that $F={\mathbb
  Q}_p$, that $\widetilde{G}$ has connected center, and that the Coxeter
number of ${\widetilde{G}}$ is $p$-small. More generally, we treat non-trivial
algebraic coefficients which 'belong to the bottom alcove', and non-trivial
central characters. See Proposition \ref{bettercond}, Corollary \ref{rank2}, Theorem
\ref{modest} for the precise (and in fact more general) statements. 

(1c) {\it Exact resolutions of $\widetilde{G}$-representations over
  $\overline{\mathbb F}_p$.} An important technique in the smooth
representation theory of $\widetilde{G}$ on complex vector spaces,
systematically introduced by Schneider and Stuhler \cite{ss}, is to 'spread
out' smooth $\widetilde{G}$-representations as coefficient systems on the
Bruhat-Tits building of $\widetilde{G}$, assigning in an economic way to each
(poly)simplex a certain (small) representation of its stabilizer. The point of
this procedure is that it typically leads to {\it acyclic}
$\widetilde{G}$-equivariant complexes and thereby provides resolutions of
smooth admissible $\widetilde{G}$-represensentations by compactly induced
finitely generated representations. Later on these techniques have been
generalized by Vign\'{e}ras \cite{vig0} to other coefficient fields in which $p$ is invertible. On the other hand, if one considers smooth representations over coefficient fields of characteristic $p$ then the complexes constructed analogously are no longer exact in general, even for ${\rm GL}_2(F)$ (although for ${\rm GL}_2(F)$ to some extent these techniques remain useful, see \cite{vig}, but the proofs, trivial for complex coefficients, become very painful). Thus, at present there is no method available to produce resolutions of smooth (admissible) $\widetilde{G}$-representations over $\overline{\mathbb F}_p$ by compactly induced finitely generated representations. In fact, for groups $\widetilde{G}$ other than ${\rm GL}_2(F)$ we are not aware of a single non-trivial example of an exact complex of compactly induced finitely generated $\widetilde{G}$-representations over $\overline{\mathbb F}_p$ which had been constructed in the literature so far. In this paper, for $\widetilde{G}$ as above we construct resolutions for certain very natural smooth $\widetilde{G}$-represensentations over $\overline{\mathbb F}_p$ (e.g. reductions modulo $p$ of the lattices appearing in (1a)). More generally, we introduce a technique for verifying the exactness of a given sequence of representations of a $p$-adic group on $\overline{\mathbb F}_p$-vector spaces which we think to be of some more general interest. See Proposition \ref{bettercond}, Corollary \ref{rank2}, Theorem
\ref{modest} for the precise (and in fact more general) statements.\\

{\bf (2) The methods.} The following techniques and constructions appear to be new.

(2a) {\it Vertex coefficient systems.} We introduce a notion of 'vertex
coefficient systems' on the (semisimple) Bruhat-Tits building $X$ of $\widetilde{G}$ (see
section 3). They bear
some similarity with the coefficient systems studied in \cite{ss}, but are
different and the relation remains unclear (this deserves further
investigation). Vertex
coefficient systems are supported on the special vertices only, while the
coefficient systems studied in \cite{ss} are supported on all the
simplices. In fact, for our purposes it is sufficient and convenient to define
vertex
coefficient systems as living on the set $A^0$ of special vertices of a chosen
apartment $A$ in $X$, but it is clear that this notion is a 'projection to $A$'
of a notion of coefficient systems on $X$. Depending on an additional set of parameters $\alpha$ there is a
natural Koszul complex $C_{\alpha}{\mathcal
  F}^{\bullet}$ associated with a
vertex coefficient system ${\mathcal F}$. Proposition \ref{fide} provides a filtration of $C_{\alpha}{\mathcal
  F}^{\bullet}$, together with a
description of the graded pieces. As these graded pieces are direct sums of complexes
supported only on small areas of $A$ we obtain a local criterion for $C_{\alpha}{\mathcal
  F}^{\bullet}$ to be exact. All this is very general, based just on combinatorial geometry of affine Coxeter complexes. It is not specific for $p$-adic or $p$-modular representation theory.

(2b) {\it Cohomology of ${\mathfrak N}$: an analog of Kostant's theorem.} Let ${\mathfrak
  N}$ denote the group of ${\mathcal O}_F$-valued points of the unipotent
radical of a Borel subgroup in $\widetilde{G}$. In our cases of interest we are dealing with certain vertex coefficient
systems consisting of $\overline{\mathbb F}_p$-vector space. The said local exactness criterion reduces us to investigating the exactness
of certain complexes which are equivariant under ${\mathfrak
  N}$. We show that this exactness can be deduced from a result, if available,
on the group cohomology of ${\mathfrak N}$ with coefficients in algebraic
representations of $\widetilde{\bf G}(k_F)$. More precisely, what is needed is
a certain analog of Kostant's famous result on the corresponding Lie algebra
cohomology. We prove it here if $F={\mathbb Q}_p$ and if the algebraic
representation has highest weight in the' bottom alcove', see Theorem \ref{prodpoloclou}. For the trivial representation, this condition becomes the above condition that the Coxeter number of $\widetilde{\bf G}$ be $p$-small. One certainly should expect that these group cohomology computations generalize to finite extensions $F$ of ${\mathbb Q}_p$, and in this way all the results stated for $F={\mathbb Q}_p$ in this paper will generalize to other $F$.\\

{\bf (3) The strategy.} Let $V$ be a representation
of $\widetilde{K}$ on a free $R$-module for some ring $R$ and
consider its compact induction to $\widetilde{G}$, i.e. the
$R$-module ${\rm
  ind}_{\widetilde{K}}^{\widetilde{G}}V$ of compactly supported functions $\widetilde{G}\to
  V$ with $f(gk)=k^{-1}f(g)$ for $g\in
  \widetilde{G}$ and $k\in\widetilde{K}$ on which $\widetilde{G}$ acts by left
  translations. The corresponding spherical Hecke algebra is the endomorphism algebra ${\mathfrak
    H}_V(\widetilde{G},\widetilde{K})={\rm End}_{R[\widetilde{G}]}{\rm
  ind}_{\widetilde{K}}^{\widetilde{G}}V$. It is desirable to understand the ${\mathfrak
    H}_V(\widetilde{G},\widetilde{K})$-module structure of ${\rm
  ind}_{\widetilde{K}}^{\widetilde{G}}V$; and notably to understand if this
module structure is free,
or at least if its specializations $$M_{\chi}(V)={\rm
  ind}_{\widetilde{K}}^{\widetilde{G}}V\otimes_{{\mathfrak
    H}_V(\widetilde{G},\widetilde{K}),\chi}B$$ are $B$-free, for (suitable) $R$-algebras $B$
and $R$-algebra homomorphisms $\chi:{\mathfrak
    H}_V(\widetilde{G},\widetilde{K})\to B$. In this
paper we consider this problem under the additional assumption that $\widetilde{G}$ has
{\it connected center} $Z$ --- which we assume from now on. Moreover, we
restrict our attention to the following cases:\\

(I) $R$ is a field or a principal ideal domain and $V$ is $R$-free of rank one, or

(II) $R=\overline{{\mathbb F}}_p$ and $V$ is an irreducible smooth
$\overline{{\mathbb F}}_p[\widetilde{K}]$-module, or

(III) $R={\mathcal O}_F$ where now $F$ is a finite extension of ${\mathbb
  Q}_p$, and $V$ is a $\widetilde{\bf G}$-stable ${\mathcal O}_F$-lattice in an irreducible rational
$\widetilde{\bf G}_{F}$-representation such that $V\otimes_{{\mathcal O}_F}\overline{{\mathbb F}}_p$ is an
irreducible rational $\widetilde{\bf G}_{k_F}$-representation.\\
 
We fix an $R$-valued character $\theta$ of the subalgebra
  ${\mathfrak H}(V)$ of ${\mathfrak
    H}_V(\widetilde{G},\widetilde{K})$ corresponding to the $Z$-cocharacters of $X_*({\widetilde{\bf
      T}})_+$ (see below). We consider only those $\chi:{\mathfrak
    H}_V(\widetilde{G},\widetilde{K})\to B$ whose restrictions to ${\mathfrak
    H}(V)$ factor through $\theta$. Using $\theta$ we extend the
  $\widetilde{K}$-action on $V$ to a $\widetilde{K}Z$-action and form the compact (modulo $Z$) induction ${\rm
  ind}_{\widetilde{K}Z}^{\widetilde{G}}V$. Write $V_B=V\otimes_RB$. Instead of just studying the common cokernel $M_{\chi}(V)={\rm
  ind}_{\widetilde{K}Z}^{\widetilde{G}}V_B/\sum_{T\in  {\mathfrak
    H}_V(\widetilde{G},\widetilde{K})}(T-\chi(T))$ of the operators $T-\chi(T)$
on ${\rm
  ind}_{\widetilde{K}Z}^{\widetilde{G}}V_B$ it is more natural to consider an entire Koszul complex $$C_{\alpha}{\mathcal F}_{V_B}^{\bullet}=[\ldots\longrightarrow C_{\alpha}{\mathcal F}_{V_B}^{2}\longrightarrow C_{\alpha}{\mathcal F}_{V_B}^{1}\longrightarrow C_{\alpha}{\mathcal
  F}_{V_B}^{0}\longrightarrow
M_{\chi}(V)\longrightarrow0]$$$$C_{\alpha}{\mathcal F}_{V_B}^j={\rm
  ind}_{\widetilde{K}Z}^{\widetilde{G}}V_B\otimes\bigwedge^j{\mathbb
  Z}[\nabla]$$
whose differentials are (sums of) the $T-\chi(T)$, for certain standard Hecke
operators $T\in {\mathfrak
    H}_V(\widetilde{G},\widetilde{K})$. (Here $\nabla$ is a
finite set in bijection with the standard Hecke
operators, and ${\mathbb Z}[\nabla]$ denotes the free ${\mathbb Z}$-module on
the set $\nabla$. The subscript $\alpha=\alpha(\chi)$ indicates the dependence
on $\chi$.) The existence of such a Koszul complex and, more precisely, that it can be
viewed as the Koszul complex of a vertex coefficient system ${\mathcal
  F}_{V_B}$ on $X$, follows from the structure theorems (Satake type isomorphism, see below) on ${\mathfrak
    H}_V(\widetilde{G},\widetilde{K})$.

Therefore, as indicated above, to prove that $M_{\chi}(V)$ is free and that
 $C_{\alpha}{\mathcal F}_{V_B}^{\bullet}$ is a resolution of it, it is enough to prove that certain
subquotient complexes $\overline{C}{\mathcal F}_{V_B,\epsilon}^{\bullet}$ of
$C_{\alpha}{\mathcal F}_{V_B}^{\bullet}$, direct sums of complexes living on small areas of $X$, are exact. Similar to constructions in the papers \cite{bali} and \cite{bo}, the trick is that the $\overline{C}{\mathcal
  F}_{V_B,\epsilon}^{\bullet}$ are independent of $\alpha$ and that for our
purposes it is enough to prove
their exactness if $B=R$, $V_B=V$.

We work out sufficient criteria for the
exactness of the complexes $\overline{C}{\mathcal F}_{V,\epsilon}^{\bullet}$. Let ${\mathfrak N}$ denote the intersection of $\widetilde{K}$ with the unipotent
radical of a Borel subgroup of $\widetilde{G}$. Then in case (II) (resp. in case (III)) these exactness
criteria for the complexes $\overline{C}{\mathcal F}_{V,\epsilon}^{\bullet}$ are given in terms of (continuous)
cohomology of the group ${\mathfrak N}$, with coefficients in $V$ (resp. in
$V\otimes_{{\mathcal O}_F}\overline{{\mathbb F}}_p$), see Propostion \ref{bettercond}. \\

{\bf (4) Outline.} To begin these investigations one has to understand the $R$-algebra structure
of ${\mathfrak H}_V(\widetilde{G},\widetilde{K})$. The complete solution of
this problem, given by Schneider and
Teitelbaum \cite{st} in case (I) and then in cases (II) and (III) by Herzig
\cite{her} and finally Henniart/Vign\'{e}ras \cite{hv}, is a Satake type isomorphism which identifies ${\mathfrak
    H}_V(\widetilde{G},\widetilde{K})$ with the commutative $R$-algebra $R[X_*({\widetilde{\bf
      T}})_+]$. Here $X_*({\widetilde{\bf
      T}})_+$ is the monoid of dominant coweights with respect to a fixed
  maximal split torus in $\widetilde{G}$. Our assumption that $Z$ be connected implies that $X_*({\widetilde{\bf
      T}})_+$ is a direct sum of copies of ${\mathbb Z}$ on the one hand --- corresponding to the
  cocharacters of $Z$ --- and the commutative monoid on the set $\nabla$ of fundamental
  coweights on the other hand. All this is explained in section 2; in
  particular, Propositions \ref{redsat} and \ref{prredsat} recall the
  statements on the Satake type isomorphism. We then explain that ${\rm
  ind}_{\widetilde{K}Z}^{\widetilde{G}}V$ can naturally be viewed as a direct
sum, indexed by the set $X^0$ of special vertices of $X$, of representations of the stabilizers of these vertices.  

In section 3 we introduce the general concept of a vertex coefficient system
${\mathcal F}$ on the affine Coxeter complex $A$ associated with a split
maximal torus in $G$. By definition such an ${\mathcal F}$ comes equipped with
certain standard Hecke operators, and specializing eigenvalues ${\alpha}$ for them leads to an associated Koszul complex $C_{\alpha}{\mathcal F}^{\bullet}={\mathcal F}(A^0)\otimes\bigwedge{\mathbb Z}[\nabla]$. We use the Euclidean structure of $A$ to construct a certain filtration $\{C_{\alpha}{\mathcal F}_{\epsilon}^{\bullet}\}_{\epsilon\in{\mathbb R}}$ of it. The support in $A$ of the graded pieces $\overline{C}{\mathcal F}_{\epsilon}^{\bullet}$ can be precisely described.

In section 4 we explain how these concepts can be applied to ${\rm
  ind}_{\widetilde{K}Z}^{\widetilde{G}}V$ and produce the Koszul complex $C_{\alpha}{\mathcal F}_{V_B}^{\bullet}={\rm
  ind}_{\widetilde{K}Z}^{\widetilde{G}}V_B\otimes\bigwedge{\mathbb Z}[\nabla]$ together with its filtration.

In section 5 we investigate the graded pieces $\overline{C}{\mathcal
  F}_{V_B,\epsilon}^{\bullet}$ for this filtration. It turns out that they are direct sums of certain very concrete ${\mathfrak N}$-equivariant complexes ${\mathcal K}_D^{\bullet}(V_B)$, for
subsets $D$ of $\nabla$. To prove freeness of $M_{\chi}(V)$ and exactness of $C_{\alpha}{\mathcal F}_{V_B}^{\bullet}$ it is enough to prove exactness of the complexes ${\mathcal K}_D^{\bullet}(V)$, and in fact exactness of ${\mathcal K}_{\nabla}^{\bullet}(V)$ (i.e. where $D=\nabla$) is enough.  

In section 6 we prove by hand (Corollary \ref{rank2}) that ${\mathcal
  K}_{\nabla}^{\bullet}(V)$ is always exact (and hence that $M_{\chi}(V)$ is free
and $C_{\alpha}{\mathcal F}_{V_B}^{\bullet}$ is exact) in the following cases:\\ 

(a) The
semisimple rank $|\nabla|$ is $1$.

 (b) The
semisimple rank $|\nabla|$ is $2$, we are in case (I) and ${\mathfrak N}$
acts trivially on $V$. \\

The result in case (a) is a theorem of Barthel and Livne \cite{bali}. The freeness of $M_{\chi}(V)$ in case (b) generalizes the main result of
\cite{bo} (which deals with $\widetilde{G}=G={\rm PGL}_3(F)$). 

For general $\widetilde{G}$ we then provide sufficient
criteria (Proposition \ref{bettercond}) for the exactness of ${\mathcal K}_{\nabla}^{\bullet}(V)$ in terms of
its ${\mathfrak N}$-cohomology, in case $R$ is a field of characterictic $p$ or a
discrete valuation ring with residue field of characterictic $p$.  

This motivates our investigations in section 7 (which is logically independent
on the preceding sections). Let ${\mathcal G}$ be a
simple, simply connected split algebraic group over ${\mathbb
  Z}_p$, let ${\mathcal T}$ be a maximal split torus. Let ${\mathcal N}$ be
the unipotent radical of a Borel subgroup in ${\mathcal G}$ containing
${\mathcal T}$, let ${\mathcal N}({\mathbb Z}_p)$ denote its group of ${\mathbb Z}_p$-valued
points. In order to apply Proposition \ref{bettercond} for $F={\mathbb Q}_p$
we compute (Theorem \ref{prodpoloclou}) the cohomology of ${\mathcal
  N}({\mathbb Z}_p)$ with values in rational ${\mathcal G}_{{\mathbb
    F}_p}$-representations $V(\mu)$ with highest weight $\mu$ lying 'in the
(closed) bottom alcove'. More precisely, the condition is that
$\langle\mu+\rho,\check{\beta}\rangle\le p$ for all positive roots $\beta$,
where $\rho$ denotes the half sum of positive roots. The result (Theorem
\ref{poloclou}) is that there is a natural isomorphism\begin{gather}H^i({\mathcal N}({\mathbb
    Z}_p),V(\mu))\cong \prod_{w\in {\mathcal W}\atop
    \ell(w)=i}H^i({\mathcal N}_w({\mathbb Z}_p),V(\mu)_{{\mathcal N}({\mathbb Z}_p)(w)})\end{gather} for each
$i\ge0$. Here ${\mathcal W}$ denotes the (finite) Weyl group of $({\mathcal
  G},{\mathcal T})$ and for $w\in
{\mathcal W}$ of length $\ell(w)=i$ the group ${\mathcal N}_w({\mathbb Z}_p)$ is an $i$-dimensional subgroup of ${\mathcal
  N}({\mathbb Z}_p)$. It acts on $V(\mu)_{{\mathcal N}({\mathbb Z}_p)(w)}$, a certain one-dimensional
quotient of $V(\mu)$. In particular, each $H^i({\mathcal N}_w({\mathbb
  Z}_p),V(\mu)_{{\mathcal N}({\mathbb Z}_p)(w)})$ is one-dimensional. The
proof of Theorem
\ref{poloclou} consists in the reduction to a result of Polo and Tilouine on the Lie algebra
cohomology of ${\mathcal N}$ with values in $V(\mu)$ (who transfered a famous result of Kostant from
characteristic zero to positive characteristic). Besides that, we found it also quite
instructive to compute $H^i({\mathcal N}({\mathbb
    Z}_p),V(\mu))$ by hand, avoiding Lie algebra
cohomology,  in the first non-trivial case: the case where ${\mathcal G}={\rm
    SL}_3$ and $V(\mu)$ is the trivial representation. This computation is
  given in the appendix section \ref{appco}. We conclude section 7 by
  extending the above theorem to non-simple ${\mathcal G}$, Theorem \ref{prodpoloclou}. 

In section 8 we then combine Proposition \ref{bettercond} and Theorem
\ref{prodpoloclou} to prove our main result. Let $J$ be an index set for the
irreducible components of a root system for $\widetilde{G}$, for $j\in J$ let $h^{(j)}$, $\rho^{(j)}$, $\Phi^{(j),+}$ denote the respective Coxeter number, half sum of positve roots, set of positive roots, respectively. In case (II) the restriction of $V$ to the product of universal coverings of the simple factors of the derived group of $\widetilde{\bf G}_{{\mathbb F}_p}$ may be regarded as a tensor product of irreducible rational representations with $p$-restricted highest weights $\mu^{(j)}$. Similarly, in case (III) the rational representation $V\otimes_{{\mathcal O}_F}\overline{\mathbb F}_p$ is a tensor product of irreducible rational representations with $p$-restricted highest weights $\mu^{(j)}$.  \\

{\bf Theorem 1.1} (= Theorem \ref{modest}) {\it Let $F={\mathbb Q}_p$ and assume that we are in (at least) one of the
  following cases (i), (ii) or (iii):

(i) We are in case (I), $R$ is a discrete valuation ring with residue field $k$ of characteristic
$p$, and we have $p\ge h^{(j)}-1$ for all $j\in J$.

(ii) We are in case (II), and $\langle\mu^{(j)}+\rho^{(j)},\check{\beta}\rangle\le
p$ for all $\beta\in\Phi^{(j),+}$, all $j$.

(iii) We are in case (III), and $\langle\mu^{(j)}+\rho^{(j)},\check{\beta}\rangle\le
p$ for all $\beta\in\Phi^{(j),+}$, all $j$.

Then the $B$-module $M_{\chi}(V)$ is
  free and $C_{\alpha}{\mathcal F}_{V_B}^{\bullet}$ is a
  $\widetilde{G}$-equivariant resolution of $M_{\chi}(V)$, for any $R$-algebra
  $B$ and any $R$-algebra
  homomorphism $\chi:{\mathfrak H}_V(\widetilde{G},\widetilde{K})\to B$ such
  that $\chi|_{{\mathfrak H}(V)}$ factors through $R$. }\\

Notice that for $\widetilde{G}={\rm GL}_2({\mathbb Q}_p)$ our hypotheses are empty (the
$p$-restricted region consists of the (closed) bottom alcove only), so in this
case the freeness of the $M_{\chi}(V)$ holds unconditionally (as already
seen --- more generally for $\widetilde{G}={\rm GL}_2(F)$ for any $F$ ---  in
Corollary 6.1, resp. in
\cite{bali}).

Section 9 contains our discussion of existence of $\widetilde{G}$-invariant norms on locally algebraic
$\widetilde{G}$-representations.\\


{\bf (5) Further remarks.}


(5a) We cannot prove that for
general $\widetilde{G}$ our conditions on $V$ imposed in Theorem 1.1 (= Theorem \ref{modest}) are necessary to guarantee the
freeness of the $M_{\chi}(V)$. However, in two of the very first cases not covered
by our hypotheses --- the case $\widetilde{G}={\rm PGL}_3({\mathbb Q}_p)$, taking as $V$ the Steinberg
representation of ${\rm PGL}_3({\mathbb F}_p)$, and the case
$\widetilde{G}={\rm PGL}_4({\mathbb Q}_p)$, $p=2$, taking as $V$ the trivial representation --- we indeed verified that the exactness of ${\mathcal
  K}^{\bullet}_{\nabla}(V)$ fails ! For this reason we are in
fact sceptical on whether the $M_{\chi}(V)$ are free for general $V$.

(5b) The representations $M_{\chi}(V)$ clearly deserve further investigation. Their
freeness immediately implies the existence of supersingular
$\widetilde{G}$-representations over $\overline{\mathbb F}_p$, just as this was
demonstrated in the case of $\widetilde{G}={\rm GL}_2(F)$ in \cite{bali}.

(5c) Theorem \ref{normexis} states the existence of $\widetilde{G}$-invariant
norms on
{\it unramified} principal series representations whenever
certain obvious necessary conditions are fulfilled. One may formulate similar such necessary
conditions for the existence of $\widetilde{G}$-invariant
norms on much more general smooth (or even locally algebraic)
$\widetilde{G}$-representations. These conditions, and that they
should in fact be {\it sufficient} conditions, seem to have
been first observed and propagated by Marie France Vign\'{e}ras; see also
e.g. \cite{bs}, \cite{dat1}, \cite{eme}, \cite{st}, \cite{vig1}. A natural next type of
example class to consider would be that of tamely ramified principal series
representations. In view of the strategy followed here in the unramified case,
one may speculate if the tamely unramified case is accessible to similar techniques involving the (pro-$p$) Iwahori Hecke algebra instead
of the spherical Hecke algebra.\\

{\bf (6) Acknowledgements:} It is a pleasure to thank Florian Herzig and Peter
Schneider for discussions on the subject matter of this paper. I thank the
members of the working group on Geometry, Topology and Group Theory and the
working group on $p$-adic Arithmetic at the University of M\"unster for the
very careful reading of this paper and numerous remarks. In particular Jan
Kohlhaase made several helpful remarks on section \ref{jache}. I am
grateful to Linus Kramer for providing a proof of Lemma
\ref{daniel} being much more elegant than my original approach, and to Peter
Schneider for his critical comments on an earlier version of section \ref{cohin} and in particular for suggesting
to use \cite{gruen} for the proof of Theorem \ref{poloclou}. I thank Niko
Naumann for the hint to consult the paper \cite{dwy}. I am very grateful
to the organizers of the $p$-adic Galois trimester, held in the spring 2010
at the Institute Henri Poincar\'{e} at Paris, for the invitation to this
highly rewarding and instructive activity and for giving me the chance to
report on the results of the present paper there. I thank the referee for suggesting improvements on the exposition.

This work was done while the Deutsche Forschungsgemeinschaft (DFG) supported my position as a Heisenberg Professor at the Humboldt University at Berlin.

\section{The universal module}

\label{defis}

Let $F$ be a local non-Archimedean field with finite residue field $k_F$ of
characteristic $p>0$, let $p_F\in{\mathcal
  O}_F$ denote a fixed uniformizer. Let $\widetilde{G}=\widetilde{\bf
  G}(F)$ be the group of $F$-rational points of a split connected reductive algebraic group $\widetilde{\bf
  G}$ over ${\mathcal O}_F$ with connected center; then $\widetilde{K}=\widetilde{\bf
  G}({\mathcal O}_F)$ is a hyperspecial maximal compact open subgroup of $\widetilde{G}$.

Let $\widetilde{\bf T}$ be a maximal split torus, let
$\widetilde{T}=\widetilde{\bf T}(F)$. Choose a Borel subgroup $N\widetilde{T}$ in $\widetilde{G}$ (containing $\widetilde{T}$ and with unipotent
radical $N$). Let $\overline{N}$ denote the unipotent radical of the Borel
subgroup opposite to $N\widetilde{T}$ and then let ${{{{\mathfrak N}}}}=\widetilde{K}\cap
\overline{N}$. [A more conventional notation for this group would be
  e.g. $\overline{N}_0$. We chose the
  non decorated symbol ${\mathfrak N}$ instead since it is used again and again.] We have the reduction map ${\rm red}:\widetilde{K}=\widetilde{\bf
  G}({\mathcal O}_F)\to \widetilde{\bf
  G}(k_F)$. Let $\widetilde{I}={\rm red}^{-1}({\rm red}({\mathfrak N}))$; this
is a pro-$p$-Iwahori subgroup in $\widetilde{G}$. Put $$\widetilde{T}_+=\{t\in \widetilde{T}\,|\,t(\widetilde{K}\cap
N)t^{-1}\subset (\widetilde{K}\cap N)\}=\{t\in \widetilde{T}\,|\,t^{-1}{{{{\mathfrak N}}}}t\subset{{{{\mathfrak N}}}}\},$$a submonoid
of $\widetilde{T}$. Let $X^*(\widetilde{\bf T})={\rm Hom}(\widetilde{\bf T},{\mathbb
  G}_m)$. Let $X_*({\widetilde{\bf T}})={\rm Hom}({\mathbb
  G}_m,{\widetilde{\bf T}})$. We have the natural isomorphism
\begin{gather}X_*({\widetilde{\bf T}})\cong {\widetilde{\bf T}(F)}/{\widetilde{\bf
    T}({\mathcal O}_F)}=\widetilde{T}/(\widetilde{T}\cap\widetilde{K}),\quad
\lambda\mapsto\mbox{ class of }\lambda(p_F).\label{unifiso}\end{gather} Let
${\Phi}^+$ (resp. ${\Phi}^-$) denote the set of positive (resp. negative) roots
with respect to $N\widetilde{T}$, let $\Phi={\Phi}^+\cup{\Phi}^-$. Let $X_*({\widetilde{\bf T}})_{+}$ (the monoid
of dominant coweights) denote the image of
$\widetilde{T}_+$ in the quotient $X_*({\widetilde{\bf T}})$ of
$\widetilde{T}$ (under the isomorphim (\ref{unifiso})). Equivalently,
$\lambda\in X_*({\widetilde{\bf T}})$ lies in $X_*({\widetilde{\bf T}})_{+}$ if
and only if for a representative $t\in \widetilde{T}$ of $\lambda$ we have $|\alpha(t)|\le1$ for all $\alpha\in\Phi^+$.

The quotient ${G}=\widetilde{G}/Z$ of
$\widetilde{G}$ by its center $Z$ is a split semisimple algebraic group over $F$ of adjoint
type. Of course, more precisely, $G$ is the group ${\bf
  G}(F)$ of $F$-valued points of a split semisimple algebraic group ${\bf
  G}$ over $F$ of adjoint
type. Moreover, ${K}=\widetilde{K}/(Z\cap\widetilde{K})$ is a special maximal compact open subgroup of
${G}$, and ${T}=\widetilde{T}/Z$ is a split maximal torus in
${G}$, i.e. ${T}={\bf T}(F)$ for a split
maximal torus ${\bf T}$ in ${{\bf
  G}}$. We let
$W=N(T)/T=N(\widetilde{T})/\widetilde{T}$ denote the Weyl group. The sets of roots for $\widetilde{G}$ and for ${G}$ can be naturally
identified. For $\alpha\in\Phi$ let $U_{\alpha}$ be the
corresponding root subgroup, regarded as a subgroup of both $G$ and
$\widetilde{G}$. Similarly we may identify $N$ and ${{{{\mathfrak N}}}}$ with
their images in $G$.
The image $I$ of $\widetilde{I}$ in $G$ is a pro-$p$-Iwahori subgroup in $G$ containing ${{{{\mathfrak N}}}}$, and
$NT$ is a Borel subgroup in $G$. What has been explained for
$\widetilde{G}$ above applies similarly to $G$. In particular, the image
${T}_+$ of $\widetilde{T}_+$ in $T$ can be characterized
as ${T}_+=\{t\in{T}\,|\,t^{-1}{{{{\mathfrak N}}}}t\subset {{{{\mathfrak N}}}}\}$, and its image in
$X_*({\bf T})$ is $X_*({{\bf T}})_{+}$. As $Z$ is connected we have the following
facts:

---  $Z/(Z\cap \widetilde{K})$ is a free abelian group. 

--- There is a unique subset $\nabla$ of $X_*({{\bf T}})_{+}$ which generates
  $X_*({{\bf T}})_{+}$ as a monoid (so $X_*({{\bf T}})_+$ is a direct sum of copies of ${\mathbb N}$) and is a ${\mathbb Z}$-basis of
  $X_*({{\bf T}})$.

We have the exact sequence of monoids$$0\longrightarrow
Z/(Z\cap\widetilde{K})\longrightarrow X_*(\widetilde{\bf T})_+\longrightarrow
X_*({{\bf T}})_+\longrightarrow0$$and we choose once
and for all a splitting\begin{gather}X_*({\widetilde{\bf T}})_+\cong Z/(Z\cap\widetilde{K})\bigoplus
X_*({{\bf T}})_+\label{monoiex}\end{gather}of this exact sequence and allow us to regard $X_*({{\bf T}})_+$ as being contained in $X_*({\widetilde{\bf T}})_+$. \\

Let $R$ be a commutative ring with unity and let $V$ be an
$R[\widetilde{K}]$-module, finitely generated as an $R$-module. Let $${\rm
  ind}_{\widetilde{K}}^{\widetilde{G}}V=\{f:\widetilde{G}\longrightarrow
  V\,|\,{\rm supp}(f)\mbox{ compact}, f(gk)=k^{-1}f(g)\mbox{ for }g\in
  \widetilde{G}, k\in\widetilde{K}\}$$denote the compactly induced
  $\widetilde{G}$-representation with $\widetilde{G}$-action
  $(gf)(g')=f(g^{-1}g')$. Let ${\mathfrak
    H}_V(\widetilde{G},\widetilde{K})={\rm End}_{\widetilde{G}}({\rm
  ind}_{\widetilde{K}}^{\widetilde{G}}V)$ denote the corresponding Hecke
algebra. By Frobenius reciprocity, we have a natural isomorphism between ${\mathfrak
    H}_V(\widetilde{G},\widetilde{K})$ and \begin{gather}\{\Psi:\widetilde{G}\to
  {\rm End}_R(V)\,|\,{\rm supp}(\Psi)\mbox{ compact}, \Psi(k_1gk_2)=k_1\Psi(g)k_2\mbox{ for }g\in
  \widetilde{G}, k_1,k_2\in\widetilde{K}\}\label{frohec}\end{gather}by sending
$\Psi$
to $$[f\mapsto[g\mapsto\sum_{x\widetilde{K}\in\widetilde{G}/\widetilde{K}}\Psi(g^{-1}x)(f(x))]]\in{\mathfrak
    H}_V(\widetilde{G},\widetilde{K}).$$Multiplication in ${\mathfrak
    H}_V(\widetilde{G},\widetilde{K})$ corresponds to convolution in (\ref{frohec}). We define the
subalgebra ${\mathfrak H}(V)$ of ${\mathfrak
    H}_V(\widetilde{G},\widetilde{K})$ as the one corresponding under this isomorphism
  to the subalgebra of all $\Psi$ in
  (\ref{frohec}) with support in $Z\widetilde{K}=\widetilde{K}Z$.\\

We will always assume that we are in (at least) one of the following three cases:\\

(I) $R$ is a field or a principal ideal domain and $V$ is $R$-free of rank one.

(II) $R=\overline{{\mathbb F}}_p$ and $V$ is an irreducible smooth $\widetilde{K}$-representation
over $R$. Thus ${\rm dim}_R(V)<\infty$ and the
$\widetilde{K}$-action on $V$ factors through the natural projection
$\widetilde{K}\to\widetilde{\bf G}(k_F)$.

(III) $R={\mathcal O}_F$ where now $F$ is a finite extension of ${\mathbb Q}_p$ and $V$ is a $\widetilde{\bf G}$-stable
${\mathcal O}_F$-lattice in an irreducible rational representation of $\widetilde{\bf
  G}_{{F}}$, and moreover, $V\otimes_{{\mathcal
  O}_F}{\overline{\mathbb F}_p}$ is an irreducible rational representation of
$\widetilde{\bf G}_{k_{F}}$.\\

In case (I) we define
  $\xi_{\lambda}:V\to V$ for $\lambda\in X_*(\widetilde{\bf T})_+$ to be the identity on $V$.
 
In case (II) we follow
\cite{her} to define $\xi_{\lambda}:V\to V$ for $\lambda\in
X_*({\widetilde{\bf T}})$. Let
$\widetilde{\bf P}_{\lambda}=\widetilde{\bf M}_{\lambda}{\bf N}_{\lambda}\subset \widetilde{\bf G}$ denote the parabolic subgroup
defined by $\lambda$ (cf. \cite{spri} 13.4.2, 15.4.4), with Levi factor $\widetilde{\bf M}_{\lambda}$ containing $\widetilde{\bf T}$ and unipotent radical ${\bf N}_{\lambda}$. Then $\widetilde{\bf P}_{\lambda}(k_F)$ is the image of
$\widetilde{K}\cap \lambda(p_F)^{-1}\widetilde{K}\lambda(p_F)$ under the
natural projection $\widetilde{K}\to\widetilde{\bf G}(k_F)$, cf. \cite{her}
Proposition 3.8. The composition of
natural maps $$V^{{\bf N}_{-\lambda}(k_F)}\longrightarrow V\longrightarrow V_{{\bf N}_{\lambda}(k_F)}$$ is an
isomorphism of $\widetilde{\bf M}_{\lambda}(k_F)$-modules. We use its inverse to
define the map $\xi_{\lambda}:V\to V$ as the composition $$\xi_{\lambda}:V\to
V_{{\bf N}_{\lambda}(k_F)}\cong V^{{\bf N}_{-\lambda}(k_F)}\to V.$$ 

Finally, in case (III) we define $\xi_{\lambda}:V\to V$ as in the proof of
\cite{her} Proposition 2.10, as follows. The irreducible rational representation $V\otimes_{{\mathcal
  O}_F}{\overline{\mathbb F}_p}$ of
$\widetilde{\bf G}_{k_{F}}$ remains irreducible as a representation of the
abstract group
$\widetilde{\bf G}(k_{F})$ (cf. \cite{her}). Therefore the construction in
case (II) yields maps $\xi_{\lambda}:V\otimes_{{\mathcal
  O}_F}{\overline{\mathbb F}_p}\to V\otimes_{{\mathcal
  O}_F}{\overline{\mathbb F}_p}$ for $\lambda\in
X_*({\widetilde{\bf T}})$. In the proof of
\cite{her} Proposition 2.10 it is explained that they lift naturally to
maps $\xi_{\lambda}:V\to V$.\\   

In either case, (I), (II) or (III), for $\lambda\in X_*({\widetilde{\bf T}})$ we denote by
$T_{\lambda}\in {\mathfrak
    H}_V(\widetilde{G},\widetilde{K})$ the element corresponding to the unique map $\Psi_{\lambda}:\widetilde{G}\longrightarrow
  {\rm End}(V)$ in (\ref{frohec}) with support in
  $\widetilde{K}\lambda(p_F)\widetilde{K}$ and sending $\lambda(p_F)$ to
  $\xi_{\lambda}$. (For the existence of such $T_{\lambda}$ see again \cite{her}.)

\begin{pro}\label{redsat} (Schneider/Teitelbaum, Herzig, Henniart/Vign\'{e}ras) ${\mathfrak
    H}_V(\widetilde{G},\widetilde{K})$ is the free commutative polynomial algebra over ${\mathfrak H}(V)$ in the variables $T_{\lambda}$ for $\lambda\in\nabla$.
\end{pro}

This is all we need to know about ${\mathfrak
    H}_V(\widetilde{G},\widetilde{K})$ in the present paper. It is a
  consequence of the following Proposition \ref{prredsat}. 

\begin{pro}\label{prredsat} (Schneider/Teitelbaum, Herzig, Henniart/Vign\'{e}ras) (a) We have a canonical injective $R$-algebra homomorphism\begin{gather}{\mathfrak S}:{\mathfrak
    H}_V(\widetilde{G},\widetilde{K})\longrightarrow R[X_*({\widetilde{\bf
      T}})].\label{satboth}\end{gather}In particular, ${\mathfrak
    H}_V(\widetilde{G},\widetilde{K})$ is commutative.

(b) Under the homomorphism (\ref{satboth}), the subalgebra ${\mathfrak H}(V)$ of ${\mathfrak H}_V(\widetilde{G},\widetilde{K})$ maps isomorphically to the subalgebra $R[Z/(Z\cap\widetilde{K})]$ of $R[X_*(\widetilde{\bf T})]$.

(c) In case (I), the image of ${\mathfrak S}$ is the subalgebra $R[X_*({\widetilde{\bf
      T}})]^W$ of $W$-invariants in $R[X_*({\widetilde{\bf T}})]$ for a
  certain twisted action of $W$ on $R[X_*({\widetilde{\bf T}})]$. 

 In case (II) the image of ${\mathfrak S}$ is the subalgebra of
functions supported on $\widetilde{T}_+$ --- under the natural identification
of $R[X_*({\widetilde{\bf T}})]$ with the $R$-algebra of
$\widetilde{T}\cap\widetilde{K}$-biequivariant, compactly supported $R$-valued
functions on $\widetilde{T}$.  
\end{pro}

{\sc Proof:} In case (II) this is shown in \cite{her} if $F$ has
characteristic zero; for $F$ of general characteristic (and in fact for more
general, not necessarily split
$\widetilde{G}$) see \cite{hv}. In case (III) this is derived from the results
in \cite{st} and \cite{her}, notably from \cite{her} Proposition 2.10 and its proof. In case (I) this is shown in \cite{st},
section 3, from which we recall the basic features. We begin by explaining the
twisted action of $W$ on $R[X_*({\widetilde{\bf T}})]$. The $W$-action on ${\widetilde{\bf T}}$ induces a $W$-action on $X_*({\widetilde{\bf
      T}})$. Let $\delta:\overline{N}{\widetilde{T}}\to {\mathbb
    Q}^{\times}\cap{\mathbb Z}[1/|k_F|]$ denote the modulus character of the
  Borel subgroup $\overline{N}{\widetilde{T}}$ of ${\widetilde{G}}$. We
  temporarily choose a square root of $|k_F|$ to
  define $$\gamma(w,\lambda)=\frac{\delta^{1/2}(w(\lambda)(p_F))}{\delta^{1/2}(\lambda(p_F))}$$for
  $w\in W$ and $\lambda\in X_*({\widetilde{\bf T}})$. In fact this is an
  element in ${\mathbb Z}$ (cf. \cite{st} section 2), hence maps to $R$. The map $\gamma(.,.):W\times X_*({\widetilde{\bf T}})\to{\mathbb Q}^{\times}$ is a cocylce (i.e. $\gamma(w,\lambda\mu)=\gamma(w,\lambda)\gamma(w,\mu)$ and $\gamma(vw,\lambda)=\gamma(v,w(\lambda))\gamma(w,\lambda)$ for any $v,w\in W$, any $\lambda,\mu\in X_*({\widetilde{\bf T}})$). Therefore we can define an action of $W$ on $R[X_*({\widetilde{\bf T}})]$ through the formula$$W\times
  R[X_*({\widetilde{\bf T}})]\longrightarrow R[X_*({\widetilde{\bf
      T}})],$$$$(w,\sum_{\lambda}c_{\lambda}\lambda)\mapsto
  \sum_{\lambda}\gamma(w,\lambda)c_{\lambda}w(\lambda).$$The Cartan decomposition says that $\widetilde{G}$ is the disjoint union
of the double cosets $\widetilde{K}\lambda(p_F)\widetilde{K}$, with $\lambda$ running through
$X_*({\widetilde{\bf T}})_{+}$, hence $\{T_{\lambda}\}_{\lambda\in
  X_*({\widetilde{\bf T}})_{+}}$ is an $R$-basis of ${\mathfrak
    H}_V(\widetilde{G},\widetilde{K})$. Consider the map$${\mathfrak S}:{\mathfrak
    H}_V({\widetilde{G}},{\widetilde{K}})\longrightarrow R[X_*({\widetilde{\bf
      T}})],$$$$\psi\mapsto\sum_{\lambda\in X_*({\widetilde{\bf
      T}})}(\sum_{n\in
  \overline{N}/{{{{{\mathfrak N}}}}}}\psi({\lambda}(p_F)n))\lambda.$$(Here $\psi\in {\mathfrak
    H}_V({\widetilde{G}},{\widetilde{K}})$ is identified with the element in (\ref{frohec}) it corresponds to.) The map ${\mathfrak S}$ differs
from the classical Satake map by the character $\delta^{1/2}$ of $X_*({\widetilde{\bf T}})$. In particular, ${\mathfrak S}$ is an $R$-algebra
homomorphism and takes values in $R[X_*({\widetilde{\bf
      T}})]^W$. We define $\sigma_{\mu}\in R[X_*({\widetilde{\bf
      T}})]$ for $\mu\in X_*({\widetilde{\bf T}})_+$ as follows: if $W(\mu)$
  denotes the stabilizer of $\mu$ in $W$, then $$\sigma_{\mu}=\sum_{w\in
    W/W(\mu)}\gamma(w,\mu)w(\mu).$$It is easy to see that the $\sigma_{\mu}$
  for all $\mu\in X_*({\widetilde{\bf T}})_+$ form an $R$-basis of $R[X_*({\widetilde{\bf
      T}})]^W$. If for $\lambda\in X_*({\widetilde{\bf T}})_{+}$ we
write\begin{gather}{\mathfrak S}(T_{\lambda})=\sum_{\mu\in X_*({\widetilde{\bf
      T}})_+}c(\mu,\lambda)\sigma_{\mu}\quad\mbox{ with }c(\mu,\lambda)\in R,\label{satcoef}\end{gather}then it
follows from \cite{bt} Proposition I.4.4.4 that for all $\lambda, \mu\in
X_*({\widetilde{\bf T}})_{+}$ we have: $c(\lambda,\lambda)=1$, and $c(\mu,\lambda)=0$ unless $\mu\le \lambda$. Here for $\mu,\lambda\in X_*({\widetilde{\bf T}})_{+}$ we write $\mu\le
  \lambda$ if and only if $\lambda-\mu$ is a non-negative real linear
  combination of the positive coroots (i.e. the coroots associated to the
  elements of $\Phi^+$). Hence ${\mathfrak S}$ is represented as an $R$-linear map to
  $R[X_*({\widetilde{\bf T}})]^W$ by an infinite block matrix with unipotent
  triangular finite blocks, hence is bijective.\hfill$\Box$\\

{\sc Proof of Proposition \ref{redsat}}: Recall that we regard $X_*({{\bf T}})_{+}$ as a subset of $X_*({\widetilde{\bf T}})_{+}$ via the splitting (\ref{monoiex}), hence we may form $T_{\lambda}\in {\mathfrak H}_V(\widetilde{G},\widetilde{K})$ for $\lambda\in\nabla\subset X_*({{\bf T}})_{+}$. To prove Proposition \ref{redsat} it is enough, by the Cartan decomposition, to
show that the $T_{\lambda}$ for any
$\lambda\in X_*({\widetilde{\bf T}})_+$ are contained in the ${\mathfrak H}(V)$-algebra
generated by the $T_{\lambda}$ for all $\lambda\in \nabla$. Now for any two
$\lambda, \lambda'\in X_*({\widetilde{\bf T}})_+$ we
claim that\begin{gather}T_{\lambda}T_{\lambda'}-T_{\lambda+\lambda'}\in
  \sum_{\lambda''\in X_*({\widetilde{\bf T}})_+\atop
    \lambda''<\lambda+\lambda'}{\mathfrak H}(V).T_{\lambda''}.\label{convolcoe}\end{gather}Indeed, first we remark that for any two
$\mu, \mu'\in X_*({\widetilde{\bf T}})_+$ we have$$\sigma_{\mu+\mu'}-\sigma_{\mu}\sigma_{\mu'}\in\sum_{\mu''<\mu+\mu'}{\mathfrak H}(V).\sigma_{\mu''}.$$This can be shown as in \cite{bou} chapter VI, sections 3.2 and 3.4; see also the proof of \cite{st} Proposition 2.6. Therefore, since ${\mathfrak S}$ is an {\it algebra} isomorphism, formula (\ref{convolcoe}) can be translated via ${\mathfrak S}$ into a similar statement in $R[X_*({\widetilde{\bf T}})]$, and this statement holds true by the informations on the coefficients $c(\mu,\lambda)$ in formula
  (\ref{satcoef}) recorded above. By \cite{her} these arguments apply in case (II) as well. Alternatively, one can prove formula (\ref{convolcoe}) more directly by using \cite{bt} Prop. 4.4.4 (iii) and the formula for multiplication in ${\mathfrak
    H}_V(\widetilde{G},\widetilde{K})$ (i.e. convolution). As the monoid
  $X_*({\widetilde{\bf T}})_+$ is generated modulo $Z/(Z\cap\widetilde{K})$ by the elements in $\nabla$, formula
  (\ref{convolcoe}) shows that we only need to observe that the non-zero minimal (with
  respect to our ordering on $X_*({\widetilde{\bf T}})_+$ introduced above)
  elements in $X_*({\widetilde{\bf T}})_+$ belong to $\nabla$. But
  this is clear.\hfill$\Box$\\

For $z\in Z$ let $T_{z}$ denote the element in ${\mathfrak H}(V)$ corresponding to the unique map $\widetilde{G}\longrightarrow
  {\rm End}_R(V)$ in (\ref{frohec}) with support in
  $\widetilde{K}z$ and sending $z$ to ${\rm id}_V$.

Suppose we are given an $R$-algebra homomorphism $\theta:{\mathfrak H}(V)\to R$. We let $Z$ act on $V$ by
  declaring that $zv=\theta(T_{z^{-1}})v$ for $z\in Z$ and $v\in V$; this is
  compatible with the $\widetilde{K}$-action, i.e. defines an action of the subgroup $\widetilde{K}Z$ of $\widetilde{G}$ on $V$ (cf. e.g. the discussion preceding
\cite{herzclas} formula
(3.5)). Consider the $\widetilde{G}$-representation$${\rm
  ind}_{\widetilde{K}Z}^{\widetilde{G}}V=\{f:\widetilde{G}\longrightarrow
  V\,|\,{\rm supp}(f)\mbox{ compact mod } Z, f(gk)=k^{-1}f(g)\mbox{ for }g\in
  \widetilde{G}, k\in\widetilde{K}Z\}.$$It depends on $\theta$ but we suppress
  this in our notation. The natural map (provided by
Frobenius reciprocity) ${\rm
  ind}_{\widetilde{K}}^{\widetilde{G}}V\to{\rm
  ind}_{\widetilde{K}Z}^{\widetilde{G}}V$ induces an isomorphism\begin{gather}{\rm
  ind}_{\widetilde{K}}^{\widetilde{G}}V\otimes_{{\mathfrak H}(V),\theta}R\cong{\rm
  ind}_{\widetilde{K}Z}^{\widetilde{G}}V,\label{isocind}\end{gather}see \cite{herzclas} formula
(3.5). Using the composition \begin{gather}\theta^*:{\mathfrak
    H}_V(\widetilde{G},\widetilde{K})\cong {\mathfrak
    H}(V)[T_{\lambda}]_{\lambda\in\nabla}={\mathfrak
    H}(V)\otimes_RR[T_{\lambda}]_{\lambda\in\nabla}\stackrel{\theta\otimes{\rm
      id}}{\longrightarrow}R[T_{\lambda}]_{\lambda\in\nabla}\label{thetastern}\end{gather}where
the first isomorphism is that of Proposition \ref{redsat}, we may also read the isomorphism (\ref{isocind}) as$${\rm ind}_{\widetilde{K}}^{\widetilde{G}}V\otimes_{{\mathfrak
    H}_V(\widetilde{G},\widetilde{K})}R[T_{\lambda}]_{\lambda\in\nabla}\cong{\rm
  ind}_{\widetilde{K}Z}^{\widetilde{G}}V.$$In this way we see the $T_{\lambda}$ acting on ${\rm
  ind}_{\widetilde{K}Z}^{\widetilde{G}}V$. These actions may simply be described as follows. Let $\overline{\Psi}_{\lambda}:\widetilde{G}\to{\rm End}_R(V)$ denote the unique function
with support in $\widetilde{K}\lambda(p_F)\widetilde{K}Z$, with
$\overline{\Psi}_{\lambda}(k_1gk_2)=k_1\overline{\Psi}_{\lambda}(g)k_2$ for
$k_1, k_2\in \widetilde{K}Z$ and $g\in \widetilde{G}$, and with
$\overline{\Psi}_{\lambda}(\lambda(p_F))=\xi_{\lambda}$. Then $T_{\lambda}$ acts on ${\rm
  ind}_{\widetilde{K}Z}^{\widetilde{G}}V$ as \begin{gather}f\mapsto T_{\lambda}(f)=[g\mapsto\sum_{x
  \widetilde{K}Z\in
  \widetilde{G}/\widetilde{K}Z}\overline{\Psi}_{\lambda}(g^{-1}x)f(x)]\label{heckebuild}\end{gather}where
$f\in{\rm
  ind}_{\widetilde{K}Z}^{\widetilde{G}}V$ and $g\in \widetilde{G}$. Writing
  $\widetilde{K}\lambda(p_F)\widetilde{K}Z=\coprod_ih_i\widetilde{K}Z$ for suitable $h_i\in \widetilde{G}$ this becomes\begin{gather}T_{\lambda}(f)(g)=\sum_i\overline{\Psi}_{\lambda}(h_i)f(gh_i).\label{hecfor}\end{gather}

We conclude this section by recalling how ${\rm
  ind}_{\widetilde{K}Z}^{\widetilde{G}}V$ can be spread out onto the
Bruhat-Tits building $X$ of ${\bf G}/F$. Let $A=X_*({\bf T})\otimes{\mathbb R}$ denote the
standard apartment with respect to $T$ in $X$. Let $X^0$, resp. $A^0$,
denote the set of special vertices of $X$, resp. of $A$. Then $G$,
resp. $X_*({{\bf T}})\cong {{\bf T}(F)}/{{\bf T}({\mathcal O}_F)}=T/T\cap K$, act on
$X^0$, resp. on $A^0$.

\begin{lem}\label{trans} The action of $X_*({{\bf T}})$ on $A^0$ is transitive. The action of $G$ on $X^0$ is transitive.

\end{lem}

{\sc Proof:} Since ${G}$ is of adjoint type we have ${\mathbb Z}\Phi=X^*({\bf
  T})$, hence ${\mathbb Z}\Phi$ is in perfect duality with $X_*({{\bf
    T}})$. For any element in ${\mathbb Z}\Phi$ we have an associated wall in
$A$, and $A^0$ may be characterized as
the set of points in $A$ contained in just as many different such walls as is the unique special vertex in $X^0$ fixed by $K$. Together the first statement follows. The second
statement then follows by $G$-equivariance. Alternatively, see \cite{tits} section 2.5.\hfill$\Box$\\

 Let $x_0\in X^0$ be the unique special vertex in $X^0$ fixed by $K$. Let $\pi$ denote the map $\widetilde{G}\to X^0,\, g\mapsto gx_0$. By Lemma
\ref{trans} it induces a bijection
$$\widetilde{G}/\widetilde{K}Z=G/K\cong X^0.$$Let $V$ be an
$R[\widetilde{K}]$-module as before, and let $\theta:{\mathfrak H}(V)\to R$ be
an $R$-algebra homomorphism, so that we may form the $\widetilde{G}$-representation ${\rm
  ind}_{\widetilde{K}Z}^{\widetilde{G}}V$. For a subset $Y$ of $X^0$ we put $${\mathcal
    C}(Y,V)=\{f\in{\rm
  ind}_{\widetilde{K}Z}^{\widetilde{G}}V\,|\,{\rm supp}(f)\subset
\pi^{-1}(Y)\}$$(Once more: this also depends on $\theta$). In particular we have \begin{gather}{\mathcal
    C}(X^0,V)={\rm
  ind}_{\widetilde{K}Z}^{\widetilde{G}}V,\label{geodarst}\end{gather} and the $T_{\lambda}$ for
$\lambda\in X_*(\widetilde{\bf
    T})_+$ (or $\lambda\in X_*({\bf
    T})_+$ --- recall the splitting (\ref{monoiex})) act on ${\mathcal
    C}(X^0,V)$ as described in formula (\ref{heckebuild}). We regard the ${\mathcal
    C}(Y,V)$ for any $Y\subset X^0$ as $R$-submodules in ${\mathcal
    C}(X^0,V)$ (extension of functions by zero); in particular we may apply the
  $T_{\lambda}$ to them (to obtain new elements in ${\mathcal
    C}(X^0,V)$). On the other hand, given any $f\in {\mathcal
    C}(X^0,V)$ we denote by $f|_{Y}$ the element in ${\mathcal
    C}(Y,V)$ obtained by restriction to $\pi^{-1}(Y)$.

\section{Vertex coefficient systems}

\label{sectcox}

Let $W$ denote a finite reflection group acting on a Euclidean vector space $V$. Fix an {\it open} Weyl chamber ${\mathcal C}$ in $V$. The intersection of the closure $\overline{\mathcal C}$ of ${\mathcal C}$ with the coweight lattice in $V$ is a free commutative monoid with a unique basis $\nabla$. In particular, $\overline{\mathcal C}=\{\sum_{\lambda\in\nabla}r_{\lambda}\lambda\,|\,(r_{\lambda})_{\lambda}\in {\mathbb R}_{\ge0}^{\nabla}\}$.

\begin{lem}\label{daniel} Let $z\in-{\mathcal C}$ and let $Q$ be a subset of $\nabla$. The assignment$$W^Q\longrightarrow{\mathbb R}_{\ge0},\quad\quad
  (w_{\lambda})_{\lambda\in Q}\mapsto \|z-\sum_{\lambda\in
    Q}w_{\lambda}\lambda\|$$attains its maximum {\it exactly} for those
  $(w_{\lambda})_{\lambda\in Q}\in W^Q$ with $\sum_{\lambda\in Q}w_{\lambda}\lambda=\sum_{\lambda\in Q}\lambda$. 
\end{lem}

{\sc Proof:} (Linus Kramer)

{\it Step 1:} Fix $x\in \overline{\mathcal C}$. We claim that the following (i), (ii) and (iii) are equivalent conditions on $w\in W$: 

(i) $wx=x$.

(ii) The map $W\to{\mathbb R}_{\ge0}, w\mapsto\langle z,wx\rangle$ attains a minimum in $w$.

(iii) The map $W\to{\mathbb R}_{\ge0}, w\mapsto \|z-wx\|^2$ attains a maximum in $w$.

Indeed, we have $\|z-w(x)\|^2=\|z\|^2+\|x\|^2-2\langle z,wx\rangle$, so we see that (ii) and (iii) are equivalent. It remains to prove that (iii) implies (i) (which then also gives the reverse implication). Suppose that $wx\ne x$. Since $\overline{\mathcal C}$ is a fundamental domain for the action of $W$ on $V$ we find a reflection hyperplane $H$ which separates $x$ and $wx$. Then $z$ and $wx$ ly on the same connected component of the complement $V-H$ of $H$ in $V$. This implies that if $h\in W$ denotes the reflection in $H$, then $\|z-hwx\|^2>\|z-wx\|^2$, so $\|z-wx\|^2$ is not maximal.

{\it Step 2:} To prove the Lemma, consider$$\|z-\sum_{\lambda\in Q}w_{\lambda}\lambda\|^2=\|z\|^2+\sum_{\lambda\in Q}\|\lambda\|^2-2\sum_{\lambda\in Q}\langle z,w_{\lambda}\lambda\rangle+2\sum_{\lambda\ne\mu\in Q}\langle w_{\lambda}\lambda, w_{\mu}\mu\rangle.$$By step 1, the terms $\langle z,w_{\lambda}\lambda\rangle$ are minimized if $w_{\lambda}\lambda\in\overline{\mathcal C}$. On the other hand, the condition that $w_{\lambda}\lambda\in\overline{\mathcal C}$ for all $\lambda\in Q$ clearly maximizes the terms $\langle w_{\lambda}\lambda, w_{\mu}\mu\rangle$. Together it follows that $\|z-\sum_{\lambda\in Q}w_{\lambda}\lambda\|^2$ and hence $\|z-\sum_{\lambda\in Q}w_{\lambda}\lambda\|$ is maximized if $w_{\lambda}\lambda\in\overline{\mathcal C}$ for all $\lambda\in Q$, or equivalently if $w_{\lambda}\lambda=\lambda$ for all $\lambda\in Q$, or equivalently if $\sum_{\lambda\in Q}w_{\lambda}\lambda=\sum_{\lambda\in Q}\lambda$.   \hfill$\Box$\\

We return to the setting and notations of section \ref{defis}. The group $W$ acts on $A$ as a reflection group, the group $X_*({\bf T})$ acts
on $A$ by translations. These group actions induce group actions on the set
$A^0$. We have $x_0\in A^0$, and inside $A^0$ this vertex is the unique point fixed by $W$. It is also the origin for the structure of $A$ as an ${\mathbb R}$-vector space
obtained from the group structure of $X_*({\bf T})$.

A {\it wall} (cf. already the proof of Lemma \ref{trans}) in the affine
Coxeter complex $A$ is a translate by an element in $X_*({\bf T})$ of a wall passing through $x_0$. A wall
through $x_0$ can be described as follows (as it is best adapted to our
purposes): it is a $1$-codimensional sub ${\mathbb R}$-vector space (in the ${\mathbb R}$-vector space $A$) spanned by a set $w(\nabla-\{\lambda\})$, for some $w\in
W$, some $\lambda\in\nabla$.

In particular, our walls must be distinguished from the $p$-walls in
\cite{jant}. On the other hand, in later sections we will use the terminology
of alcoves in $X^*({\bf T})$: those are then meant to be $p$-alcoves as in
\cite{jant}.

Using the structure of $A$ as an ${\mathbb R}$-vector space we define the
closed cone (vector chamber)$$\overline{C}=\{\sum_{\lambda\in\nabla}(r_{\lambda}\lambda) x_0\in A=X_*({{\bf T}})\otimes{\mathbb R}  \,;\,(r_{\lambda})_{\lambda}\in{\mathbb
  R}_{\ge0}^{\nabla}\}.$$Let $\sigma_0$ denote the unique chamber (connected
component of the complement
of all walls in $A$) contained in $\overline{C}$ and whose closure contains $x_0$. Choose a point $z_0\in \sigma_0$ such that no two different elements of $A^0$ are at the same Euclidean distance from $z_0$. In other words, for any $\epsilon\in{\mathbb R}$ the set $\{z\in A^0\,;\,\|z-z_0\|=\epsilon\}$ is either empty or contains a single element $z_{\epsilon}$.

[To see that such a $z_0$ exists observe that for a pair of two different elements $z\ne z'$ of $A^0$ the locus of equidistance from $z$ and $z'$ is an affine hyperplane in $A$. But there are only countably many pairs $z\ne z'$ in $A^0$, therefore the open subset $\sigma_0$ of $A$ can not be contained in the union of all these affine hyperplanes, i.e. there is a point $z_0\in \sigma_0$ avoiding all these affine hyperplanes.]

Let ${{\mathcal E}}$ denote the set of all $\epsilon\in{\mathbb R}$ such that
$z_{\epsilon}$ exists. Thus ${{\mathcal E}}$ is in bijection with $A^0$.

Our choice of $z_0$ implies in particular that $z_0$ is not contained in any
wall of $A$. Therefore there is for each $z\in A^0$ a
unique $w_z\in W$ with $$z\in z_0-w_z\overline{C}=\{z_0-w_z(c)\,|\,c\in\overline{C}\}.$$For $\epsilon\in{{\mathcal E}}$
let $$Q(\epsilon)=\{\lambda\in\nabla\,;\,(w_{z_{\epsilon}}\lambda)z_{\epsilon}\in
z_0-w_{z_{\epsilon}}\overline{C}\}.$$

Put \begin{align}{\mathfrak P}&=A^0\times\{Q\,;\,Q\subset\nabla\}\notag\\{}&=\{(z,Q)\,;\,z\in
 A^0\mbox{ and }Q\subset \nabla\}.\notag\end{align}For $(z,Q)\in{\mathfrak P}$
put$$\epsilon(z,Q)={\rm
    max}\{\|x-z_0\|\,;\,x=\mu z\mbox{ for some
}\mu=-\sum_{\lambda\in Q}w_{\lambda}\lambda\mbox{ with }w_{\lambda}\in
W\}.$$For $\epsilon\in{{\mathcal E}}$ define the subset$${\mathfrak P}(\epsilon)=\{(\sum_{\lambda\in Q}(w_{z_{\epsilon}}\lambda)z_{\epsilon},Q)\,;\,Q\subset Q(\epsilon)\}$$of ${\mathfrak P}$. Notice that for two pairs $(z,Q), (z',Q')\in {\mathfrak P}(\epsilon)$ with $z=z'$ we necessarily have $Q=Q'$, but we find it convenient to build in this redundancy in the definition of ${\mathfrak P}(\epsilon)$.

\begin{lem}\label{elemgeo} (a) For any $\epsilon\in{{\mathcal E}}$ and any $(z,Q)\in {\mathfrak P}(\epsilon)$ we have $\epsilon(z,Q)=\epsilon$.

(b) ${\mathfrak P}=\coprod_{\epsilon\in{{\mathcal E}}}{\mathfrak P}(\epsilon)$
(disjoint union). 

\end{lem}

{\sc Proof:} Let $(z,Q)\in {\mathfrak P}$. In the vector space
$A$ with origin $x_0$ we may regard the $x=\mu z$ for $\mu=-\sum_{\lambda\in Q}w_{\lambda}\lambda$ (some $w_{\lambda}\in
W$) simply as the sums of the vectors $\mu$ and $z$ (since $\mu\in X_*({{\bf
    T}})$ acts by translation on $A$). We claim that $\|x-z_0\|$ is maximized if
and only if all the $w_{\lambda}\lambda$ belong to $w_z\overline{C}$, i.e. if
and only if $w_{\lambda}=w_z$ for all $\lambda\in Q$, i.e. if and only if
$x=x_{(z,Q)}$ where we put $x_{(z,Q)}=(-\sum_{\lambda\in Q}w_z\lambda)z$.

Indeed, this follows from Lemma \ref{daniel}: our data $A$, resp. $z_0$, resp. $z$ under discussion here play the role of the following data in Lemma \ref{daniel}: $V$, resp. $z$, resp. the origin of $V$.

Clearly we have $w_{x_{(z,Q)}}=w_z$ for any $(z,Q)\in {\mathfrak P}$. Similarly, for $(z,Q)\in
 {\mathfrak P}(\epsilon)$ we have $w_{x_{(z,Q)}}=w_{z}=w_{z_{\epsilon}}$ and one deduces that $x_{(z,Q)}=z_{\epsilon}$ and $\epsilon(z,Q)=\epsilon$. This proves (a). From (a) we immediately obtain that the ${\mathfrak P}(\epsilon)$ are pairwise disjoint, and from the above argument we also get that their union is ${\mathfrak P}$: statement (b).\hfill$\Box$\\

\begin{lem}\label{drehkey} Let $(z_1,Q_1)$, $(z_2,Q_2)\in{\mathfrak P}(\epsilon)$ with $Q_1\subset Q_2\subset Q(\epsilon)$. Let $H$ be a wall in $A$ such that $z_2\in H$ but
  $z_1\notin H$. Then $\sigma_0$ and $z_1$ ly on different
  components of $A-H$.
\end{lem}

{\sc Proof:} This is clear from our definitions.\hfill$\Box$\\ 

For $x\in A^0$ and $\lambda\in
X_*({\bf T})$ let ${\mathcal B}_{\lambda}(x)$ be the set of all $(w(-\lambda))x$ with $w\in W$, together with the set of all $s((w(-\lambda))x)$ with $w\in W$ and with $s$ denoting a reflection in an (affine) wall of $A$ which separates $x$ both from $(w(-\lambda))x$ and from the central chamber $\sigma_0$.\\

{\bf Definition:} A vertex coefficient system on $A^0$ over a ring $R$ is a
collection of data $$({\mathcal F}(x)_{x\in
  A^0},(T_{\lambda})_{{\lambda}\in\nabla})$$ as follows:

(a) for each $x\in
  A^0$ an $R$-module ${\mathcal F}(x)$. For a subset $Y$ of $A^0$ put
  ${\mathcal F}(Y)=\oplus_{y\in Y}{\mathcal F}(y)$.

(b) for each $\lambda\in\nabla$ an $R$-module endomorphism $T_{\lambda}$ of
${\mathcal F}(A^0)$ such that

(b1) for any two $\lambda_1,\lambda_2\in\nabla$ we have $T_{\lambda_1}\circ
T_{\lambda_2}=T_{\lambda_2}\circ T_{\lambda_1}$, and

(b2) for any $x\in A^0$ we have \begin{gather}T_{\lambda}({\mathcal F}(x))\subset{\mathcal
  F}({\mathcal B}_{{\lambda}}(x)).\label{hecdif}\end{gather}

Let ${\mathbb Z}[\nabla]$ denote the free ${\mathbb Z}$-module on the set
$\nabla$. Let $({\mathcal F}(x)_{x\in
  A^0},(T_{\lambda})_{{\lambda}\in\nabla})$ be a vertex coefficient system. Suppose we are given a collection $\alpha=(\alpha_{\lambda})_{{\lambda}\in\nabla}$
with $\alpha_{\lambda}\in R$. Consider for $0\le j\le |\nabla|-1$ the
map$${\mathcal F}(A^0)\otimes \bigwedge^{j+1}{\mathbb Z}[\nabla]\stackrel{{\delta}^j_{\alpha}}{\longrightarrow}{\mathcal F}(A^0)\otimes \bigwedge^{j}{\mathbb
  Z}[\nabla]$$$$f\otimes
(\lambda_1\wedge\cdots\wedge\lambda_{j+1})\mapsto\sum_{s=1}^{j+1}(-1)^s(T_{\lambda_s}(f)-\alpha_{\lambda_s}f)\otimes(\lambda_1\wedge\cdots\wedge\widehat{\lambda}_s\wedge\cdots\wedge\lambda_{j+1}).$$ As the
  $T_{\lambda}$ for $\lambda\in\nabla$ commute with each other we see that
    $\delta^j_{\alpha}\circ\delta^{j+1}_{\alpha}=0$. Writing $$C_{\alpha}{\mathcal
    F}^j={\mathcal F}(A^0)\otimes \bigwedge^{j}{\mathbb
  Z}[\nabla]$$for $j\ge0$ we obtain a Koszul complex $C_{\alpha}{\mathcal F}^{\bullet}$. (So only the differentials ${\delta}^j_{\alpha}$ depend on $\alpha$, not the individual terms $C_{\alpha}{\mathcal
    F}^j$.) Note that $C_{\alpha}{\mathcal
    F}^j=0$ for $j>|\nabla|$. For $\epsilon\in{\mathbb R}$ let$$C_{\alpha}{\mathcal
    F}^{\bullet}_{\epsilon}=\bigoplus_{(z,Q)\in{\mathfrak P}\atop\epsilon(z,Q)\le \epsilon}{\mathcal F}(z)\otimes
  \bigwedge_{\lambda\in Q}\lambda,$$$$C_{\alpha}{\mathcal
    F}^{\bullet}_{\epsilon-}=\bigoplus_{(z,Q)\in{\mathfrak P}\atop\epsilon(z,Q)<\epsilon}{\mathcal F}(z)\otimes
  \bigwedge_{\lambda\in Q}\lambda.$$The $C_{\alpha}{\mathcal
    F}_{\epsilon}^{\bullet}$ form a separated and exhausting filtration of $C_{\alpha}{\mathcal
    F}^{\bullet}$. Of course, the jumps in this
  filtration occur only at the discrete subset ${{\mathcal E}}$ of ${\mathbb R}_{\ge0}$, i.e. the graded pieces$$\overline{C}{\mathcal F}_{\epsilon}^{\bullet}=\frac{C_{\alpha}{\mathcal
    F}_{\epsilon}^{\bullet}}{C_{\alpha}{\mathcal
    F}_{\epsilon-}^{\bullet}}$$are zero if $\epsilon\notin {{\mathcal E}}$.

\begin{pro}\label{fide} Let $\epsilon\in{\mathbb R}$.

 (a) $C_{\alpha}{\mathcal
    F}^{\bullet}_{\epsilon}$ and $C_{\alpha}{\mathcal
    F}^{\bullet}_{\epsilon-}$ are subcomplexes of $C_{\alpha}{\mathcal
    F}^{\bullet}$, i.e. they are respected
 by the differential of $C_{\alpha}{\mathcal F}^{\bullet}$.

(b) The differential on the quotient
$\overline{C}{\mathcal F}_{\epsilon}^{\bullet}$ is independent of the collection ${\alpha}=(\alpha_{\lambda})_{{\lambda}\in\nabla}$.

(c) If $\epsilon\in{{\mathcal E}}$ we have $$\overline{C}{\mathcal F}_{\epsilon}^{\bullet}=\bigoplus_{(z,Q)\in{\mathfrak P}(\epsilon)}{\mathcal F}(z)\otimes
  \bigwedge_{\lambda\in Q}\lambda$$where the differentials on this complex are induced from those of $C_{\alpha}{\mathcal
    F}^{\bullet}$ by composing with the inclusion maps ${\mathcal F}(z)\to{\mathcal F}(A^0)$ on the one hand, and the  projection maps ${\mathcal F}(A^0)\to{\mathcal F}(z)$ on the other hand. 

\end{pro}

{\sc Proof:} Let $(z,Q)\in{\mathfrak P}$. For $f\otimes\wedge_{\lambda\in
  Q}\lambda\in{\mathcal F}(z)\otimes\wedge_{\lambda\in
  Q}\lambda$ with $|Q|=j+1$ write $$\delta^j_{\alpha}(f\otimes\wedge_{\lambda\in
  Q}\lambda)=\sum_{(z',Q')\in{\mathfrak P}}f_{(z',Q')}\otimes\wedge_{\lambda\in
  Q'}\lambda,$$listing only those $(z',Q')$ with $f_{(z',Q')}\ne 0$. For $z'=z$ we
  have $\epsilon(z',Q')<\epsilon(z,Q)$ since $Q'\subset Q$. For $z'\ne z$,
  formula (\ref{hecdif}) shows $z'\in{\mathcal B}_{{\lambda}}(x)$, so the
  definition of $\epsilon(z,Q)$ gives
  $\epsilon(z',Q')\le\epsilon(z,Q)$. Together this gives statement (a), while
  statement (b) follows already from $\epsilon(z',Q')<\epsilon(z,Q)$ for
  $z'=z$. For statement (c) just observe Lemma \ref{elemgeo}.\hfill$\Box$\\

\section{The Koszul complex ${C}{\mathcal F}_{V}^{\bullet}$ of the vertex
  coefficient system ${\mathcal F}_{V}$ }

Notice that $\sigma_0$ is the unique chamber in $X$ fixed by the
pro-$p$-Iwahori subgroup $I$. Let $\rho:X\to A$ denote the retraction from $X$ to $A$ centered in $\sigma_0$
(see \cite{bro} p.86), i.e. the unique map of poly-simplicial complexes $\rho:X\to A$
restricting to the identity on $A$ and with $\rho^{-1}(\sigma_0)=\sigma_0$. We also write $\rho:X^0\to A^0$ for the induced map between sets of special vertices. Notice that the sets $\rho^{-1}(x)$ for $x\in A^0$ are exactly the orbits for the action of $I$ on $X^0$. 

We fix $V$ and $\theta:{\mathfrak H}(V)\to R$ as before and we are going to assign to ${\mathcal
    C}(X^0,V)={\rm ind}_{\widetilde{K}Z}^{\widetilde{G}}V$ (formula (\ref{geodarst})) a vertex coefficient system ${\mathcal F}_{V}$ on $A^0$ (depending also on $\theta$, but we suppress this in our notation). For $x\in A^0$ put $${\mathcal F}_{V}(x)=\bigoplus_{y\in X^0\atop\rho(y)=x}{\mathcal C}(\{y\},V).$$

\begin{lem}\label{bspvertcoe} The endomorphisms $T_{\lambda}$ of ${\rm ind}_{\widetilde{K}Z}^{\widetilde{G}}V={\mathcal
    C}(X^0,V)={\mathcal F}_{V}(A^0)$ for $\lambda\in \nabla$ provide $({\mathcal F}_{V}(x)_{x\in
  A^0},(T_{\lambda})_{{\lambda}\in\nabla})$ with the structure of a vertex coefficient system on $A^0$.
\end{lem}

{\sc Proof:} We already know that the $T_{\lambda}$ for $\lambda\in\nabla$ commute with each other on $${\rm ind}_{\widetilde{K}Z}^{\widetilde{G}}V={\mathcal
    C}(X^0,V)=\oplus_{y\in X^0}{\mathcal C}(\{y\},V)=\oplus_{x\in
    A^0}{\mathcal F}_{V}(x)={\mathcal F}_{V}(A^0).$$We verify formula
  (\ref{hecdif}). From formula (\ref{hecfor}) it follows that\begin{gather}T_{\lambda}({\mathcal
    C}(\{gx_0\},V))\subset{\mathcal
    C}(gK(-\lambda)x_0,V)\label{hebu}\end{gather}for all $\lambda\in\nabla$, all
$g\in G$. Any $x\in A^0$ is of the form $x=tx_0$ for some $t\in T$ (Lemma
\ref{trans}). Consider the action of the stabilizer $tKt^{-1}$ of $x=tx_0$ in
$G$ on the set $tKt^{-1}\{(-\mu)x\,|\,\mu\in \nabla\}$. The intersection of an orbit $tKt^{-1}(-\mu)x$ with $A^0$ is
  $(W(-\mu))x$. Therefore\begin{gather}\rho(ItK(-\lambda)x_0)=\rho(tKt^{-1}(-\lambda)x)={\mathcal
    B}_{{\lambda}}(x).\label{jaeinf}\end{gather}Together\begin{align}T_{\lambda}({\mathcal F}_{V}(x))=T_{\lambda}({\mathcal C}(\rho^{-1}(x),V))&=T_{\lambda}({\mathcal C}(Itx_0,V))\\\notag{}&\stackrel{(\ref{hebu})}{\subset}{\mathcal C}(ItK(-\lambda)x_0,V))\notag\\{}&\stackrel{(\ref{jaeinf})}{=}{\mathcal
  F}_{V}({\mathcal B}_{{\lambda}}(x))\notag\end{align}where we used the formulae (\ref{hebu})
and (\ref{jaeinf}) as indicated.\hfill$\Box$\\

Suppose that we are given an $R$-algebra homomorphism $\chi:{\mathfrak H}_V(\widetilde{G},\widetilde{K})\to B$ such that the restriction $\chi|_{{\mathfrak H}(V)}:{\mathfrak H}(V)\to B$
of $\chi$ factors through an $R$-algebra homomorphism $\theta:{\mathfrak
  H}(V)\to R$. We put$$M_{\chi}(V)={\rm
  ind}_{\widetilde{K}}^{\widetilde{G}}V\otimes_{{\mathfrak
    H}_V(\widetilde{G},\widetilde{K}),\chi}B.$$Let us write $(.)_B=(.)\otimes_RB$. For example, $V_B$ is naturally a $B[\widetilde{K}]$-module. Let $\overline{\chi}$ denote the composition \begin{gather}R[T_{\lambda}]_{\lambda\in\nabla}\stackrel{1\otimes{\rm id}}{\longrightarrow}{\mathfrak H}(V)\otimes_RR[T_{\lambda}]_{\lambda\in\nabla}\cong {\mathfrak
    H}_V(\widetilde{G},\widetilde{K})\stackrel{\chi}{\longrightarrow}B.\label{chiquer}\end{gather} Put $\alpha=(\alpha_{\lambda})_{\lambda\in\nabla}$ with
  $\alpha_{\lambda}=\overline{\chi}(T_{\lambda})\in B$. Then we may alternatively write
  $$M_{\chi}(V)=\frac{{\rm
    ind}_{\widetilde{K}Z}^{\widetilde{G}}V_B}{\sum_{\lambda}(T_{\lambda}-\alpha_{\lambda})}$$(where ${\rm ind}_{\widetilde{K}Z}^{\widetilde{G}}V_B\cong({\rm ind}_{\widetilde{K}Z}^{\widetilde{G}}V)\otimes_RB$ is formed with respect to $\theta$, formula (\ref{isocind})).

Base extension $R\to B$ in Lemma \ref{bspvertcoe} provides us with a vertex
coefficient system ${\mathcal
  F}_{V_B}$ on $A^0$, so we may consider the associated Koszul complex $C_{\alpha}{\mathcal
  F}_{V_B}^{\bullet}$ defined in section
\ref{sectcox}.\\ 

{\bf Definition:} We say that the complex $C_{\alpha}{\mathcal
  F}_{V_B}^{\bullet}$ is exact if it is exact in $C_{\alpha}{\mathcal
  F}_{V_B}^{j}$ for all $j>0$, and similarly for the subcomplexes ${C}_{\alpha}{\mathcal
  F}_{V_B,\epsilon}^{\bullet}:={C}_{\alpha}({\mathcal
  F}_{V_B})_{\epsilon}^{\bullet}$ resp. the subquotient complexes $\overline{C}{\mathcal
  F}_{V_B,\epsilon}^{\bullet}:=\overline{C}({\mathcal
  F}_{V_B})_{\epsilon}^{\bullet}$ (some $\epsilon$) considered in section
\ref{sectcox}, and base extensions of them. We write $\overline{C}{\mathcal
  F}_{V,\epsilon}^{\bullet}=\overline{C}{\mathcal
  F}_{V_R,\epsilon}^{\bullet}$ (i.e. in the case $B=R$).

\begin{pro}\label{bellotwtrick} Suppose that the complexes $\overline{C}{\mathcal
  F}_{V,\epsilon}^{\bullet}\otimes_RR/{\mathfrak p}$ are exact, for any prime ideal ${\mathfrak p}$
  of $R$ (including ${\mathfrak p}=0$) and for any $\epsilon$. Then the $B$-module $M_{\chi}(V)$ is free. 
\end{pro}

{\sc Proof:} {\it Step 0:} Suppose that $B$ is a flat $R$-algebra or an $R/{\mathfrak p}$-algebra for a prime ideal ${\mathfrak p}$
  of $R$. Since formation of $\overline{C}{\mathcal
  F}_{V_B,\epsilon}^{\bullet}$ commutes with base changes in $B$ it follows from our hypothesis that $\overline{C}{\mathcal
  F}_{V_B,\epsilon}^{\bullet}$ is exact for all such $B$.

{\it Step 1:} For $B$ as in Step 0 we claim $${C}_{\alpha}{\mathcal
  F}_{V_B,\epsilon}^{0} \cap{\rm Im}[C_{\alpha}{\mathcal
  F}_{V_B}^{1}\to C_{\alpha}{\mathcal
  F}_{V_B}^{0}]={\rm Im}[{C}_{\alpha}{\mathcal
  F}_{V_B,\epsilon}^{1}\to C_{\alpha}{\mathcal
  F}_{V_B}^{0}]$$for all $\epsilon\in{\mathbb R}$. For this it is enough to
prove$${C}_{\alpha}{\mathcal
  F}_{V_B,\epsilon}^{0} \cap{\rm Im}[C_{\alpha}{\mathcal
  F}_{V_B,\epsilon'}^{1}\to C_{\alpha}{\mathcal
  F}_{V_B}^{0}]={\rm Im}[{C}_{\alpha}{\mathcal
  F}_{V_B,\epsilon}^{1}\to C_{\alpha}{\mathcal
  F}_{V_B}^{0}]$$for all $\epsilon'>\epsilon$. Since the set of
$\epsilon$ where the filtration $({C}_{\alpha}{\mathcal
  F}_{V_B,\epsilon}^{\bullet})_{\epsilon}$ of $C_{\alpha}{\mathcal
  F}_{V_B}^{\bullet}$ jumps is discrete, an induction shows that it is enough
to prove\begin{gather}{C}_{\alpha}{\mathcal
  F}_{V_B,\epsilon}^{0} \cap{\rm Im}[C_{\alpha}{\mathcal
  F}_{V_B,\epsilon'}^{1}\to C_{\alpha}{\mathcal
  F}_{V_B}^{0}]={\rm Im}[{C}_{\alpha}{\mathcal
  F}_{V_B,\epsilon}^{1}\to C_{\alpha}{\mathcal
  F}_{V_B}^{0}]\label{bello}\end{gather}for any $\epsilon'>\epsilon$ such that $\overline{C}{\mathcal
  F}_{V_B,\epsilon'}^{\bullet}={C}_{\alpha}{\mathcal
  F}_{V_B,\epsilon'}^{\bullet}/{C}_{\alpha}{\mathcal
  F}_{V_B,\epsilon}^{\bullet}$.

 The containment of the right hand side in the left hand side
  of (\ref{bello}) is clear. To see the converse, let $f\in C_{\alpha}{\mathcal
  F}_{V_B,\epsilon'}^{1}$ and suppose that the differential ${C}_{\alpha}{\mathcal
  F}_{V_B,\epsilon'}^{1}\to C_{\alpha}{\mathcal
  F}_{V_B}^{0}$ sends $f$ to an element in ${C}_{\alpha}{\mathcal
  F}_{V_B,\epsilon}^{0}$. This means that this image element vanishes in $\overline{C}{\mathcal
  F}_{V_B,\epsilon'}^{0}$. Therefore exactness of $\overline{C}{\mathcal
  F}_{V_B,\epsilon'}^{\bullet}$ in $\overline{C}{\mathcal
  F}_{V_B,\epsilon'}^{1}$, which holds true by Step 0, shows that the class $\overline{f}$ of
$f$ in $\overline{C}{\mathcal
  F}_{V_B,\epsilon'}^{1}$ is the image, under the differential $\overline{C}{\mathcal
  F}_{V_B,\epsilon'}^{2}\to \overline{C}{\mathcal
  F}_{V_B,\epsilon'}^{1}$, of an element $\overline{g}\in\overline{C}{\mathcal
  F}_{V_B,\epsilon'}^{2}$. Lift $\overline{g}$ to an element $g\in C_{\alpha}{\mathcal
  F}_{V_B,\epsilon'}^{2}$ and then subtract from $f$ the image of $g$ under the differential ${C}_{\alpha}{\mathcal
  F}_{V_B,\epsilon'}^{2}\to {C}_{\alpha}{\mathcal
  F}_{V_B,\epsilon'}^{1}$. Replacing $f$ by this difference we may now assume that the class $\overline{f}$ of
$f$ in $\overline{C}{\mathcal
  F}_{V_B,\epsilon'}^{1}$ vanishes. But this means that $f$ lies in ${C}_{\alpha}{\mathcal
  F}_{V_B,\epsilon}^{1}$. Our claim follows.

 {\it Step 2:} We may factor $\chi$ as $${\mathfrak
    H}_V(\widetilde{G},\widetilde{K})\stackrel{{\theta}^*}{\longrightarrow}R[T_{\lambda}]_{\lambda\in\nabla}\stackrel{\overline{\chi}}{\longrightarrow}B$$with
  $\overline{\chi}$ as in formla (\ref{chiquer}) and with ${\theta}^*$ as in formula (\ref{thetastern}). For $\epsilon\ge0$
put$$M_{\chi,\epsilon}(V)=\frac{{C}_{\alpha}{\mathcal
    F}_{V_B,\epsilon}^{0}}{{C}_{\alpha}{\mathcal
    F}_{V_B,\epsilon}^{0}\cap{\rm Im}[C_{\alpha}{\mathcal
  F}_{V_B}^{1}\to C_{\alpha}{\mathcal
  F}_{V_B}^{0}]}.$$Since ${\mathfrak
    H}_V(\widetilde{G},\widetilde{K})$ is generated as an ${\mathfrak
    H}(V)$-algebra by the $T_{\lambda}$ for $\lambda\in\nabla$ (Proposition
  \ref{redsat}) it follows that $M_{\chi,\epsilon}(V)$
is the image in $M_{\chi}(V)$ of the composition $${C}_{\alpha}{\mathcal
    F}_{V_B,\epsilon}^{0}\longrightarrow {C}_{\alpha}{\mathcal
    F}_{V_B}^{0}={\rm
  ind}_{\widetilde{K}Z}^{\widetilde{G}}V_B\cong$$$${\rm ind}_{\widetilde{K}}^{\widetilde{G}}V\otimes_{{\mathfrak
    H}_V(\widetilde{G},\widetilde{K}),{\theta}^*}R[T_{\lambda}]_{\lambda\in\nabla}\otimes_RB\stackrel{{\rm id}\otimes\overline{\chi}\otimes {\rm id}}{\longrightarrow} {\rm
  ind}_{\widetilde{K}}^{\widetilde{G}}V\otimes_{{\mathfrak
    H}_V(\widetilde{G},\widetilde{K}),\chi}B=M_{\chi}(V).$$Again, since the set of
$\epsilon$ where the filtration $({C}_{\alpha}{\mathcal
  F}_{V_B,\epsilon}^{\bullet})_{\epsilon}$ of $C_{\alpha}{\mathcal
  F}_{V_B}^{\bullet}$ jumps is discrete, to show that
$M_{\chi}(V)$ is free over $B$ it is therefore enough to show that the $B$-modules
$M_{\chi,\epsilon'}(V)/M_{\chi,\epsilon}(V)$ are free, for any
$\epsilon'>\epsilon$ such that $\overline{C}{\mathcal
  F}_{V,\epsilon'}^{\bullet}={C}_{\alpha}{\mathcal
  F}_{V_B,\epsilon'}^{\bullet}/{C}_{\alpha}{\mathcal
  F}_{V_B,\epsilon}^{\bullet}$. 

{\it Step 3:} Let again $B$ be as in Step 0. We have$$M_{\chi,\epsilon'}(V)=\frac{{C}_{\alpha}{\mathcal
    F}_{V_B,\epsilon'}^{0}}{{C}_{\alpha}{\mathcal
    F}_{V_B,\epsilon'}^{0}\cap{\rm Im}[C_{\alpha}{\mathcal
  F}_{V_B}^{1}\to C_{\alpha}{\mathcal
  F}_{V_B}^{0}]}=\frac{{C}_{\alpha}{\mathcal
    F}_{V_B,\epsilon'}^{0}}{{\rm Im}[C_{\alpha}{\mathcal
  F}_{V_B,\epsilon'}^{1}\to C_{\alpha}{\mathcal
  F}_{V_B}^{0}]}$$where the second equality results from Step 1. Hence \begin{align}\frac{M_{\chi,\epsilon'}(V)}{M_{\chi,\epsilon}(V)}&\cong \frac{{C}_{\alpha}{\mathcal
    F}_{V_B,\epsilon'}^{0}}{{C}_{\alpha}{\mathcal
    F}_{V_B,\epsilon}^{0}+ {\rm Im}[C_{\alpha}{\mathcal
  F}_{V_B,\epsilon'}^{1}\to C_{\alpha}{\mathcal
  F}_{V_B}^{0}]}\notag\\{}&\cong\frac{\overline{C}{\mathcal
    F}_{V_B,\epsilon'}^{0}}{{\rm Im}[\overline{C}{\mathcal
  F}_{V_B,\epsilon'}^{1}\to \overline{C}{\mathcal
  F}_{V_B,\epsilon'}^{0}]}\notag\\{}&\cong\frac{\overline{C}{\mathcal
    F}_{V,\epsilon'}^{0}}{{\rm Im}[\overline{C}{\mathcal
  F}_{V,\epsilon'}^{1}\to \overline{C}{\mathcal
  F}_{V,\epsilon'}^{0}]}\otimes_RB.\notag\end{align}In particular, for any $R$-algebra morphism
$b:B\to B'$ with both $B$ and $B'$ as in Step 0, and for any character
$\chi':{\mathfrak H}_{V}(G,K)\to B'$ such that $\chi'|_{{\mathfrak H}(V)}=(b\circ\chi)|_{{\mathfrak H}(V)}$ we have a natural
isomorphism of
$B'$-modules$$\frac{M_{\chi,\epsilon'}(V)}{M_{\chi,\epsilon}(V)}\otimes_BB'\cong\frac{M_{\chi',\epsilon'}(V)}{M_{\chi',\epsilon}(V)}.$$We
claim that the $B$-modules $M_{\chi,\epsilon'}(V)/M_{\chi,\epsilon}(V)$ are free, for any $\epsilon'>\epsilon$ such that $\overline{C}{\mathcal
  F}_{V,\epsilon'}^{\bullet}={C}_{\alpha}{\mathcal
  F}_{V_B,\epsilon'}^{\bullet}/{C}_{\alpha}{\mathcal
  F}_{V_B,\epsilon}^{\bullet}$. By what we just said, it is enough to show this if $R=B$. Thus, if $R$ is a
field we are done. If $R$ is a principal ideal domain we argue as follows. Since $M_{\chi,\epsilon'}(V)/M_{\chi,\epsilon}(V)\cong{\overline{C}{\mathcal
    F}_{V,\epsilon'}^{0}}/{{\rm Im}[\overline{C}{\mathcal
  F}_{V,\epsilon'}^{1}\to \overline{C}{\mathcal
  F}_{V,\epsilon'}^{0}}]$ is finitely generated it
is enough to show$${\rm dim}_{{\rm Q}(R)}(\frac{\overline{C}{\mathcal
    F}_{V,\epsilon'}^{0}}{{\rm Im}[\overline{C}{\mathcal
  F}_{V,\epsilon'}^{1}\to \overline{C}{\mathcal
  F}_{V,\epsilon'}^{0}]}\otimes_R{\rm Q}(R))={\rm dim}_{R/{\mathfrak p}}(\frac{\overline{C}{\mathcal
    F}_{V,\epsilon'}^{0}}{{\rm Im}[\overline{C}{\mathcal
  F}_{V,\epsilon'}^{1}\to \overline{C}{\mathcal
  F}_{V,\epsilon'}^{0}]}\otimes_R R/{\mathfrak p})$$for all ${\mathfrak p}$, where ${\rm Q}(R)$
denotes the quotient field of $R$. Since $\overline{C}{\mathcal
  F}_{V,\epsilon'}^{\bullet}\otimes_R{\rm Q}(R)$ and $\overline{C}{\mathcal
  F}_{V,\epsilon'}^{\bullet}\otimes_RR/{\mathfrak p}$ are exact this boils
down to showing $$\sum_{j\ge0}(-1)^j{\rm dim}_{{\rm Q}(R)}(\overline{C}{\mathcal
  F}_{V,\epsilon'}^{j}\otimes_R{\rm Q}(R))=\sum_{j\ge0}(-1)^j{\rm dim}_{R/{\mathfrak p}}(\overline{C}{\mathcal
  F}_{V,\epsilon'}^{j}\otimes_RR/{\mathfrak p}).$$ But for all $j$ the $\overline{C}{\mathcal
  F}_{V,\epsilon'}^{j}$ are $R$-free, hence $${\rm dim}_{{\rm Q}(R)}(\overline{C}{\mathcal
  F}_{V,\epsilon'}^{j}\otimes{\rm Q}(R))={\rm dim}_{R/{\mathfrak p}}(\overline{C}{\mathcal
  F}_{V,\epsilon'}^{j}\otimes R/{\mathfrak p})$$and this implies what we want.

{\it Step 4:} For $B$ as in Step 0 it follows already that $M_{\chi}(V)$ is
$B$-free. In particular this follows for $B=R[T_{\lambda}]_{\lambda\in\nabla}$ and $\chi=\theta^*:{\mathfrak
    H}_V(\widetilde{G},\widetilde{K})\to R[T_{\lambda}]_{\lambda\in\nabla}$ for any $\theta:{\mathfrak H}(V)\to R$ satisfying our hypotheses. But $M_{\chi}(V)$ for a general $B$ and $\chi$ is obtained from such $B$ and $\chi$ by base change. \hfill$\Box$\\

{\bf Remark:} The proof shows that if $R$ is a field then we need exactness of
the $\overline{C}{\mathcal F}_{V,\epsilon}^{\bullet}$ only in $\overline{C}{\mathcal F}_{V,\epsilon}^{1}$.  

\section{The subquotient complexes  ${\mathcal K}^{\bullet}_D(V)$ of  ${C}{\mathcal F}_{V}^{\bullet}$}

Choose once and for all an enumeration
$\lambda_1,\ldots,\lambda_{|\nabla|}$ of $\nabla$. For a subset $E$ of
$\nabla$ and $\lambda\in \nabla-E$ define $\sigma(\lambda,E)\in\{\pm 1\}$ as follows. Write
$\nabla-E=\{\lambda_{i_1},\ldots,\lambda_{i_{i+1}}\}$ with
${i_1}<\ldots<{i_{i+1}}$ and if $\lambda=\lambda_{i_s}$ put $\sigma(\lambda,E)=(-1)^s$. \\

Recall that via the
splitting (\ref{monoiex}) (or rather its $\otimes{\mathbb Q}$-version) we read $X_*({{\bf T}})$ as a subset
of $X_*({\widetilde{\bf T}})$. For a subset $E$ of $\nabla (\subset X_*({{\bf T}}))$ we may
thus define$$\gamma_{[E]}=\prod_{\lambda\in
  E}\lambda(p_F)\in \widetilde{T}.$$The projection
$\widetilde{T}\to X_*({{\bf T}})$ sends $\gamma_{[E]}$ to$$\lambda_{E}=\sum_{\lambda\in
  E}\lambda\in X_*({{\bf T}}).$$For a subset $E$ of $\nabla$ we
put $${{{{\mathfrak N}}}}[E]=\gamma_{[E]}^{-1}{{{{\mathfrak
        N}}}}\gamma_{[E]}.$$Notice that $\gamma_{[E]}\in\widetilde{T}_+$ and hence ${{{{\mathfrak N}}}}[E]\subset {{{{\mathfrak N}}}}[E']\subset
{{{{\mathfrak N}}}}[\emptyset]={{{{\mathfrak N}}}}$ for $E'\subset E\subset \nabla$. For $\lambda\in
X_*({{\bf T}})$ we write $\gamma_{[\lambda]}=\gamma_{[\{\lambda\}]}$ and ${\mathfrak N}[\lambda]={\mathfrak N}[\{\lambda\}]$. For a subset $E$ of $\nabla$ we define the $R[{{{{\mathfrak N}}}}]$-module$${\rm
  ind}^{{{{{\mathfrak N}}}}}_{{{{{\mathfrak N}}}}[E]}(V)=\{f:{{{{\mathfrak N}}}}\longrightarrow
  V\,\,|\,\,f(i_1i_2)=(\gamma_{[E]}i_2^{-1}\gamma_{[E]}^{-1})f(i_1)\mbox{
    for }i_1\in {{{{\mathfrak N}}}}, i_2\in {{{{\mathfrak N}}}}[E]\}$$with ${{{{\mathfrak N}}}}$-action
  given by $(if)(i')=f(i^{-1}i')$ for $i, i'\in{{{{\mathfrak N}}}}$.

\begin{lem}\label{petersfrage} (a) Suppose that $\lambda\in\nabla$. We then have ${\bf
  N}_{-\lambda}({\mathcal O}_F)\subset {\mathfrak N}$, the composition ${\bf
  N}_{-\lambda}({\mathcal O}_F)\to{\mathfrak N}\to{\mathfrak N}/{\mathfrak N}[\lambda]$ is surjective and $${\mathfrak N}[\lambda]\cap{\bf
  N}_{-\lambda}({\mathcal O}_F)\subset{\rm Ker}[{\bf
  N}_{-\lambda}({\mathcal O}_F)\longrightarrow{\bf
  N}_{-\lambda}(k_F)].$$In particular, ${\mathfrak N}/{\mathfrak N}[\lambda]$ naturally surjects onto ${\bf
  N}_{-\lambda}(k_F)$. 

(b) For $\lambda\in\nabla-E$ the map  
$$T_{\lambda}:{\rm
  ind}^{{{{{\mathfrak N}}}}}_{{{{{\mathfrak N}}}}[E]}(V)\longrightarrow{\rm
  ind}^{{{{{\mathfrak N}}}}}_{{{{{\mathfrak N}}}}[E\cup\{\lambda\}]}(V)$$$$f\mapsto[i\mapsto
 \sigma(\lambda,E)\xi_{\lambda}(f(i))]$$is well defined and ${{{{\mathfrak
         N}}}}$-equivariant.
\end{lem} 

{\sc Proof:} (a) follows easily from the definitions (that of ${\bf
  N}_{-\lambda}$ has been given along with the definition of the map
$\xi_{\lambda}$). (Warning: in general, ${\mathfrak N}[\lambda]$ is {\it
  not} normal in ${\mathfrak N}$.) In (b) we need to check
that $\xi_{\lambda}\circ f$ for $f\in {\rm
  ind}^{{{{{\mathfrak N}}}}}_{{{{{\mathfrak N}}}}[E]}(V)$ is indeed an element
of ${\rm
  ind}^{{{{{\mathfrak N}}}}}_{{{{{\mathfrak N}}}}[E\cup\{\lambda\}]}(V)$. This may be reduced to the well
definedness of the operators $T_{\lambda}$ on ${\mathcal C}(X^0,V)$, as
follows. The assignment $i\mapsto i\gamma_{[E]}^{-1}$ induces a
bijection ${{{{\mathfrak N}}}}/{{{{\mathfrak N}}}}[E]\cong {{{{\mathfrak
        N}}}}\gamma_{[E]}^{-1}\widetilde{K}/\widetilde{K}$ and then an isomorphism$${\rm
  ind}^{{{{{\mathfrak N}}}}}_{{{{{\mathfrak N}}}}[E]}(V)\cong {\mathcal C}({{{{\mathfrak
        N}}}}\lambda_E^{-1}x_0,V)$$(of course $\lambda_E^{-1}x_0=\gamma_{[E]}^{-1}x_0$), and similarly for $E\cup\{\lambda\}$
instead of $E$. Tracing back our definitions we then see that the map in
question is nothing but the map $${\mathcal C}({{{{\mathfrak
        N}}}}\lambda_E^{-1}x_0,V)\longrightarrow {\mathcal C}({{{{\mathfrak
        N}}}}\lambda_{E\cup\{\lambda\}}^{-1}x_0,V)$$ induced from the
endomorphism $T_{\lambda}$ of ${\mathcal C}(X^0,V)$ via the inclusion of ${\mathcal C}({{{{\mathfrak
        N}}}}\lambda_E^{-1}x_0,V)$ into ${\mathcal C}(X^0,V)$ (extension of
functions by zero) and the projection of ${\mathcal C}(X^0,V)$ onto ${\mathcal C}({{{{\mathfrak
        N}}}}\lambda_{E\cup\{\lambda\}}^{-1}x_0,V)$ (restriction of
functions). \hfill$\Box$\\

{\bf Remark:} Statement Lemma \ref{petersfrage} (b) is in fact equivalent
with the following statement: For $\lambda\in\nabla-E$ the $R$-linear endomorphism $\xi_{\lambda}$ of $V$ is ${{{{\mathfrak
        N}}}}[E\cup\{\lambda\}]$-equivariant if ${{{{\mathfrak
        N}}}}[E\cup\{\lambda\}]$ acts on the source $V$ through the map\begin{gather}{{{{\mathfrak
        N}}}}[E\cup\{\lambda\}]\subset {{{{\mathfrak N}}}}[E]=\gamma_{[E]}^{-1}{{{{\mathfrak N}}}}\gamma_{[E]}\stackrel{\cong}{\longrightarrow}
{{{{\mathfrak N}}}}\subset\widetilde{K},\quad\quad \gamma_{[E]}^{-1}i\gamma_{[E]}\mapsto
i\end{gather}composed with the given $\widetilde{K}$-action on $V$, and on the
target $V$ through the map\begin{gather}{{{{\mathfrak
        N}}}}[E\cup\{\lambda\}]\stackrel{\cong}{\longrightarrow}
{{{{\mathfrak N}}}}\subset\widetilde{K},\quad\quad \gamma_{[E\cup\{\lambda\}]}^{-1}i\gamma_{[E\cup\{\lambda\}]}\mapsto
i\end{gather} composed with the given $\widetilde{K}$-action on $V$.

Let us give the direct argument in case (II). We may assume
$E=\emptyset$. By construction, $\xi_{\lambda}$ is equivariant for the natural
action of the subgroup $\widetilde{\bf M}_{\lambda}({\mathcal O}_F)\cdot
\widetilde{K}_1=\widetilde{\bf M}_{-\lambda}({\mathcal
  O}_F)\cdot\widetilde{K}_1$ of $\widetilde{K}$ on $V$, where
$\widetilde{K}_1={\rm ker}[\widetilde{K}\to \widetilde{\bf G}(k_F)]$. Now ${{{{\mathfrak
        N}}}}[\lambda]$ is contained in this subgroup (cf.  Lemma \ref{petersfrage} (a)). Therefore, as the image of $\xi_{\lambda}$ is contained in $V^{{\bf
  N}_{-\lambda}(k_F)}$, it is enough to show that for any $i\in {{{{\mathfrak
        N}}}}$ the elements $\lambda(p_F)^{-1}i\lambda(p_F)$ and $i$ of $\widetilde{\bf
  M}_{-\lambda}({\mathcal O}_F){\bf
  N}_{-\lambda}({\mathcal O}_F)\cdot\widetilde{K}_1$ map to the same element in $\widetilde{\bf
  M}_{-\lambda}(k_F)$. This is shown in \cite{her} Proposition
3.8.\hfill$\Box$\\ 

For subset $D$ of $\nabla$ and for $j\ge0$ define the $R[{{{{\mathfrak N}}}}]$-module$${\mathcal
  K}^j_D(V)=\bigoplus_{E\subset D\atop|E|=j}{\rm
  ind}^{{{{{\mathfrak N}}}}}_{{{{{\mathfrak N}}}}[E]}(V).$$Using as differentials the sums of the maps $T_{\lambda}:{\rm
  ind}^{{{{{\mathfrak N}}}}}_{{{{{\mathfrak N}}}}[E]}(V)\to{\rm
  ind}^{{{{{\mathfrak N}}}}}_{{{{{\mathfrak N}}}}[E\cup\{\lambda\}]}(V)$ of
Lemma \ref{petersfrage} we obtain a
complex of $R[{{{{\mathfrak N}}}}]$-modules\begin{gather}{\mathcal K}^{\bullet}_D(V)=[0\longrightarrow
  {\mathcal K}^{0}_D(V)\longrightarrow{\mathcal
    K}^{1}_D(V)\longrightarrow{\mathcal
    K}^{2}_D(V)\longrightarrow\ldots\longrightarrow {\mathcal
    K}^{|D|}_D(V)].\label{redgopad}\end{gather}It will be formally convenient to put ${\mathcal
  K}^j_D(V)=0$ for $j<0$.\\

{\bf Definition:} We say that the complex ${\mathcal
  K}^{\bullet}_D(V)$ is exact if it is exact in ${\mathcal
  K}^j_D(V)$ for all $j<|D|$. 

\begin{lem}\label{borelred} If the complex ${\mathcal K}^{\bullet}_{\nabla}(V)$ is exact, then
  the complexes ${\mathcal K}^{\bullet}_D(V)$ are exact for any $D\subset\nabla$.
\end{lem}

{\sc Proof:} We argue by descending induction on $|D|$. It is enough to
show that if ${\mathcal K}^{\bullet}_D(V)$ is exact and if $D'=D-\{\lambda\}$ with
$\lambda\in D$, then ${\mathcal K}^{\bullet}_{D'}(V)$ is exact. Consider the decomposition\begin{align}{\mathcal
  K}^j_D(V)=\bigoplus_{E\subset D\atop|E|=j}{\rm
  ind}^{{{{{\mathfrak N}}}}}_{{{{{\mathfrak N}}}}[E]}V&=(\bigoplus_{\lambda\notin E\subset
  D\atop|E|=j}{\rm
  ind}^{{{{{\mathfrak N}}}}}_{{{{{\mathfrak N}}}}[E]}V)\bigoplus(\bigoplus_{\lambda\in E\subset
  D\atop|E|=j}{\rm
  ind}^{{{{{\mathfrak N}}}}}_{{{{{\mathfrak N}}}}[E]}V)\notag\\{}&={\mathcal
  K}^j_{D'}(V)\bigoplus(\bigoplus_{\lambda\in E\subset
  D\atop|E|=j}{\rm
  ind}^{{{{{\mathfrak N}}}}}_{{{{{\mathfrak N}}}}[E]}V).\notag\end{align}If $j$ varies the direct summands $\oplus_{\lambda\in E\subset
  D\atop|E|=j}{\rm
  ind}^{{{{{\mathfrak N}}}}}_{{{{{\mathfrak N}}}}[E]}V$ define a subcomplex of ${\mathcal
  K}^{\bullet}_{D}(V)$ whose quotient complex is isomorphic with ${\mathcal
  K}^{\bullet}_{D'}(V)$. We may identify this subcomplex as follows. We choose representatives $u\in
{{{{\mathfrak N}}}}$ for the cosets $u{{{{\mathfrak N}}}}[\lambda]\in {{{{\mathfrak N}}}}/{{{{\mathfrak N}}}}[\lambda]$. For
$\lambda\in E\subset D$ we then have
the bijections of left cosets \begin{align}u{{{{\mathfrak N}}}}[\lambda]/{{{{\mathfrak N}}}}[E]&\cong
{{{{\mathfrak N}}}}[\lambda]/{{{{\mathfrak N}}}}[E]\notag\\{}&\cong {{{{\mathfrak N}}}}/{{{{\mathfrak N}}}}[E']\quad\quad\mbox{ with }E'=E-\{\lambda\},\notag\end{align}the first one
induced by multiplication with $u^{-1}$, the second one by conjugation with $\lambda(p_F)$. As we may view ${{{{\mathfrak N}}}}/{{{{\mathfrak N}}}}[E]$ as the
disjoint union of the $u{{{{\mathfrak N}}}}[\lambda]/{{{{\mathfrak N}}}}[E]$ for all
$u{{{{\mathfrak N}}}}[\lambda]\in {{{{\mathfrak N}}}}/{{{{\mathfrak N}}}}[\lambda]$ this provides isomorphisms \begin{align}\bigoplus_{\lambda\in E\subset
  D\atop|E|=j}{\rm
  ind}^{{{{{\mathfrak N}}}}}_{{{{{\mathfrak N}}}}[E]}V&\cong  \bigoplus_{E'\subset
  D'\atop|E'|=j-1}\bigoplus_{{{{{\mathfrak N}}}}/{{{{\mathfrak N}}}}[\lambda]}{\rm
  ind}^{{{{{\mathfrak N}}}}}_{{{{{\mathfrak N}}}}[E']}V\notag\\{}&=\bigoplus_{{{{{\mathfrak N}}}}/{{{{\mathfrak N}}}}[\lambda]}{\mathcal
  K}^{j-1}_{D'}(V).\notag\end{align}These identifications are compatible
with the differentials in ${\mathcal
  K}^{\bullet}_{D}(V)$ and ${\mathcal
  K}^{\bullet-1}_{D'}(V)$. The conclusion is that we have an exact sequence of
complexes$$0\longrightarrow \bigoplus_{{{{{\mathfrak N}}}}/{{{{\mathfrak N}}}}[\lambda]}{\mathcal
  K}^{\bullet-1}_{D'}(V)\longrightarrow {\mathcal
  K}^{\bullet}_{D}(V)\longrightarrow {\mathcal
  K}^{\bullet}_{D'}(V)\longrightarrow0.$$In cohomology we obtain the long exact
sequence$$\ldots\longrightarrow h^{s-1}({\mathcal
  K}^{\bullet}_{D}(V))\longrightarrow h^{s-1}({\mathcal
  K}^{\bullet}_{D'}(V))\longrightarrow h^s(\bigoplus_{{{{{\mathfrak N}}}}/{{{{\mathfrak N}}}}[\lambda]}{\mathcal
  K}^{\bullet-1}_{D'}(V))\longrightarrow h^{s}({\mathcal
  K}^{\bullet}_{D}(V))\longrightarrow\ldots.$$Exactness of ${\mathcal
  K}^{\bullet}_D(V)$ means $ h^{s}({\mathcal
  K}^{\bullet}_{D}(V))=0$ for all $0\le s< |D|$. The long exact cohomology
sequence therefore implies $$h^{s}({\mathcal
  K}^{\bullet}_{D'}(V))\cong h^{s+1}(\bigoplus_{{{{{\mathfrak N}}}}/{{{{\mathfrak N}}}}[\lambda]}{\mathcal
  K}^{\bullet-1}_{D'}(V))\cong\bigoplus_{{{{{\mathfrak N}}}}/{{{{\mathfrak N}}}}[\lambda]} h^s({\mathcal
  K}^{\bullet}_{D'}(V))$$ for all $0\le s< |D|-1=|D'|$. As
$|{{{{\mathfrak N}}}}/{{{{\mathfrak N}}}}[\lambda]|>1$ this implies $h^{s}({\mathcal
  K}^{\bullet}_{D'}(V))=0$ for all $0\le s< |D'|$, i.e. exactness of ${\mathcal
  K}^{\bullet}_{D'}(V)$.\hfill$\Box$\\

 For $x\in X^0$ let $K_x\subset G$ denote the maximal compact
open subgroup fixing $x$. If more specifically $x\in A^0$ then for each $\alpha\in\Phi$ the set of special vertixes in $A$ which are fixed by $U_{\alpha}\cap K_x$ is the set of special vertices of a halfspace in $A$; its boundary is a wall containing $x$ (and each wall through $x$ is obtained in this way). 

For $y\in A^0$ let $I_y=K_y\cap I$, the stabilizer of $y$ in $I$. It is generated by its subgroups $I_y\cap
 U_{\alpha}$ for $\alpha\in\Phi$. For $y,z\in A^0$ let $\Phi_y^z$ denote the
 set of all $\alpha\in\Phi$ such that $U_{\alpha}\cap K_y$ is {\it not}
 contained in $K_z$, i.e. does not fix $z$. Let $U_y^z$ denote the subgroup of
 $G$ generated by the $U_{\alpha}\cap K_y$ for all $\alpha\in \Phi_y^z$. 

\begin{lem}\label{lemfix} Suppose that for each wall $H$ in
  $A$ with $y\in H$ but $z\notin H$ we have: $z$ and
  $\sigma_0$ do not ly in the same component of $A-H$. Then $I_y.z=U_y^z.z$.
\end{lem}

{\sc Proof:} This is clear as $\sigma_0$ is the unique chamber fixed by $I$.\hfill$\Box$\\  

\begin{pro}\label{marialichtmessneu} For any $\epsilon\in{\mathcal E}$ there exists a subset $D$ of $\nabla$ such that the complex $\overline{C}{\mathcal
  F}_{V,\epsilon}^{\bullet}$ is isomorphic with a direct sum of complexes of the form ${\mathcal K}^{|D|-\bullet}_D(V)$.
\end{pro}

{\sc Proof:} (1) Let $E\subset\nabla$. Each wall $H$ in $A$ containing $x_0$ but not containing
$(-\lambda_E)x_0$ separates $(-\lambda_E)x_0$ from $\sigma_0$. Moreover, for each $\alpha\in \Phi_{x_0}^{(-\lambda_E)x_0}$ we
have $U_{\alpha}\cap K_{x_0}=U_{\alpha}\cap K\subset{\mathfrak N}$, hence
$U_{x_0}^{(-\lambda_E)x_0}\subset {\mathfrak N}$. It therefore follows from Lemma \ref{lemfix} that \begin{gather}I.(-\lambda_E)x_0=I_{x_0}.(-\lambda_E)x_0=U_{x_0}^{(-\lambda_E)x_0}.(-\lambda_E)x_0={\mathfrak N}.(-\lambda_E)x_0.\label{fixzei}\end{gather} Since ${\mathfrak
   N}[E]$ is the stabilizer of $(-\lambda_E)x_0$ in ${\mathfrak
   N}$ we obtain an isomorphism of ${\mathfrak N}$-representations\begin{gather}{\rm
  ind}^{{{{{\mathfrak N}}}}}_{{{{{\mathfrak N}}}}[E]}(V)\cong {\mathcal
  C}(I.(-\lambda_E)x_0,V).\label{sprung}\end{gather}Now let
$D\subset\nabla$. Our definitions are made in such a way that the isomorphisms
(\ref{sprung}), tensored on the right hand side with (the free rank one module generated by) $\bigwedge_{\lambda\in D-E}\lambda$, induce an isomorphism of complexes$${\mathcal
  K}^{\bullet}_D(V)\cong \bigoplus_{E\subset D\atop |E|=\bullet} {\mathcal C}(I.(-\lambda_E)x_0,V)\otimes\bigwedge_{\lambda\in D-E}\lambda,$$the right hand side viewed as a subquotient complex of ${C}_{\alpha}{\mathcal
    F}_V^{|D|-\bullet}$ (i.e. with differentials obtained from those of ${C}_{\alpha}{\mathcal
    F}_V^{|D|-\bullet}$ by composing them with the natural inclusions ${\mathcal C}(I.(-\lambda_E)x_0,V)\to {\mathcal C}(A^0,V)$ (extension of functions by $0$) resp. projections $ {\mathcal C}(A^0,V)\to {\mathcal C}(I.(-\lambda_E)x_0,V)$ (restriction of
  functions)). 

(2) Put $y_{\epsilon}=\sum_{\lambda\in
  Q(\epsilon)}(w_{z_{\epsilon}}\lambda)z_{\epsilon}$. For $i_1, i_2\in I$ and
$(z,Q)\in {\mathfrak P}(\epsilon)$ with $i_1z=i_2z$ we also have
$i_1y_{\epsilon}=i_2y_{\epsilon}$. Indeed, one may deduce from Lemma
\ref{drehkey} that each minimal gallery in the affine Coxeter complex $A$ connecting the chamber $\sigma_0$ with the vertex $i_1z=i_2z$ must contain $i_1y_{\epsilon}$ and $i_2y_{\epsilon}$, but then necessarily $i_1y_{\epsilon}=i_2y_{\epsilon}$, as desired. It follows that $$\overline{C}{\mathcal
  F}_{V,\epsilon}^{\bullet}\cong \bigoplus_{i\in I/I_{y_{\epsilon}}}(\bigoplus_{(z,Q)\in{\mathfrak P}(\epsilon)} {\mathcal C}(iI_{y_{\epsilon}}.z,V)\otimes\bigwedge_{\lambda\in Q}\lambda).$$The proof of Lemma \ref{bspvertcoe} (i.e. formula (\ref{hebu})) shows that the differential of $\overline{C}{\mathcal
  F}_{V,\epsilon}^{\bullet}$ respects this $I/I_{y_{\epsilon}}$-indexed decomposition. Therefore the complex $\overline{C}{\mathcal
  F}_{V,\epsilon}^{\bullet}$ is isomorphic with the direct sum, indexed by $i\in
  I/I_{y_{\epsilon}}$, of its subcomplexes\begin{gather}\bigoplus_{(z,Q)\in{\mathfrak
        P}(\epsilon)} {\mathcal C}(iI_{y_{\epsilon}}.z,V)\otimes\bigwedge_{\lambda\in
      Q}\lambda.\label{lomo}\end{gather}These complexes are all isomorphic
  (via $i$), so it is enough to show that the complex (\ref{lomo}) with $i=1$ is
  isomorphic with ${\mathcal K}^{|D|-\bullet}_D(V)$ for some $D$.

(3) Put $D=Q(\epsilon)$. Sending $(\sum_{\lambda\in Q}(w_{z_{\epsilon}}\lambda)z_{\epsilon},Q)$ to $(-\lambda_{D-Q})x_0$ is a bijection \begin{gather}{\mathfrak P}(\epsilon)\cong\{(-\lambda_E)x_0\,|\,E\subset
D\}.\label{baniso}\end{gather} It sends $(z_{\epsilon},\emptyset)$ to
$(-\lambda_D)x_0$ and sends $(y_{\epsilon},Q(\epsilon))$ to $x_0$. Let $t\in
\widetilde{T}$ such that $ty_{\epsilon}=x_0$. Then the bijection (\ref{baniso}) can be described by the formula $w_{z_{\epsilon}}^{-1}tz=(-\lambda_{D-Q})x_0$ for all $(z,Q)\in{\mathfrak P}(\epsilon)$. Let ${n}\in
N(\widetilde{T})$ denote a lift of $w_{z_{\epsilon}}^{-1}\in W=N(\widetilde{T})/\widetilde{T}$
and put $g=nt\in\widetilde{G}$. Then the bijection (\ref{baniso}) can also be described by the formula $gz=(-\lambda_{D-Q})x_0$ for all $(z,Q)\in{\mathfrak P}(\epsilon)$. For any $(z,Q)\in {\mathfrak P}(\epsilon)$ we infer from Lemma \ref{drehkey} and Lemma \ref{lemfix} that $I_{y_{\epsilon}}.z=U_{y_{\epsilon}}^{z}.z$. On the other hand we have $I.(-\lambda_{D-Q})x_0=U_{x_0}^{(-\lambda_{D-Q})x_0}.(-\lambda_{D-Q})x_0$ as noticed
in formula (\ref{fixzei}) above. Together it follows that $g$ induces for any $(z,Q)\in {\mathfrak P}(\epsilon)$ a bijection of orbits
$I_{y_{\epsilon}}.z\cong I.(-\lambda_{D-Q})x_0$. But then $g$ also induces an isomorphism ${\mathcal C}(I_{y_{\epsilon}}.z,V)\cong {\mathcal C}(I.(-\lambda_{D-Q})x_0,V)$. It is clear that these fit together into an isomorphism of complexes
$$\bigoplus_{(z,Q)\in{\mathfrak
        P}(\epsilon)} {\mathcal C}(I_{y_{\epsilon}}.z,V)\otimes\bigwedge_{\lambda\in
      Q}\lambda\cong\bigoplus_{Q\subset D} {\mathcal C}(I.(-\lambda_{D-Q})x_0,V)\otimes\bigwedge_{\lambda\in Q}\lambda.$$

Combining with what we saw in (1) and (2) we are done.\hfill$\Box$\\

 Written out, the natural map $C_{\alpha}{\mathcal F}_{V_B}^{\bullet}\longrightarrow
    M_{\chi}(V)$ is the $\widetilde{G}$-equivariant
    complex\begin{gather}0\longrightarrow{\rm
        ind}_{\widetilde{K}Z}^{\widetilde{G}}V_B\otimes\bigwedge^{|\nabla|}{\mathbb
        Z}[\nabla]\longrightarrow\ldots\longrightarrow{\rm
        ind}_{\widetilde{K}Z}^{\widetilde{G}}V_B\otimes\bigwedge^{0}{\mathbb
        Z}[\nabla]\longrightarrow M_{\chi}(V)\longrightarrow0.\label{indkos}\end{gather}

\begin{kor}\label{endredkom}

If the complex ${\mathcal
    K}^{\bullet}_{\nabla}(V)\otimes_RR/{\mathfrak p}$ is exact for any prime ideal ${\mathfrak p}$
  of $R$ (including ${\mathfrak p}=0$) then the $B$-module $M_{\chi}(V)$ is
  free and the Koszul complex (\ref{indkos}) is an exact resolution of
  $M_{\chi}(V)$.

\end{kor}

{\sc Proof:} Proposition \ref{marialichtmessneu} shows that the complexes $\overline{C}{\mathcal F}_{V,\epsilon}^{\bullet}$ are direct sums of complexes
  isomorphic with ${\mathcal K}^{|D|-\bullet}_D(V)$ for
  suitable $D\subset \nabla$. Now notice that formation of all the complexes
  involved commutes with base changes $R\to R/{\mathfrak p}$. Therefore Lemma
  \ref{borelred} and our hypothesis show that the complexes
  $\overline{C}{\mathcal F}_{V,\epsilon}^{\bullet}\otimes_RR/{\mathfrak p}$ are exact, thus the $B$-module $M_{\chi}(V)$ is
  free by Proposition \ref{bellotwtrick}.

For the second statement, assume first that
$B=R[T_{\lambda}]_{\lambda\in\nabla}$ and that $\chi=\theta^*:{\mathfrak
    H}_V(\widetilde{G},\widetilde{K})\to
  R[T_{\lambda}]_{\lambda\in\nabla}$ (with ${\theta}^*$ as in formula (\ref{thetastern})). Since $B$ is $R$-flat, the exactness of the complexes ${\mathcal K}^{|D|-\bullet}_D(V)$ --- here
  again we invoke Lemma
  \ref{borelred} and our hypothesis --- implies the exactness of the complexes ${\mathcal K}^{|D|-\bullet}_D(V_B)$ for
  all $D\subset \nabla$. Proposition \ref{fide} and Proposition \ref{marialichtmessneu} show that the complex $C_{\alpha}{\mathcal F}_{V_B}^{\bullet}$ admits a filtration with graded pieces
  isomorphic with complexes of the form ${\mathcal K}^{|D|-\bullet}_D(V_B)$ for
  suitable $D\subset \nabla$. Together we obtain the exactness of the complexes
  $C_{\alpha}{\mathcal F}_{V_B}^{\bullet}$. These are therefore $B$-free
  resolutions of the $B$-free module $M_{\chi}(V)$. As all the terms are
  $B$-free, tensoring with an arbitrary
  $B$-algebra does not spoil the exactness. Since a general $B$ and $\chi$
  (with presribed $\chi|_{{\mathfrak H}(V)}=\theta$) is obtained from the specific $B=R[T_{\lambda}]_{\lambda\in\nabla}$ and
  $\chi=\theta^*$ by base change, we obtain that the Koszul complex (\ref{indkos}) is an exact resolution of
  $M_{\chi}(V)$, for any $B$ and $\chi$.\hfill$\Box$\\

{\it Remark:} The complexes $\overline{C}{\mathcal
  F}_{V,\epsilon}^{\bullet}$ truly depend on $\theta:{\mathfrak
  H}(V)\to R$, but the complexes ${\mathcal K}^{\bullet}_D(V)$ visibly do
not. Moreover, the ${\mathcal K}^{\bullet}_D(V)$ only depend on the restriction to ${\mathfrak N}$ of the given
$\widetilde{K}$-action on $V$.\\

\section{Exactness criteria for ${\mathcal
    K}^{\bullet}_{\nabla}(V)$}

\begin{kor}\label{rank2} Suppose that

--- either $|\nabla|=1$,

--- or $|\nabla|=2$, we are in case (I) and the action of ${\mathfrak N}$ on $V$ (induced by the action of $\widetilde{K}$ on $V$)
is trivial.

Then the $B$-module $M_{\chi}(V)$ is free and the Koszul complex (\ref{indkos}) is a resolution of $M_{\chi}(V)$, for any $R$-algebra
  $B$ and any $\chi:{\mathfrak H}_V(\widetilde{G},\widetilde{K})\to B$ such
that $\chi|_{{\mathfrak H}(V)}$ factors through $R$.
\end{kor}

{\sc Proof:} If $|\nabla|=1$, say $\nabla=\{\lambda\}$, it is enough to show, by Corollary \ref{endredkom}, the injectivity of $$V\longrightarrow {\rm
  ind}_{{\mathfrak N}[\lambda]}^{\mathfrak N}V,\quad v\mapsto [{\mathfrak N}\to V, i\mapsto
\xi_{\lambda}(i^{-1}v)].$$[We identified ${\rm
  ind}_{\mathfrak N}^{\mathfrak N}V\cong V$, $f\mapsto f(1)$.] In case (I) (where $\xi_{\lambda}={\rm
  id}_V$) this injectivity is trivial. In case (III) it
follows from case (II). In case (II) we need to show that $\cap_{i\in {\mathfrak N}}i{\rm ker}(\xi_{\lambda})=0$. Now ${\rm ker}(\xi_{\lambda})$ is the kernel of the natural projection $V\to V_{{\bf N}_{\lambda}(k_F)}$, and from this we see that the collection of subspaces $i{\rm ker}(\xi_{\lambda})$ for ${i\in {\mathfrak N}}$ can be identified with the collection of the kernels of $|k_F|$ of the $|k_F|+1$ many functionals discussed in Lemma 21 of \cite{bali}. From that Lemma 21 of \cite{bali} --- which in a slightly different formulation had been proved earlier by Teitelbaum --- we indeed obtain $\cap_{i\in {\mathfrak N}}i{\rm ker}(\xi_{\lambda})=0$ as desired. 

If $|\nabla|=2$, say $\nabla=\{\lambda_1, \lambda_2\}$, and if we are in case (I) and ${\mathfrak N}$ acts trivially on $V$, then we need to show, by
Corollary \ref{endredkom}, the exactness of $$0\longrightarrow {\rm
  ind}_{\mathfrak N}^{\mathfrak N}V \longrightarrow    {\rm
  ind}_{{\mathfrak N}[\lambda_1]}^{\mathfrak N}V\oplus{\rm
  ind}_{{\mathfrak N}[\lambda_2]}^{\mathfrak N}V\longrightarrow {\rm
  ind}_{{\mathfrak N}[\{\lambda_1,\lambda_2\}]}^{\mathfrak N}V.$$Injectivity even of ${\rm
  ind}_{\mathfrak N}^{\mathfrak N}V \longrightarrow    {\rm
  ind}_{{\mathfrak N}[\lambda_1]}^{\mathfrak N}V$ is obvious. The exactness in ${\rm
  ind}_{{\mathfrak N}[\lambda_1]}^{\mathfrak N}V\oplus{\rm
  ind}_{{\mathfrak N}[\lambda_2]}^{\mathfrak N}V$ follows
easily from the fact that ${\mathfrak N}$ is generated by ${\mathfrak N}[\lambda_1]$ and ${\mathfrak N}[\lambda_2]$. [We have ${\mathfrak N}\cap
U_{\alpha_1}\subset {\mathfrak N}[\lambda_1]$ and $ {\mathfrak N}\cap
U_{\alpha_2}\subset{\mathfrak N}[\lambda_2]$ if $\alpha_1, \alpha_2\in\Phi^-$ denote the negatives of the two simple
roots, and ${\mathfrak N}$ is even generated by ${\mathfrak N}\cap
U_{\alpha_1}$ and ${\mathfrak N}\cap
U_{\alpha_2}$.]\hfill$\Box$\\ 

{\it Remark:} Corollary \ref{rank2} for $|\nabla|=1$ has been proved first in the paper
\cite{bali} by Barthel and Livn\'{e}. The freeness of the $B$-module
$M_{\chi}(V)$ in case (I) for $\widetilde{G}=G={\rm
  PGL}_3(F)$ is the main result of the paper \cite{bo} by Bella${\rm \ddot{\i}}$che and
Otwinowska (where only flatness is claimed, but in fact freeness is proved). However, it seems that
already Corollary \ref{rank2} (b) for $\widetilde{G}={\rm GL}_3(F)$ does
not easily follow from that result for ${\rm
  PGL}_3(F)$.\\

\begin{lem}\label{stjosef} Let ${\mathfrak G}$ be a pro-$p$-group, let$$\ldots\stackrel{d_{i+1}}{\longrightarrow}
  C_{i+1}\stackrel{d_i}{\longrightarrow}
  C_{i}\stackrel{d_{i-1}}{\longrightarrow}C_{i-1}\stackrel{d_{i-2}}{\longrightarrow}\ldots$$ be a
  complex of ${\mathbb F}_p[{\mathfrak G}]$-modules with $C_i=0$ for $i>>0$. Fix
  $i\in{\mathbb Z}$. If for all $m\ge s\ge 0$ the
  complexes$$\ldots\longrightarrow H^{s}({{\mathfrak G}},
  C_{i+m+2})\longrightarrow H^{s}({{\mathfrak G}},C_{i+m+1})\longrightarrow H^{s}({{\mathfrak G}},C_{i+m})$$are exact, then $\ldots\longrightarrow C_{i+2}\stackrel{d_{i+1}}{\longrightarrow}
  C_{i+1}\stackrel{d_i}{\longrightarrow}
  C_{i}$ is exact.
\end{lem}

{\sc Proof:} We argue by descending induction on $i$, beginning with the
maximal $i$ with $C_i\ne0$ (if $C_i=0$ for all $i$ then there is nothing to
prove). In the induction step we may assume
that $$\ldots\stackrel{d_{i+3}}{\longrightarrow}
C_{i+3}\stackrel{d_{i+2}}{\longrightarrow}
C_{i+2}\stackrel{d_{i+1}}{\longrightarrow} C_{i+1}$$is exact, and we need to
show ${\rm Ker}(d_i)={\rm Im}(d_{i+1})$. If ${\rm Ker}(d_i)/{\rm Im}(d_{i+1})\ne0$ then $H^0({\mathfrak G},{\rm Ker}(d_i)/{\rm Im}(d_{i+1}))\ne0$ since ${\mathfrak G}$
is a pro-$p$-group and ${\rm Ker}(d_i)/{\rm Im}(d_{i+1})$ is an ${\mathbb F}_p$-vector
space. Therefore we only need to show $H^0({\mathfrak G},{\rm Ker}(d_i)/{\rm
  Im}(d_{i+1}))=0$, or equivalently, the injectivity of $$d_i:H^0({\mathfrak G},C_{i+1}/{\rm Im}(d_{i+1}))\to H^0({\mathfrak G},C_i).$$By descending induction on $j$ we prove that, more generally,$$d_{i+j}:H^j({\mathfrak G},C_{i+j+1}/{\rm Im}(d_{i+j+1}))\to H^j({\mathfrak G},C_{i+j})$$is injective, for all $j\ge0$. For $j>>0$ there is nothing to show. Now assume that we know this claim for some $j>0$. By the (outer) induction hypothesis we know ${\rm Im}(d_{i+1+j})={\rm Ker}(d_{i+j})$, hence that $d_{i+j}$ induces an isomorphism $C_{i+j+1}/{\rm Im}(d_{i+j+1})\cong  {\rm Im}(d_{i+j})$. Therefore, the (inner) induction hypothesis and the long exact cohomology sequence associated with the short exact sequence $$0\longrightarrow{\rm Im}(d_{i+j})\longrightarrow C_{i+j}\longrightarrow \frac{C_{i+j}}{{\rm Im}(d_{i+j})}\longrightarrow0$$show that the map $H^{j-1}({\mathfrak G},C_{i+j})\to H^{j-1}({\mathfrak G},{C_{i+j}}/{{\rm Im}(d_{i+j})})$ is surjective. On the other hand, by hypothesis the sequence$$H^{j-1}({\mathfrak G},C_{i+j+1})\stackrel{d_{i+j}}{\longrightarrow} H^{j-1}({\mathfrak G},C_{i+j})\stackrel{d_{i+j-1}}{\longrightarrow} H^{j-1}({\mathfrak G},C_{i+j-1})$$is exact. Together it follows that $$d_{i+j-1}:H^{j-1}({\mathfrak G},C_{i+j}/{\rm Im}(d_{i+j}))\to H^{j-1}({\mathfrak G},(C_{i+j-1})$$ is injective, as desired.\hfill$\Box$\\

\begin{lem}\label{exmodp} Let $R$ be a discrete valuation ring with residue
  field $k$. A complex $$\ldots\longrightarrow C_{i+1}{\longrightarrow}
  C_{i}{\longrightarrow}C_{i-1}{\longrightarrow}\ldots$$of
  finite free $R$-modules with $C_i=0$ for $i>>0$ is exact if the reduced
  complex$$\ldots\longrightarrow C_{i+1}\otimes k{\longrightarrow}
  C_{i}\otimes k{\longrightarrow}C_{i-1}\otimes k{\longrightarrow}\ldots$$is exact.
\end{lem}

{\sc Proof:} Let $i_0$ be such that $C_i=0$ for $i>i_0$. By Nakayama's Lemma,
that $C_{i_0}\otimes k\to
C_{i_0-1}\otimes k$ is injective implies that the image of $C_{i_0}\to
C_{i_0-1}$ is $R$-free of the same rank as $C_{i_0}$ and is a
direct summand of $C_{i_0-1}$. Hence $C_{i_0}\to
C_{i_0-1}$ is injective and $C_{i_0-1}/{\rm Im}[C_{i_0}\to C_{i_0-1}]$ is
$R$-free. Replacing $C_{i_0}$ by $0$ and $C_{i_0-1}$ by $C_{i_0-1}/{\rm
  Im}[C_{i_0}\to C_{i_0-1}]$ we obtain a new complex satisfying the hypotheses
of our Lemma. We conclude by induction. \hfill$\Box$\\  

 Suppose that $R$ is a discrete valuation
  ring with residue field $k$ of characteristic $p$, or $R=k$ is itself a
  field of characteristic $p$ (in the following read $V\otimes_Rk=V$ if $R=k$). For $\lambda\in\nabla$ and $i\ge0$ let ${{{{\mathfrak N}}}}[\lambda]$ act on $V\otimes_Rk$ through the isomorphism\begin{gather}{{{{\mathfrak N}}}}[\lambda]=\gamma_{[\lambda]}^{-1}{{{{\mathfrak N}}}}\gamma_{[\lambda]}\longrightarrow
{{{{\mathfrak N}}}},\quad\quad \gamma_{[\lambda]}^{-1}i\gamma_{[\lambda]}\mapsto
i\label{jan}\end{gather}composed with the inclusion of ${\mathfrak N}$ into
$\widetilde{K}$ and the given $\widetilde{K}$-action on $V$. Consider the two maps$$\iota_{\lambda}^i,\rho_{\lambda}^i:H^i({{{{\mathfrak N}}}},V\otimes_Rk)\to
  H^i({{{{\mathfrak N}}}}[\lambda],V\otimes_Rk)$$defined as follows. The map $\rho_{\lambda}^i$ is the one induced in cohomology by the inclusion of
groups ${{{{\mathfrak N}}}}[\lambda]\subset {{{{\mathfrak N}}}}$ together with the
${{{{\mathfrak N}}}}[\lambda]$-equivariant map
$\xi_{\lambda}:V\otimes_Rk\longrightarrow V\otimes_Rk$ (cf. Lemma
\ref{petersfrage} and the remark following it). The map $\iota_{\lambda}^i$ is the isomorphism induced by conjugation with $\gamma_{[\lambda]}$.

\begin{pro}\label{bettercond} Suppose that for any $0\le i<|\nabla|$ the $k$-vector space $H^i({\mathfrak N},V\otimes_Rk)$ is generated by elements $x$ with the following property: for any $\lambda\in\nabla$ we have $\rho_{\lambda}^i(x)=\theta_{\lambda,x}\iota_{\lambda}^i(x)$ for some $\theta_{\lambda,x}\in k$, and moreover $\theta_{\lambda,x}\ne0$ for at least one $\lambda$ (depending on $x$). 

Then the complexes ${\mathcal
    K}^{\bullet}_{\nabla}(V)$ and ${\mathcal
    K}^{\bullet}_{\nabla}(V)\otimes_R k$ are exact. In particular (Corollary
  \ref{endredkom}), the $B$-module $M_{\chi}(V)$ is free and the Koszul complex (\ref{indkos}) is a resolution of $M_{\chi}(V)$, for any $R$-algebra
  $B$ and any $R$-algebra homomorphism $\chi:{\mathfrak H}_V(\widetilde{G},\widetilde{K})\to B$ such
  that $\chi|_{{\mathfrak H}(V)}$ factors through $R$.
\end{pro}

We emphasize that our hypothesis concerns only small cohomology degrees $i$, i.e. cohomology degrees $i$ strictly smaller than $|\nabla|$.\\

{\sc Proof:} Suppose first that $R$ is a discrete valuation
  ring with residue field $k$ of characteristic $p$. (Then we must be in case (I) or in case (III).) Since formation of the complex of $R$-modules ${\mathcal
  K}^{\bullet}_{\nabla}(V)$ commutes with base changes in $R$, Lemma
\ref{exmodp} tells us that we may reduce this case to the case where $R=k$, so from now on in
this proof $R$ is a field of characteristic $p$. Then we must be in case (I) or in case (II); in case (I) the action of the pro-$p$-group ${\mathfrak N}$ on the one-dimensional $k$-vector space $V$ is necessarily trivial.

{\it Step 1:} We appeal to Lemma
\ref{stjosef} which tells us that, since $R$ is a field of characteristic $p$, to show exactness of the complex ${\mathcal K}^{\bullet}_{\nabla}(V)$ it is enough to show exactness of the complexes\begin{gather}0\longrightarrow
  H^i({{{{\mathfrak N}}}},V)\longrightarrow\bigoplus_{\lambda\in
    \nabla}H^i({{{{\mathfrak N}}}},{\rm ind}_{{{{{\mathfrak N}}}}[\lambda]}^{{\mathfrak N}}V)\longrightarrow\ldots\longrightarrow \bigoplus_{E\subset
    \nabla\atop|E|=|\nabla|-i} H^i({{{{\mathfrak N}}}},{\rm ind}_{{{{{\mathfrak N}}}}[E]}^{{\mathfrak N}}V)\label{ufpex}\end{gather}for all $0\le i<|\nabla|$. Using Shapiro's
  Lemma this means proving exactness of the complexes\begin{gather}0\longrightarrow
  H^i({{{{\mathfrak N}}}},V)\longrightarrow\bigoplus_{\lambda\in
    \nabla}H^i({{{{\mathfrak N}}}}[\lambda],V)\longrightarrow\ldots\longrightarrow \bigoplus_{E\subset
    \nabla\atop|E|=|\nabla|-i} H^i({{{{\mathfrak N}}}}[E],V)\label{gernex}\end{gather}for
all $0\le i<|\nabla|$. Here ${{{{\mathfrak N}}}}[E]$ acts on $V$ through the isomorphism$${{{{\mathfrak N}}}}[E]=\gamma_{[E]}^{-1}{{{{\mathfrak N}}}}\gamma_{[E]}\longrightarrow
{{{{\mathfrak N}}}},\quad\quad \gamma_{[E]}^{-1}i\gamma_{[E]}\mapsto
i$$composed with the inclusion of ${\mathfrak N}$ into
$\widetilde{K}$ and the given $\widetilde{K}$-action on $V$. The differentials
  are the sums of the maps$$\rho_{E,E\cup\{\lambda\}}^i:H^i({{{{\mathfrak N}}}}[E],V)\to
  H^i({{{{\mathfrak N}}}}[E\cup\{\lambda\}],V)$$induced in cohomology by the inclusions of
groups ${{{{\mathfrak N}}}}[E\cup\{\lambda\}]\subset {{{{\mathfrak N}}}}[E]$ together with the
${{{{\mathfrak N}}}}[E\cup\{\lambda\}]$-equivariant maps
$\sigma(\lambda,E)\xi_{\lambda}:V\longrightarrow V$ (cf. Lemma
\ref{petersfrage} and the remark following it). We have $\rho_{\lambda}^i=\pm\rho_{\emptyset,\{\lambda\}}^i$.\\

{\it Step 2:} For any $E'\subset E\subset \nabla$ and $i\ge0$ conjugation with
$\gamma_{[E-E']}$ induces an isomorphism
$$\iota_{E',E}^i:H^i({{{{\mathfrak N}}}}[E'],V)\to H^i({{{{\mathfrak
        N}}}}[E],V).$$ We have
$\iota_{\lambda}^i=\iota_{\emptyset,\{\lambda\}}^i$. Moreover, for
$E\subset\nabla-\{\lambda\}$ we have $\rho_{E,E\cup\{\lambda\}}^i\circ\iota_{\emptyset,E}^i=\pm\iota_{\{\lambda\},E\cup\{\lambda\}}^i\circ\rho^i_{\emptyset,\{\lambda\}}$.

For an element $x\in H^i({\mathfrak N},V)$ as described in the hypothesis of Proposition \ref{bettercond} and for any $E\subset\nabla$ and any $\lambda\in\nabla-E$
we have \begin{align}\rho_{E,E\cup\{\lambda\}}^i(\iota_{\emptyset,E}^i(x))&=\pm\iota_{\{\lambda\},E\cup\{\lambda\}}^i(\rho^i_{\emptyset,\{\lambda\}}(x))\notag\\{}&=\pm\theta_{\lambda,x}\iota_{\{\lambda\},E\cup\{\lambda\}}^i(\iota^i_{\emptyset,\{\lambda\}}(x))\notag\\{}&=\pm\theta_{\lambda,x} \iota^i_{\emptyset,E\cup\{\lambda\}}(x).\label{reio}\end{align}

Now choose an $R$-basis of $H^i({{{{\mathfrak N}}}},V)$ consisting of elements $x$ as described in the hypothesis of Proposition \ref{bettercond}. Consider the
direct sum decomposition of $H^i({{{{\mathfrak N}}}},V)$ into its one dimensional subspaces $R.x$ spanned by these basis elements. By means of the isomorphisms $\iota_{\emptyset,E}^i:H^i({{{{\mathfrak N}}}},V)\cong H^i({{{{\mathfrak N}}}}[E],V)$ this decomposition induces a decomposition of $H^i({{{{\mathfrak N}}}}[E],V)$ into one dimensional subspaces $R.\iota_{\emptyset,E}^i(x)$, for each $E\subset\nabla$. By formula (\ref{reio}) these decompositions in fact define a direct sum decomposition of the complex (\ref{gernex}). For each $x$ in our chosen $R$-basis of $H^i({{{{\mathfrak N}}}},V)$ there is some $\lambda\in\nabla$ such that $\theta_{\lambda,x}$ is a unit in $k$. By formula (\ref{reio}) this means that the maps $\rho_{E,E\cup\{\lambda\}}^i:R.\iota_{\emptyset,E}^i(x)\to R.\iota_{\emptyset,E\cup\{\lambda\}}^i(x)$ are isomorphisms for any
$E\subset\nabla-\{\lambda\}$ (for this fixed $\lambda$ depending on $x$). From this we easily deduce that the direct summands of the complex (\ref{gernex}) spanned by the $\iota_{\emptyset,E}^i(x)$ are exact. Thus the complex (\ref{gernex}) must be exact, and we are done. 

[In fact this argument shows that not just the complexes (\ref{gernex}) for
$0\le i<|\nabla|$ are exact, but even the continued complexes up to the last
term $H^i({{{{\mathfrak N}}}}[\nabla],V)$.]
\hfill$\Box$\\ 

\section{Group Cohomology of integral unipotent radicals}

\label{cohin}

In this section, which is independent of the preceding ones, we transfer a result of Dwyer, Friedlander/Parshall and Polo/Tilouine on Lie algebra cohomology to a result on the cohomology of the group of ${\mathbb Z}_p$-valued points of the unipotent radical of a Borel subgroup in a split semisimple algebraic group ${\mathcal G}$ over ${\mathbb Z}_p$.

Let ${\mathcal T}$ denote a split maximal torus in ${\mathcal G}$, let $\Phi=\Phi^+\coprod\Phi^-$ denote the choice of a system of positive resp. negative roots in the corresponding root system $\Phi$. Let $\ell:{\mathcal W}\to{\mathbb Z}_{\ge0}$ denote the corresponding length function on the Weyl group ${\mathcal W}={\mathcal W}({\mathcal T})$. For $\alpha\in\Phi$ let ${\mathcal U}_{\alpha}\subset{\mathcal G}$ denote the corresponding root subgroup. Let $${\mathcal N}({\mathbb Z}_p)=\prod_{\alpha\in\Phi^-} {\mathcal U}_{\alpha}({\mathbb Z}_p)$$and for $w\in {\mathcal W}$ let$${\mathcal N}_w({\mathbb Z}_p)=\prod_{\alpha\in\Phi^-\cap w{\Phi^+}}{\mathcal U}_{\alpha}({\mathbb Z}_p),$$$${\mathcal N}({\mathbb Z}_p)(w)=\prod_{\alpha\in w{\Phi^+}}{\mathcal U}_{\alpha}({\mathbb Z}_p).$$Thus ${\mathcal N}({\mathbb Z}_p)$, resp. ${\mathcal N}_w({\mathbb Z}_p)$, is the group of ${\mathbb Z}_p$-valued points of the unipotent radical ${\mathcal N}$ of a Borel subgroup in ${\mathcal G}$, resp. of an algebraic subgroup ${\mathcal N}_w$ of ${\mathcal N}$, and ${\mathcal N}({\mathbb Z}_p)(w)$ is a subgroup of ${\mathcal G}({\mathbb Z}_p)$ obtained from ${\mathcal N}({\mathbb Z}_p)$ by conjugation with an element in ${\mathcal W}$. The dimension of the ${\mathbb Q}_p$-analytic group ${\mathcal N}_w({\mathbb Z}_p)$ is $\ell(w)$. 

The set of $p$-restricted weights is$$X_1^*({\mathcal T})=\{\lambda\in
X^*({\mathcal T})={\rm Hom}({\mathcal T},{\mathbb
  G}_m)\,|\,0\le\langle\lambda,\check{\alpha}\rangle<p\mbox{ for all simple
  roots }\alpha\}.$$For any $\mu\in X_1^*({\mathcal T})$ there is a unique
irreducible rational ${\mathcal G}_{{\mathbb F}_p}$-representation $V(\mu)$
over $\overline{\mathbb F}_p$ such that the space $V(\mu)_{{\mathcal
    N}({\mathbb Z}_p)}$ of coinvariants for the action of ${\mathcal
  N}({\mathbb Z}_p)$ on $V(\mu)$ is one dimensional and such that
${\mathcal T}({\mathbb F}_p)$ acts on $V(\mu)_{{\mathcal N}({\mathbb Z}_p)}$ through the character $\mu$. (We regard $V(\mu)$ as a representation of ${\mathcal N}({\mathbb
  Z}_p)$ via the natural surjection ${\mathcal N}({\mathbb
  Z}_p)\to{\mathcal N}({\mathbb F}_p)$.)  In \cite{herwe}, Appendix,
Proposition 1.3 it is shown that this $V(\mu)$ is irreducible even as
a representation of the abstract group ${\mathcal G}({{\mathbb F}_p})$, and
that, if
${\mathcal G}_{{\mathbb F}_p}$ has simply connected derived group, then each
irreducible $\overline{\mathbb F}_p[{\mathcal G}({{\mathbb F}_p})]$-module
arises in this way.

\begin{satz}\label{poloclou} Assume that the algebraic group ${\mathcal G}$ is simple. Let $\rho$ denote the half sum of the positive roots. Suppose that $\mu\in X_1^*({\mathcal T})$ satisfies $\langle\mu+\rho,\check{\beta}\rangle\le p$ for all $\beta\in\Phi^+$. Then the natural map\begin{gather}H^i({\mathcal N}({\mathbb Z}_p),V(\mu))\longrightarrow\prod_{w\in {\mathcal W}\atop \ell(w)=i}H^i({\mathcal N}_w({\mathbb Z}_p),V(\mu)_{{\mathcal N}({\mathbb Z}_p)(w)})\label{polotilgal}\end{gather}is bijective, for any $i\in{\mathbb Z}$.
\end{satz}

The map in (\ref{polotilgal}) is the product of the natural maps induced by the projections of modules $V(\mu)\to V(\mu)_{{\mathcal N}({\mathbb Z}_p)(w)}$ and the inclusion of groups ${\mathcal N}_w({\mathbb Z}_p)\to{\mathcal N}({\mathbb Z}_p)$. (The action of ${\mathcal N}_w({\mathbb Z}_p)$ on $V(\mu)_{{\mathcal N}({\mathbb Z}_p)(w)}$ is trivial, of course.)

Let $\alpha_0$ denote the highest short root. The conditions on $\mu$ imply
$p\ge h-1$, where $h=\langle\rho,\check{\alpha}_0\rangle+1$ is the Coxeter
number of the root system $\Phi$. Indeed, we
have\begin{gather}p\ge\langle\mu+\rho,\check{\alpha}_0\rangle=h-1+\langle\mu,\check{\alpha}_0\rangle\ge
h-1\label{pcox}\end{gather}(cf. \cite{vigre} p. 14). Using the terminology of alcoves the conditions on $\mu$ are equivalent with
saying that $\mu$ belongs to the '(closed) bottom alcove' (with respect to $p$) in the cone of
dominant weights. (Observe that it corresponds to (the intersection
with $X^*({\mathcal T})$ of) the {\it closure} of the region referred to as the 'lowest
alcove' e.g. in \cite{herwe}.) \\

We are going to prove Theorem \ref{poloclou} via a comparison with the analogous statement in Lie algebra cohomology. Let ${\mathfrak n}$ denote the Lie algebra of the algebraic group ${\mathcal N}_{{\mathbb F}_p}$ over ${\mathbb F}_p$. More generally, for $w\in {\mathcal W}$ let ${\mathfrak n}_{w}$ denote the Lie algebra of ${\mathcal N}_{w,{\mathbb F}_p}$.     

\begin{satz}\label{liepoloclou} (Dwyer, Friedlander/Parshall, Polo/Tilouine) Assume that ${\mathcal G}$ is simple and that $\mu\in X_1^*({\mathcal T})$ satisfies $\langle\mu+\rho,\check{\beta}\rangle\le p$ for all $\beta\in\Phi^+$. Then the natural map\begin{gather}H^i({{{{\mathfrak n}}}},V(\mu))\longrightarrow\prod_{w\in {\mathcal W}\atop \ell(w)=i}H^i({\mathfrak n}_w,V(\mu)_{{\mathcal N}({\mathbb Z}_p)(w)})\label{polotillie}\end{gather}is bijective, for any $i\in{\mathbb Z}$.
\end{satz}

{\sc Proof:} Recall the dot action of ${\mathcal W}$ on $X^*({\mathcal T})$
defined by $w\cdot\eta=w(\eta-\rho)+\rho$ (where $w\in {\mathcal W}$, $\eta\in
X^*({\mathcal T})$). For $\eta\in X^*({\mathcal T})$ let $L(\eta)$ denote the
one dimensional ${\mathcal T}_{{\mathbb F}_p}$-representation over
$\overline{\mathbb F}_p$ which $\eta$ defines. Let $w_0\in {\mathcal W}$
denote the longest element. In \cite{pt} it is shown that we have an
isomorphism of ${\mathcal T}$-modules\begin{gather}H^i({{{{\mathfrak
          n}}}},V(\mu))\cong \prod_{w\in {\mathcal W}\atop \ell(w)=i} L(w\cdot
  (w_0\mu)).\label{lievig}\end{gather}In \cite{friedpar} Theorem 2.5 this
result had been proved under the slightly stronger hypothesis
$\langle\mu+\rho,\check{\beta}\rangle\le p-1$ for all $\beta\in\Phi^+$. In
\cite{dwy} it had been proved for ${\mathcal G}={\rm SL}_n$ and $\mu=0$. In
\cite{vigre} Theorem 4.2.1 the isomorphism (\ref{lievig}) is reproved (in
general), and it is made explicit. Namely, the following is shown. In the
standard chain complex $\wedge^{\bullet}({{{\mathfrak n}}}^*)\otimes V(\mu)$
computing the Lie algebra cohomology $H^i({{{{\mathfrak n}}}},V(\mu))$ the
weight space of $w\cdot(w_0\mu)$ is one dimensional, spanned by a
vector \begin{gather}f_{\alpha_1}\wedge\ldots\wedge f_{\alpha_i}\otimes
  v_{w(w_0\mu)}\label{spanweig}\end{gather} where $\Phi^-\cap
w\Phi^+=\{\alpha_1,\ldots,\alpha_i\}$, where
$\{f_{\alpha}\,|\,\alpha\in\Phi^-\}$ is a basis of ${{{\mathfrak n}}}^*$
consisting of weight vectors, and where $v_{w(w_0\mu)}\in V(\mu)$ is a non zero vector of weight $w(w_0\mu)$ (no dot here). Now ${\mathcal T}$ acts on $V(\mu)_{{\mathcal N}({\mathbb Z}_p)(w)}$ through $w(w_0\mu)$ and we have a ${\mathcal T}$-equivariant surjection$$\wedge^{\bullet}({{{\mathfrak n}}}^*)\otimes V(\mu)\longrightarrow\wedge^{\bullet}({{{\mathfrak n}}}_w^*)\otimes V(\mu)_{{\mathcal N}({\mathbb Z}_p)(w)}.$$Hence also the weight space of $w\cdot(w_0\mu)$ in the standard chain complex $\wedge^{\bullet}({{{\mathfrak n}}}_w^*)\otimes V(\mu)_{{\mathcal N}({\mathbb Z}_p)(w)}$ computing the Lie algebra cohomology $H^i({\mathfrak n}_w,V(\mu)_{{\mathcal N}({\mathbb Z}_p)(w)})$ is one dimensional, spanned by the image of the vector (\ref{spanweig}). On the other hand, as ${\rm dim}({\mathfrak n}_w)=i$ we know that $H^i({\mathfrak n}_w,V(\mu)_{{\mathcal N}({\mathbb Z}_p)(w)})$ is isomorphic with the ${\mathfrak n}_w$-coinvariants of the one dimensional (trivial) ${\mathfrak n}_w$-representation $V(\mu)_{{\mathcal N}({\mathbb Z}_p)(w)}$. Thus $L(w\cdot(w_0\mu))\cong H^i({\mathfrak n}_w,V(\mu)_{{\mathcal N}({\mathbb Z}_p)(w)})$ and we are done.  \hfill$\Box$\\ 

For $\lambda\in X^*({\mathcal T})$ let $V(\mu)_{\lambda}$ denote the corresponding
weight space in $V(\mu)$ for the action of ${\mathcal T}$. Let $\Pi(\mu)=\{\lambda\in X^*({\mathcal T})\,|\,V(\mu)_{\lambda}\ne0\}$.

\begin{lem}\label{sana} Let ${\mathcal G}$ and $\mu\in X_1^*({\mathcal T})$ be as in Theorem \ref{liepoloclou}. For any $\lambda\in \Pi(\mu)$, any $\alpha\in\Delta$ and any $s\ge0$
such that also $\lambda-s\alpha\in \Pi(\mu)$ we have $s<p$.
\end{lem}

{\sc Proof:} By the strong linkage principle
\cite{jant}, the condition $\langle\mu+\rho,\check{\beta}\rangle\le
p$ (all $\beta\in\Phi^{+}$) on the $p$-restricted weight $\mu$ implies that
$\Pi(\mu)$ coincides with the corresponding set defined with respect to the
irreducible rational ${\mathcal G}_{\mathbb C}$-module of highest weight
$\mu$. Therefore it follows from \cite{bou} chapter VIII, section 7.2, Proposition 3
that for any $\lambda\in \Pi(\mu)$ and any $s\ge0$
such that also $\lambda-s\alpha\in \Pi(\mu)$ we have $s\le\langle
\lambda,\check{\alpha}\rangle$. Thus we need to show $\langle
\lambda,\check{\alpha}\rangle<p$ for any $\lambda\in \Pi(\mu)$, any $\alpha\in \Delta$. Let $$X^{\rm dom}=\{\lambda\in X^*({\mathcal T})\otimes{\mathbb R}\,\,|\,\,0\le\langle
\lambda,\check{\alpha}\rangle\le p\mbox{ for all }\alpha\in\Delta\}$$denote
the (closed) dominant cone. It contains $\rho$ in its interior. Let$$\overline{C}=\{\lambda\in X^{\rm
  dom}\,|\,\langle\lambda+\rho,\check{\alpha_0}\rangle\le p\}$$(recall that $\alpha_0$ denote the highest short root), the closure in $X^*({\mathcal
  T})\otimes{\mathbb R}$ of the bottom alcove in $X^*({\mathcal
  T})\otimes{\mathbb R}$. We have$$\overline{C}=\{\lambda\in X^{\rm
  dom}\,|\,\langle\lambda+\rho,\check{\beta}\rangle\le p\mbox{ for all
}\beta\in\Phi^+\}$$and this implies\begin{gather}\langle\lambda,\check{\alpha}\rangle< p\quad\mbox{ for any
  }\lambda\in\overline{C},\mbox{ any }\alpha\in\Phi.\label{invabsc}\end{gather}For any $\lambda\in
\mu-{\mathbb Z}_{\ge0}\Phi^+$ we have
$\langle\lambda+\rho,\check{\alpha_0}\rangle\le\langle\mu+\rho,\check{\alpha_0}\rangle$. As $\Pi(\mu)\subset
\mu-{\mathbb Z}_{\ge0}\Phi^+$ and
$\mu\in\overline{C}$ it follows that $\Pi(\mu)\cap X^{\rm dom}\subset\overline{C}$. Now it is known that
${\mathcal W}(\Pi(\mu)\cap X^{\rm dom})=\Pi(\mu)$. It therefore follows from
formula (\ref{invabsc}) that $\langle
\lambda,\check{\alpha}\rangle<p$ for any $\lambda\in \Pi(\mu)$ and $\alpha\in \Delta$, as desired. \hfill$\Box$\\

Let $\widetilde{\mathfrak a}$ denote the augmentation ideal in the group ring
${\mathbb F}_p[{\mathcal N}({\mathbb Z})]$, i.e. $\widetilde{\mathfrak a}$ is the kernel
of the ring homomorphism ${\mathbb F}_p[{\mathcal N}({\mathbb Z})]\to{\mathbb F}_p$ sending $n\in{\mathcal N}({\mathbb Z})$ to
$1\in{\mathbb F}_p$. Then$$B=\lim_{\leftarrow\atop i}{\mathbb F}_p[{\mathcal
  N}({\mathbb Z})]/\widetilde{\mathfrak a}^i$$is a complete Noetherian ring,
${\mathfrak a}=\widetilde{\mathfrak a}B$ is a two sided ideal in $B$, and
$B/{\mathfrak a}={\mathbb F}_p$. If we
filter $B$ by putting $F^iB={\mathfrak a}^i$ for $i\ge 0$, then also the associated
graded ring ${\bf gr} B$ is Noetherian.

We may regard $V(\mu)$ in a natural way as a $B$-module. If for $i\ge 0$ we filter $V(\mu)$ by
putting $F^iV(\mu)={\mathfrak a}^iV(\mu)$ then $V(\mu)$ becomes a filtered
$B$-module, and its associated graded ${\bf gr} V(\mu)$ becomes a ${\bf gr}
B$ module.

Let $U({\mathfrak n})$ denote the (unrestricted) universal enveloping algebra of ${\mathfrak n}$.

Recall that the prime number $p$ is called excellent with respect to the root
system $\Phi$ if $p$ is at least the length of all root strings. Concretely,
any $p>3$ is excellent; if $\Phi$ contains no subroot system of type $G_2$
then any $p>2$ is excellent; if $\Phi$ contains no subroot system of type $G_2$, $B_n$, $C_n$ or $F_4$ (i.e. if $\Phi$ has only one root length) then any $p$ is excellent. 

\begin{pro}\label{quilaz} Let ${\mathcal G}$ and $\mu\in X_1^*({\mathcal T})$ be as in Theorem \ref{liepoloclou} and suppose that $p$ is excellent with respect to $\Phi$. 

(a) (Quillen, Lazard, Friedlander/Parshall) We have a natural isomorphism of ${\mathbb F}_p$-algebras
  \begin{gather}U({\mathfrak n})\cong {\bf gr} B.\label{neuqui}\end{gather}(b) We have a natural isomorphism of $\overline{\mathbb F}_p$-vector spaces \begin{gather}V(\mu)\cong {\bf gr} V(\mu).\label{neulika}\end{gather}The isomorphism (\ref{neulika}) is
  compatible with the respective actions of $U({\mathfrak n})$ and ${\bf gr}
  B$, if $U({\mathfrak n})$ and ${\bf gr}
  B$ are identified via the isomorphism (\ref{neuqui}).
\end{pro} 

{\sc Proof:} (a) Consider the graded ${\mathbb Z}$-Lie algebra $${\bf gr} {\mathcal N}({\mathbb Z})=\coprod_{n\ge0}{\mathcal N}({\mathbb Z})^n/{\mathcal N}({\mathbb Z})^{n+1}$$where ${\mathcal N}({\mathbb Z})^0={\mathcal N}({\mathbb Z})$ and ${\mathcal N}({\mathbb Z})^{n+1}=({\mathcal N}({\mathbb Z})^{n},{\mathcal N}({\mathbb Z}))$ and where the Lie bracket is induced by commutation in ${\mathcal N}({\mathbb Z})$. Let $U({\bf gr} {\mathcal N}({\mathbb Z}))$ denote its universal enveloping algebra. It follows from results of Quillen (cf. the proof of \cite{friedpar} Theorem 4.4; compare also with \cite{ln} Theorem 2.3) that there is a natural map $U({\bf gr} {\mathcal N}({\mathbb Z}))\to {\bf gr}B$ sending (the class of) $x\in {\mathcal N}({\mathbb Z})/({\mathcal N}({\mathbb Z}),{\mathcal N}({\mathbb Z}))$ to (the class of) $x-1\in {\mathfrak a}/{\mathfrak a}^2$ and inducing an isomorphism $$U({\bf gr} {\mathcal N}({\mathbb Z}))\otimes{\mathbb F}_p\cong {\bf gr}B.$$Next, we have the isomorphism\begin{align}U({\bf gr} {\mathcal N}({\mathbb
  Z}))\otimes{\mathbb F}_p&\cong U({\bf gr} {\mathcal N}({\mathbb
  Z})\otimes{\mathbb F}_p)\notag\\{}&\cong U({\mathfrak n})\notag\end{align}induced by the natural isomorphism
${\bf gr} {\mathcal N}({\mathbb Z})\otimes{\mathbb F}_p\cong {\mathfrak n}$, cf. e.g. \cite{friedpar} Proposition 4.2 (it is here where we need the excellence of $p$). Together we obtain the isomorphism (\ref{neuqui}).

(b) (Compare also with \cite{ln} Proposition 2.4.) Let $\Delta\subset\Phi^+$
denote the set of simple roots. We extend the height function ${\rm ht}$ to
the set ${\mathbb Z}_{\ge0}\Phi^+$ by putting $${\rm
  ht}(\sum_{\alpha\in\Delta}c_{\alpha}\alpha)=\sum_{\alpha\in\Delta}c_{\alpha}.$$ We first claim that
for all $i\ge0$ we have \begin{gather}F^iV(\mu)=\bigoplus_{\lambda\in{\mathbb
      Z}_{\ge0}\Phi^+\atop{{\rm ht}(\lambda)\ge
      i}}V(\mu)_{\mu-\lambda}.\label{peterfilt}\end{gather} We argue by
induction on $i$. For $i=0$ this follows from $V(\mu)=\sum_{\lambda\in\Pi(\mu)}V(\mu)_{\lambda}$ and $\Pi(\mu)\subset\mu-{\mathbb
  Z}_{\ge0}\Phi^+$. For the induction step we need to introduce some more
notations. For $\alpha\in \Delta$ let $x_{-\alpha}:{\mathbb A}\to{\mathcal N}$
denote the root subgroup of $-\alpha$, and let $e_{-\alpha}$ denote the
corresponding basis vector in ${\mathfrak n}$. The action of $e_{-\alpha}$ on
$V(\mu)$ satisfies $e_{-\alpha}(V(\mu)_{\lambda})\subset
V(\mu)_{\lambda-\alpha}$ for all $\lambda\in X^*({\mathcal T})$. It is shown
in \cite{jant} I, 7.8 and II, 1.12 that for any $s\ge0$ there is an element
$\frac{e_{-\alpha}^s}{s!}$ in the distribution algebra of the root subgroup of
$-\alpha$, defined over ${\mathbb Z}$, such that $s!\frac{e_{-\alpha}^s}{s!}$
acts as the $s$-fold iterate $e_{-\alpha}^s$ of ${e_{-\alpha}}$.  By
\cite{jant} II, 1.19(6), for $t\in{\mathbb F}_p$ the element $x_{-\alpha}(t)$
acts on $V(\mu)$ as \begin{gather}x_{-\alpha}(t)=1+t
e_{-\alpha}+\sum_{s>1}t^s\frac{e_{-\alpha}^s}{s!}.\label{jafoful}\end{gather}Together with Lemma \ref{sana} it follows that the operator $\frac{e_{-\alpha}^s}{s!}$, which shifts
weights by $-s\alpha$, acts trivially on $V(\mu)$ for each $s\ge p$. Therefore
formula (\ref{jafoful}) now reads\begin{gather}x_{-\alpha}(t)=1+t e_{-\alpha}+\sum_{s=2}^{p-1}t^s\frac{e_{-\alpha}^s}{s!}.\label{jafo}\end{gather}The ideal ${\mathfrak a}$ in $B$ is topologically generated by the set of all $x_{-\alpha}(t)-1$ with $\alpha\in \Delta$ and $t\in{\mathbb F}_p$. 

Now assume that we know formula (\ref{peterfilt}) with $i-1$ instead of
$i$. Formula (\ref{jafo}) easily gives
$F^iV(\mu)\subset\bigoplus_{\lambda\in{\mathbb Z}_{\ge0}\Phi^+\atop{{\rm
      ht}(\lambda)\ge i}}V(\mu)_{\mu-\lambda}$. To see the reverse inclusion
we begin by recalling that $V(\mu)$ is generated as an $\overline{\mathbb
  F}_p[{\mathcal N}({\mathbb F}_p)]$-module by its subspace $V(\mu)_{\mu}$, see \cite{hum} Proposition 2.15. Since the $x_{-\alpha}(t)$ generate ${\mathcal N}({\mathbb F}_p)$ this means that $V(\mu)$ coincides with its minimal $\overline{\mathbb F}_p$-sub vector space containing $V(\mu)_{\mu}$ and stable under the $x_{-\alpha}(t)$ for all $\alpha\in \Delta$. From formula (\ref{jafo}) and more specifically from the fact that all $s!$ in that formula are invertible in ${\mathbb F}_p$ it then follows that $V(\mu)$ coincides with its minimal $\overline{\mathbb F}_p$-sub vector space containing $V(\mu)_{\mu}$ and stable under the $e_{-\alpha}$ for all $\alpha\in \Delta$. It follows that for all $\lambda\in  {\mathbb Z}_{\ge0}\Phi^+$ with $h={\rm ht}(\lambda)$ we have $$V(\mu)_{\mu-\lambda}=\sum_{(\alpha_1,\ldots,\alpha_h)\in\Delta^h\atop\lambda=\sum_{j=1}^{h}\alpha_j}e_{-\alpha_1}\cdots e_{-\alpha_h} V(\mu)_{\mu}.$$Pick such a sequence $(\alpha_1,\ldots,\alpha_h)\in\Delta^h$ with $h\ge i$ and put \begin{align}M_0&= e_{-\alpha_2}\cdots e_{-\alpha_h} V(\mu)_{\mu},\notag\\M_-&=\sum_{m\ge0}e_{-\alpha_1}^mM_0.\notag\end{align}By induction hypothesis we have $M_-\subset F^{i-1}V(\mu)$. Therefore is is enough to show that $e_{-\alpha_1}M_0\subset (x_{-\alpha_1}(1)-1)M_-'$. By descending induction on $m$ we claim that, more generally, we have $e_{-\alpha_1}^mM_0\subset (x_{-\alpha_1}(1)-1)^mM_-$. Again the easy proof is based on formula (\ref{jafo}) and more specifically on the fact that all $s!$ in that formula are invertible in ${\mathbb F}_p$. Thus we can begin the descending induction with $m=p-1$ so that $m!$ is invertible in ${\mathbb F}_p$ for all $m$ under consideration. This concludes the proof of formula (\ref{peterfilt}). Together with formula (\ref{jafo}) we easily deduce all our claims in (b).\hfill$\Box$\\ 

{\sc Proof of Theorem \ref{poloclou}:} We first notice that if $p=2$ it
follows from (\ref{pcox}) that either $\Phi$ is of type $A_1$ and $V(\mu)$ is one- or two dimensional,
or $\Phi$ is of type $A_2$ and $V(\mu)$ is one dimensional. 

If $\Phi$ is of type $A_1$ then Theorem \ref{poloclou} is a statement for the
cohomology in degrees $i=0$ and $i=1$. For $i=0$ it is the statement that the
natural map $$H^0({\mathcal N}({\mathbb Z}_p),V(\mu))\longrightarrow
H^0({\mathcal N}({\mathbb Z}_p)_1,V(\mu)_{{\mathcal N}({\mathbb Z}_p)(1)})$$
is an isomorphism, where $1\in W$ is the trivial element. As $H^0({\mathcal
  N}({\mathbb Z}_p),.)$ is the functor of taking invariants under ${\mathcal
  N}({\mathbb Z}_p)$, and as ${\mathcal N}({\mathbb Z}_p)$ and ${\mathcal
  N}({\mathbb Z}_p)(1)$ are the unipotent radicals of Borel subgroups opposite
to each other, and as ${\mathcal N}({\mathbb Z}_p)_1$ is the trivial group,
the statement is clear. For $i=1$ the statement is that the natural
map $$H^1({\mathcal N}({\mathbb Z}_p),V(\mu))\longrightarrow H^1({\mathcal
  N}({\mathbb Z}_p)_{w_0},V(\mu)_{{\mathcal N}({\mathbb Z}_p)(w_0)})$$ is an
isomorphism, where $w_0\in W$ is the nontrivial element. Since ${\mathcal
  N}({\mathbb Z}_p)\cong{\mathbb Z}_p$ we know that $H^1({\mathcal N}({\mathbb
  Z}_p),.)$ is isomorphic with the functor of taking coinvariants under
${\mathcal N}({\mathbb Z}_p)$. As ${\mathcal N}({\mathbb Z}_p)={\mathcal
  N}({\mathbb Z}_p)(w_0)$ the statement is clear also in this case. If $\Phi$
is of type $A_2$ and $V(\mu)$ is one dimensional, we may just as well assume
that $V(\mu)$ is trivial, $V(\mu)=k$. In this case ${\mathcal N}({\mathbb Z}_p)$ is isomorphic with the group of
upper triangular unipotent $3\times3$ matrices with entries in ${\mathbb
  Z}_p$. Therefore our claim follows from the explicit
computation given in the appendix section \ref{appco}.


If $p=3$ then it follows
from (\ref{pcox}) that $\Phi$ is of not of type $G_2$. Therefore, we may now
assume that $p$ is excellent with respect to $\Phi$. 

We are going to identify the map (\ref{polotillie}), which by
Theorem \ref{liepoloclou} is an isomorphism, with a graded version of the map
(\ref{polotilgal}) in question. For $\mu=0$ this had been done in \cite{friedpar}. 

Let us first look at the respective targets. For $w\in {\mathcal W}$ with $\ell(w)=i$ the
${\mathbb Q}_p$-analytic group ${\mathcal N}_w({\mathbb Z}_p)$ is torsion free
and of dimension $i$, hence $H^i({\mathcal N}_w({\mathbb Z}_p),-)$ is isomorphic with the functor of taking coinvariants in a topological ${\mathcal
  N}_w({\mathbb Z}_p)$-module. On the other hand, as the unipotent algebraic
group ${\mathcal N}_{w,{\mathbb F}_p}$ has dimension $i$ the functor
$H^i({\mathfrak n}_w,-)$ is isomorphic with that of taking coinvariants under ${\mathfrak
  n}_w$. The trivial one-dimensional ${\mathcal N}_w$-representation $V(\mu)_{{\mathcal
    N}({\mathbb Z}_p)(w)}$ coincides with its coinvariants under ${\mathcal N}_w({\mathbb
  Z}_p)$ as well as with its coinvariants under ${\mathfrak
  n}_w$. We thus obtain an isomorphism\begin{gather}\prod_{w\in {\mathcal W}\atop \ell(w)=i}H^i({\mathfrak
    n}_w,V(\mu)_{{\mathcal N}({\mathbb
      Z}_p)(w)})\cong \prod_{w\in {\mathcal W}\atop
    \ell(w)=i}H^i({\mathcal N}_w({\mathbb Z}_p),V(\mu)_{{\mathcal N}({\mathbb
      Z}_p)(w)}).\label{easypolti}\end{gather}

 It is shown in \cite{gruen} Theorem 3.3' (ii) that the filtrations on $B$ and $V(\mu)$ induce a natural filtration of ${\rm
  Tor}_*^{ B}({\mathbb F}_p,V(\mu))$ such that for the associated graded ${\bf gr} {\rm
  Tor}_*^{ B}({\mathbb F}_p,V(\mu))$ we have a natural isomorphism$${\rm Tor}_*^{{\bf gr} B}({\mathbb F}_p,{\bf gr} V(\mu))\cong{\bf gr} {\rm
  Tor}_*^{ B}({\mathbb F}_p,V(\mu)).$$ Dually we obtain an isomorphism\begin{gather}{\rm Ext}^{*}_{{\bf gr} B}({\mathbb F}_p,{\bf gr} V(\mu))\cong {\bf gr}{\rm
  Ext}^{*}_{ B}({\mathbb F}_p,V(\mu)).\label{luzi}\end{gather}To identify the right hand side of formula (\ref{luzi}) observe that\begin{align}{\rm
  Ext}^{*}_{ B}({\mathbb F}_p,V(\mu)) &\cong H^{*}({\mathcal N}({\mathbb Z}),V(\mu))\notag \\
{} & \cong H^{*}({\mathcal N}({\mathbb Z}_p),V(\mu)).\notag\end{align}Here we used \cite{sergal} chapter I, par. 2.6, Exercise 2) (d) for the second isomorphism: we are comparing the cohomology of a discrete group which is a succesive extensions of free finitely generated groups with the continuous cohomology of its $p$-adic completion. 

[The argument can be carried out, for example, as follows. Consider a filtration of ${\mathcal N}({\mathbb Z})$ with subquotients
  isomorphic with ${\mathbb Z}$, and its parallel filtration of ${\mathcal N}({\mathbb Z}_p)$ with subquotients
 isomorphic with ${\mathbb Z}_p$. Look at the
  Hochschild-Serre spectral sequence and use that on the
  modules in question we have: both $H^1({\mathbb Z},.)$
  and $H^1({\mathbb Z}_p,.)$ are the functor of taking coinvariants (look
  e.g. at an explicit resolution of ${\mathbb Z}$), both $H^0({\mathbb Z},.)$
  and $H^0({\mathbb Z}_p,.)$ are the functor of taking invariants, and both $H^t({\mathbb Z},.)$
  and $H^t({\mathbb Z}_p,.)$ vanish for $t\ne 0,1$.]

 On the other hand, by Proposition \ref{quilaz} the left hand side of (\ref{luzi}) can be identified as \begin{align}{\rm Ext}^{*}_{{\bf gr}
  B}({\mathbb F}_p,{\bf gr} V(\mu)) &\cong {\rm Ext}^{*}_{U({\mathfrak n})}({\mathbb F}_p,V(\mu))\notag \\
{} & \cong H^{*}({\mathfrak n},V(\mu)).\notag\end{align}Therefore the isomorphim (\ref{luzi}) now reads\begin{gather}H^{*}({\mathfrak
    n},V(\mu))\cong{\bf gr} H^{*}({\mathcal N}({\mathbb
    Z}_p),V(\mu)).\label{altneuss}\end{gather}This is an isomorphism between
the source of the map (\ref{polotillie}) and a
graded version of the source of the map
(\ref{polotilgal}). It is clear from its construction that we may regard the
right hand side of the isomorphism (\ref{easypolti}) as being the graded version of the target of the map
(\ref{polotilgal}), in such a way that the isomorphisms (\ref{altneuss}) and
(\ref{easypolti}) are compatible. Our Theorem is proven.\hfill$\Box$\\

We now generalize Theorem \ref{poloclou} to the case where ${\mathcal G}$ is not necessarily simple but is a product ${\mathcal G}=\prod_{j\in J}{\mathcal G}^{(j)}$ with split simple factors ${\mathcal G}^{(j)}$. Accordingly we then have decompositions $${\mathcal T}=\prod_{j\in J}{\mathcal T}^{(j)},\quad\quad \Phi=\coprod_{j\in J}\Phi^{(j)},\quad \quad\Phi^+=\coprod_{j\in J}\Phi^{+,{(j)}},\quad\quad {\mathcal W}=\prod_{j\in j}{\mathcal W}^{(j)}.$$Let $\rho^{(j)}$ denote the half sum of the positive roots in $\Phi^{(j)}$. Any $\mu\in X_1^*({\mathcal T})$ is of the form $\mu=\prod\mu^{(j)}$ with $\mu^{(j)}\in X_1^*({\mathcal T}^{(j)})$.

\begin{satz}\label{prodpoloclou} Suppose that $\mu=\prod_{j\in J}\mu^{(j)}\in X_1^*({\mathcal T})$ satisfies $\langle\mu^{(j)}+\rho^{(j)},\check{\beta}\rangle\le p$ for all $\beta\in\Phi^{+,{(j)}}$. Then the natural map\begin{gather}H^i({\mathcal N}({\mathbb Z}_p),V(\mu))\longrightarrow\prod_{w\in {\mathcal W}\atop \ell(w)=i}H^i({\mathcal N}_w({\mathbb Z}_p),V(\mu)_{{\mathcal N}({\mathbb Z}_p)(w)})\label{prodpolotilgal}\end{gather}is bijective.
\end{satz}

{\sc Proof:} We may write $V(\mu)=\otimes_{j\in J}V(\mu^{(j)})$ as a rational
${\mathcal G}_{{\mathbb F}_p}=\prod_{j\in J}{\mathcal G}_{{\mathbb
    F}_p}^{(j)}$-representation. Consider the decomposition ${\mathcal N}=\prod_{j\in J}{\mathcal N}^{(j)}$ with ${\mathcal N}^{(j)}={\mathcal N}\cap {\mathcal G}^{(j)}$. More generally, for $w=\prod_jw^{(j)}\in {\mathcal W}=\prod_{j\in J}{\mathcal W}^{(j)}$ consider the decompositions $${\mathcal N}_w=\prod_{j\in J}{\mathcal N}_w\cap {\mathcal N}^{(j)},\quad\quad {\mathcal N}({\mathbb Z}_p)(w)\cong\prod_{j\in J}{\mathcal N}^{(j)}({\mathbb Z}_p)(w^{(j)}).$$Here ${\mathcal N}_w\cap {\mathcal N}^{(j)}=({\mathcal N}^{(j)})_{w^{(j)}}$ if we define $({\mathcal N}^{(j)})_{w^{(j)}}$ inside ${\mathcal N}^{(j)}$ by the same recipe as we defined ${\mathcal N}_w$ inside ${\mathcal N}$.

 For $i^{(j)}\in{\mathbb Z}$ with $\ell(w^{(j)})={\rm dim}(({\mathcal N}^{(j)})_{w^{(j)}}({\mathbb Z}_p))<i^{(j)}$ we have $H^{i^{(j)}}(({\mathcal N}^{(j)})_{w^{(j)}}({\mathbb Z}_p),V(\mu^{(j)}))=0$ because $({\mathcal N}^{(j)})_{w^{(j)}}({\mathbb Z}_p)$ is a torsion free ${\mathbb Q}_p$-analytic group. Thus, given $(i^{(j)})_{j\in J}$ and $w\in {\mathcal W}$ with $\ell(w)=\sum_j\ell(w^{(j)})=\sum_ji^{(j)}$ we have$$\bigotimes_{j\in J}H^{i^{(j)}}(({\mathcal N}^{(j)})_{w^{(j)}}({\mathbb Z}_p),V(\mu^{(j)}))=0\quad\quad\mbox{ unless }\ell(w^{(j)})=i^{(j)}\mbox{ for all }j.$$Therefore the K\"unneth formula gives a natural isomorphism$$\prod_{w\in {\mathcal W}\atop \ell(w)=i}H^i({\mathcal N}_w({\mathbb Z}_p),V(\mu)_{{\mathcal N}({\mathbb Z}_p)(w)})\cong \prod_{(i^{(j)})_{j\in J}\atop\sum_ji^{(j)}=i}\prod_{w\in {\mathcal W}\atop \ell(w^{(j)})=i^{(j)}}\bigotimes_{j\in J}H^{i^{(j)}}(({\mathcal N}^{(j)})_{w^{(j)}}({\mathbb Z}_p),V(\mu^{(j)})_{{\mathcal N}^{(j)}({\mathbb Z}_p)(w^{(j)})})$$(in the second product symbol on the right hand side we require $\ell(w^{(j)})=i^{(j)}$ for all $j$). On the other hand, the K\"unneth formula gives a natural isomorphism $$H^i({\mathcal N}({\mathbb Z}_p),V(\mu))\cong \prod_{(i^{(j)})_{j\in J}\atop\sum_ji^{(j)}=i}\bigotimes_{j\in J}H^{i^{(j)}}({\mathcal N}^{(j)}({\mathbb Z}_p),V(\mu^{(j)})).$$Together it follows that the map (\ref{prodpolotilgal}) is the direct sum, taken over all tuples $(i^{(j)})_{j\in J}$ with $\sum_ji^{(j)}=i$, of the maps\begin{align}\bigotimes_{j\in J}H^{i^{(j)}}({\mathcal N}^{(j)}({\mathbb Z}_p),V(\mu^{(j)}))&\longrightarrow \prod_{w\in {\mathcal W}\atop \ell(w^{(j)})=i^{(j)}}\bigotimes_{j\in J}H^{i^{(j)}}(({\mathcal N}^{(j)})_{w^{(j)}}({\mathbb Z}_p),V(\mu^{(j)})_{{\mathcal N}^{(j)}({\mathbb Z}_p)(w^{(j)})})\notag\\{}&\cong \bigotimes_{j\in J}(\prod_{w\in {\mathcal W}\atop \ell(w^{(j)})=i^{(j)}}H^{i^{(j)}}(({\mathcal N}^{(j)})_{w^{(j)}}({\mathbb Z}_p),V(\mu^{(j)})_{{\mathcal N}^{(j)}({\mathbb Z}_p)(w^{(j)})})).\notag\end{align}It is enough to show that all the tensor factors of these maps are bijective, i.e. that the natural maps $$H^{i^{(j)}}({\mathcal N}^{(j)}({\mathbb Z}_p),V(\mu^{(j)}))\longrightarrow \prod_{w\in {\mathcal W}\atop \ell(w^{(j)})=i^{(j)}}H^{i^{(j)}}(({\mathcal N}^{(j)})_{w^{(j)}}({\mathbb Z}_p),V(\mu^{(j)})_{{\mathcal N}^{(j)}({\mathbb Z}_p)(w^{(j)})})$$are bijective. But as the ${\mathcal G}_{{\mathbb F}_p}^{(j)}$ are simple this has been done in Theorem \ref{poloclou}. \hfill$\Box$\\ 

\begin{lem}\label{emerhilf} Let ${\mathfrak N}_2\subsetneq{\mathfrak N}_1$ be
  subgroups of ${\mathcal N}({\mathbb Q}_p)$ of the form ${\mathfrak
    N}_j=\prod_{\alpha\in\Psi}{\mathfrak
    N}_j\cap {\mathcal U}_{\alpha}({\mathbb
    Q}_p)$ with free rank one ${\mathbb Z}_p$-submodules ${\mathfrak
    N}_j\cap {\mathcal U}_{\alpha}({\mathbb
    Q}_p)$ of ${\mathcal U}_{\alpha}({\mathbb
    Q}_p)$, for a subset $\Psi$ of $\Phi^-$. Let $i=|\Psi|$ and let
  $k$ be a field of characteristic $p$. Then the
  restriction map$$H^i({\mathfrak N}_1,k)\longrightarrow H^i({\mathfrak
    N}_2,k)$$is the zero map.
\end{lem}

{\sc Proof:} We find a filtration
$$\emptyset=\Psi_0\subset\Psi_1\subset\ldots\subset\Psi_{i-1}\subset\Psi_i=\Psi$$
with $|\Psi_{s+1}-\Psi_{s}|=1$ and such that, if we put ${\mathfrak
    N}_{j,s}=\prod_{\alpha\in\Psi_s}{\mathfrak
    N}_j\cap {\mathcal U}_{\alpha}({\mathbb
    Q}_p)$, then ${\mathfrak
    N}_{j,s}$ is a normal subgroup in ${\mathfrak
    N}_{j,s+1}$ for all $1\le s<i$. As the ${\mathfrak
    N}_{j}$ are torsion free ${\mathbb Z}_p$-analytic groups of dimension $i$
  we have $H^s({\mathfrak N}_j,.)=0$ for any $s>i$. Therefore the Hochschild-Serre
  spectral sequences provide isomorphisms$$H^i({\mathfrak N}_j,k)\cong H^1({\mathfrak
    N}_{j,i}/{\mathfrak
    N}_{j,i-1},H^1({\mathfrak
    N}_{j,i-1}/{\mathfrak
    N}_{j,i-2},\ldots H^1({\mathfrak
    N}_{j,1},k)\ldots )).$$This reduces our claim to the statement that for a proper inclusion ${\mathfrak M}_2\subsetneq{\mathfrak M}_1$ of free rank one ${\mathbb Z}_p$-modules the restriction map$$H^1({\mathfrak M}_1,k)={\rm Hom}({\mathfrak M}_1,k)\longrightarrow H^1({\mathfrak M}_2,k)={\rm Hom}({\mathfrak M}_2,k)$$is the zero map. \hfill$\Box$\\

\section{ $M_{\chi}(V)$ is free for $F={\mathbb Q}_p$ and $p$-small weights}

\label{jache}

\begin{lem}\label{prodco} (i) For any split reductive group $\widehat{\mathcal G}$ over
  ${\mathbb Z}_p$ there is a product ${\mathcal G}=\prod_{j\in J}{\mathcal
    G}^{(j)}$ of split, simple and simply connected groups ${\mathcal G}^{(j)}$ over
  ${\mathbb Z}_p$ together with a morphism of algebraic
  groups$$\gamma:{\mathcal G}=\prod_{j\in J}{\mathcal
    G}^{(j)}\longrightarrow\widehat{\mathcal G}$$ 

(a) such that the pull back
  ${\mathcal T}=\gamma^{-1}\widehat{\mathcal T}$ of a maximal split
  torus $\widehat{\mathcal T}$ in $\widehat{\mathcal G}$ is a maximal split
  torus in ${\mathcal G}$, 

(b) such that $\gamma$ induces an
  isomorphism between the root systems of $\widehat{\mathcal G}$ (with respect to $\widehat{\mathcal T}$) and ${\mathcal G}$
 (with respect to ${\mathcal T}$), and

(c) such that $\gamma$ induces an
  isomorphism between the respective root subgroups, for all roots in $\Phi$.  

(ii) Let $V$ be an irreducible $\overline{\mathbb F}_p[\widehat{\mathcal G}({\mathbb F}_p)]$-module. The restriction of $V$ to ${\mathcal G}({\mathbb F}_p)$ uniquely extends to an irreducible rational ${\mathcal G}_{{\mathbb F}_p}$-representation with $p$-restricted highest weight. 
\end{lem}

{\sc Proof:} (i) Decompose the root system $\Phi$ of $(\widehat{\mathcal
  G},\widehat{\mathcal T})$ into its irreducible components,
$\Phi=\coprod_{j\in J}\Phi^{(j)}$. For $j\in J$ define $\widehat{\mathcal
    G}^{(j)}$ as the subgroup of $\widehat{\mathcal G}$ generated by the root
  subgroups corresponding to the roots in $\Phi^{(j)}$. This is a simple algebraic group over ${\mathbb Z}_p$. If ${\mathcal
    G}^{(j)}\to\widehat{\mathcal
    G}^{(j)}$ denotes its universal covering (an isogeny with simple simply connected ${\mathcal
    G}^{(j)}$) then the product $\gamma:{\mathcal G}=\prod_{j\in J}{\mathcal
    G}^{(j)}\to\widehat{\mathcal G}$ has the desired
  properties.

(ii) First we claim that $V$ is irreducible as a $\overline{\mathbb
  F}_p[{\mathcal G}({\mathbb F}_p)]$-module. Indeed, let ${\mathcal
  N}^-({\mathbb F}_p)$ and ${\mathcal N}^+({\mathbb F}_p)$ denote the groups
of ${\mathbb F}_p$-rational points of the unipotent radicals of two opposite
Borel subgroups in ${\mathcal G}_{{\mathbb F}_p}$. For any irreducible
$\overline{\mathbb F}_p[{\mathcal G}({\mathbb F}_p)]$-module the subspace of
${\mathcal N}^+({\mathbb F}_p)$-invariants (resp. the quotient space of
${{\mathcal N}^-({\mathbb F}_p)}$-coinvariants) is one dimensional. From this
one easily deduces that a $\overline{\mathbb F}_p[{\mathcal G}({\mathbb
  F}_p)]$-module $W$ is irreducible if and only if the following two
conditions hold true: ${\rm  dim}(W^{{\mathcal
    N}^+({\mathbb F}_p)})=1$ and the natural map $W^{{\mathcal N}^+({\mathbb
    F}_p)}\to W_{{\mathcal N}^-({\mathbb F}_p)}$ is bijective. The same is true for $\widehat{\mathcal G}({\mathbb F}_p)$ instead of ${\mathcal G}({\mathbb F}_p)$. As we may read ${\mathcal N}^-({\mathbb F}_p)$ and ${\mathcal N}^+({\mathbb F}_p)$ also as unipotent radicals of Borel subgroups in $\widehat{\mathcal G}({\mathbb F}_p)$ our claim follows. But now, as ${\mathcal G}$ is semisimple and simply connected, it follows from a theorem of Steinberg that $V$ uniquely extends to an irreducible rational ${\mathcal G}_{{\mathbb F}_p}$-representation with $p$-restricted highest weight, cf. e.g. \cite{herwe}, Appendix, Proposition 1.3.
\hfill$\Box$\\

If the derived group $\widehat{\mathcal G}^{\rm der}$ of $\widehat{\mathcal G}$ is simply connected then $\gamma$ is an isogeny onto $\widehat{\mathcal G}^{\rm der}$. For example, if $\widehat{\mathcal G}={\rm GL}_n$ then ${\mathcal G}={\rm SL}_n$ and this is simple, $|J|=1$.\\

We now return to our earlier setting and assume that $F={\mathbb Q}_p$ there. We
apply Lemma \ref{prodco} to $\widehat{\mathcal G}=\widetilde{\bf G}$ to obtain a morphism $\gamma:{\mathcal G}=\prod_{j\in J}{\mathcal
    G}^{(j)}\longrightarrow\widetilde{\bf G}$ as described in that Lemma, with
  maximal torus ${\mathcal T}=\gamma^{-1}\widetilde{\bf T}$ in ${\mathcal G}$. Let
  $\Phi^{(j)}$ for $j\in J$ denote the irreducible components of our root
  system $\Phi$. Our fixed choice of positive roots in $\Phi$ determines
  choices of positive roots $\Phi^{(j),+}$ in all $\Phi^{(j)}$. Let
  $\rho^{(j)}$ denote the half sum of the elements of $\Phi^{(j),+}$. Let $h^{(j)}$ denote the Coxeter number of the root
 system $\Phi^{(j)}$. Any $\mu\in
 X_1^*({\mathcal T})$ is of the form $\mu=\prod\mu^{(j)}$ with $\mu^{(j)}\in
 X_1^*({\mathcal T}\cap {\mathcal
    G}^{(j)})$.

\begin{satz}\label{modest} Assume that $F={\mathbb Q}_p$ and that we are in (at least) one of the
  following cases (i), (ii) or (iii):

(i) We are in case (I), $R$ is a discrete valuation ring with residue field $k$ of characteristic
$p$, and we have $p\ge h^{(j)}-1$ for all $j\in J$.

(ii) We are in case (II), and the unique $\mu\in X_1^*({\mathcal T})$
with $V\cong V(\mu)$ (as rational ${\mathcal G}_{{\mathbb
    F}_p}$-representations) satisfies $\langle\mu^{(j)}+\rho^{(j)},\check{\beta}\rangle\le
p$ for all $\beta\in\Phi^{(j),+}$ for all $j\in J$.

(iii) We are in case (III), and the unique $\mu\in X_1^*({\mathcal T})$
with $V\otimes_{{\mathcal O}_F}\overline{\mathbb F}_p\cong V(\mu)$ (as rational ${\mathcal G}_{{\mathbb
    F}_p}$-representations) satisfies $\langle\mu^{(j)}+\rho^{(j)},\check{\beta}\rangle\le
p$ for all $\beta\in\Phi^{(j),+}$ for all $j\in J$.

Then the $B$-module $M_{\chi}(V)$ is
  free and the Koszul complex (\ref{indkos}) is a $\widetilde{G}$-equivariant resolution of $M_{\chi}(V)$, for any $R$-algebra
  $B$ and any $R$-algebra
  homomorphism $\chi:{\mathfrak H}_V(\widetilde{G},\widetilde{K})\to B$ such
  that $\chi|_{{\mathfrak H}(V)}$ factors through $R$. 
\end{satz}

{\sc Proof:} By Proposition \ref{bettercond} it is enough to show that for any
$0\le i<|\nabla|$ the $k$-vector space $H^i({\mathfrak N},V\otimes_Rk)$ is
generated by elements $x$ with the following property: for any
$\lambda\in\nabla$ we have
$\rho_{\lambda}^i(x)=\theta_{\lambda,x}\iota_{\lambda}^i(x)$ for some
$\theta_{\lambda,x}\in k$, and moreover $\theta_{\lambda,x}\ne0$ for at least
one $\lambda$ (depending on $x$).

In case (i) the $k[{\mathfrak N}]$-module $V\otimes_Rk$ is the trivial
$k[{\mathfrak N}]$-module $k$ (as ${\mathfrak N}$ is a pro-$p$-group). Moreover, the
  condition $p\ge h^{(j)}-1$ is equivalent with the condition $\langle\rho^{(j)},\check{\beta}\rangle\le
p$ for all $\beta\in\Phi^{(j),+}$, i.e. the condition on $\mu$ in case (ii) if
$\mu=0$. Thus case (i) is reduced to case (ii). Similarly, case (iii) is obviously reduced to case (ii).

Therefore, it is enough to concentrate on the case (ii) and to prove the above
statement with $R=k$, $V\otimes_Rk=V$.

{\it Step 1:} Consider the natural projection map\begin{gather}H^i({{{{\mathfrak
        N}}}},V)\stackrel{\varphi}{\longrightarrow}\prod_{w\in W\atop
  \ell(w)=i}H^i({\mathfrak N}_w,V_{{\mathfrak N}(w)})\label{hiercoin}\end{gather}where ${\mathfrak
N}_w$ and ${{\mathfrak N}(w)}$ are defined as in section \ref{cohin}. We may
regard the subgroups ${{{{\mathfrak N}}}}$, ${\mathfrak N}(w)$ and
 ${\mathfrak N}_w$ of $\widetilde{K}=\widetilde{\bf G}({\mathbb Z}_p)$ also as
 subgroups of ${\mathcal G}({\mathbb Z}_p)$. Therefore Theorem
 \ref{prodpoloclou} implies that $\varphi$ is bijective. For $x\in H^i({{{{\mathfrak
        N}}}},V)$ write $\varphi(x)=\prod_w\varphi(x)_w$ with $\varphi(x)_w\in
H^i({\mathfrak N}_w,V_{{\mathfrak N}(w)})$ for $w\in W$ with
$\ell(w)=i$. Consider elements $x\in H^i({{{{\mathfrak
        N}}}},V)$ for which there is some $w(x)\in W$, $\ell(w(x))=i$, such that $\varphi(x)_w=0$ for all $w\ne
w(x)$ and such that $\varphi(x)_{w(x)}\ne0$. To prove the above
 statement it is enough to prove that for any such $x$ we have $\rho_{\lambda}^i(x)=\theta_{\lambda,x}\iota_{\lambda}^i(x)$ in $H^i({{{{\mathfrak
        N}}}}[\lambda],V)$ for some
$\theta_{\lambda,x}\in k$, and moreover $\theta_{\lambda,x}\ne0$ for at least
one $\lambda$ (depending on $x$).

{\it Step 2:} Here we introduce some more notation. Recall (cf. the definition of $\xi_{\lambda}$) that for $\lambda\in\nabla$ we have the maximal 
parabolic subgroup $\widetilde{\bf P}_{-\lambda}$ of $\widetilde{\bf G}$ with Levi
factor $\widetilde{\bf M}_{-\lambda}$ such that $\widetilde{\bf T}\subset\widetilde{\bf M}_{-\lambda}$. The Weyl group $W^{\lambda}$ of
$\widetilde{\bf M}_{-\lambda}$ can be identified with a subgroup of $W$: there is a simple reflection $s_{\lambda}$ in $W$ such that $W^{\lambda}$ is the subgroup of $W$ generated by the simple reflections different from $s_{\lambda}$. We put \begin{gather}{\mathfrak N}^{\lambda}=\prod_{\alpha\in \Phi^{-}\cap W^{\lambda}\Phi^{+}}\widetilde{K}\cap
  U_{\alpha}.\label{ilam}\end{gather}This is the group of ${\mathbb Z}_p$-valued points of the unipotent radical of a Borel subgroup in
$\widetilde{\bf M}_{-\lambda}$. Notice that ${\mathfrak N}^{\lambda}\subset{\mathfrak N}[\lambda]$. On the other hand we have $${\bf
      N}_{-\lambda}({\mathbb Z}_p)=\prod_{\alpha\in \Phi^--(\Phi^-\cap W^{\lambda}\Phi^{+})}\widetilde{K}\cap
  U_{\alpha}.$$Thus ${\bf
      N}_{-\lambda}({\mathbb Z}_p)$ is a normal subgroup in ${\mathfrak N}=\prod_{\alpha\in \Phi^{-}}\widetilde{K}\cap
  U_{\alpha}$ with \begin{gather}{\mathfrak N}^{\lambda}\cong {\mathfrak N}/{\bf
      N}_{-\lambda}({\mathbb Z}_p).\label{nsubquo}\end{gather}For $w\in
    W$ we have the equivalences (for the second one observe Lemma \ref{petersfrage})\begin{align}w\in
      W^{\lambda}\quad&\Leftrightarrow\quad{\mathfrak N}_w\subset{\mathfrak
        N}^{\lambda}\notag\\{}&\Leftrightarrow\quad {\mathfrak N}_w\subset{\mathfrak
        N}[\lambda].\label{lambdaw}\end{align}

{\it Step 3:} Let $\lambda\in \nabla$. In view of formula (\ref{nsubquo}) we may regard ${{{{\mathfrak N}}}}^{\lambda}$ as acting on $V^{{{\bf
     N}_{-\lambda}({\mathbb Z}_p)}}$ and consider the natural map \begin{gather}H^i({{{{\mathfrak N}}}}^{\lambda},V^{{{\bf
     N}_{-\lambda}({\mathbb Z}_p)}})\longrightarrow\prod_{w\in W^{\lambda}\atop \ell(w)=i}H^i({\mathfrak N}_w,(V^{{{\bf
     N}_{-\lambda}({\mathbb Z}_p)}})_{{\mathfrak
   N}(w)}).\label{lambdapolotilgal}\end{gather}Take $\gamma:{\mathcal M}_{-\lambda}=\prod_{j\in J_{-\lambda}}{\mathcal
    M}_{-\lambda}^{(j)}\longrightarrow\widetilde{\bf M}_{-\lambda}$ as in Lemma \ref{prodco} (applied to $\widehat{\mathcal G}=\widetilde{\bf M}_{-\lambda}$). By \cite{her} Lemma 2.5 the $\widetilde{\bf
 M}_{-\lambda,{\mathbb F}_p}({\mathbb F}_p)$-representation $V^{{{\bf
     N}_{-\lambda}({\mathbb Z}_p)}}$ is irreducible. As an irreducible rational ${\mathcal
 M}_{-\lambda,{\mathbb F}_p}$-representation with $p$-restricted highest
weight (Lemma \ref{prodco}), $V^{{{\bf
     N}_{-\lambda}({\mathbb Z}_p)}}$ satisfies the conditions of Theorem \ref{prodpoloclou}, as
is easily deduced from our hypotheses on $V$. Therefore Theorem
\ref{prodpoloclou} implies the bijectivity of the map (\ref{lambdapolotilgal}).

{\it Step 4:} Let $\lambda\in \nabla$, let \begin{gather}H^i({\mathfrak N}/{\bf
      N}_{-\lambda}({\mathbb Z}_p),V^{{\bf
      N}_{-\lambda}({\mathbb F}_p)})\longrightarrow H^i({\mathfrak
    N},V)\label{bildim}\end{gather}denote the map induced by the natural projection $\nu_{\lambda}:{\mathfrak N}\to{\mathfrak N}/{\bf
      N}_{-\lambda}({\mathbb Z}_p)$ and the inclusion $\omega_{\lambda}:V^{{\bf
      N}_{-\lambda}({\mathbb F}_p)}\to V$. Let
  $$\widetilde{\rho}_{\lambda}^i,\widetilde{\iota}_{\lambda}^i:H^i({\mathfrak N}/{\bf
      N}_{-\lambda}({\mathbb Z}_p),V^{{\bf
      N}_{-\lambda}({\mathbb F}_p)})\longrightarrow H^i({\mathfrak
    N}[\lambda],V)$$ denote the maps obtained by composing
  $\rho_{\lambda}^i(x)$ resp. $\iota_{\lambda}^i$ with the map
  (\ref{bildim}). We claim that
  $\widetilde{\rho}_{\lambda}^i=\widetilde{\iota}_{\lambda}^i$. 

Indeed,
  $\widetilde{\iota}_{\lambda}^i$ is given by the pair $(\nu_{\lambda}\circ
  (\gamma_{[\lambda]}(.)\gamma_{[\lambda]}^{-1}),\omega_{\lambda})$, while
  $\widetilde{\rho}_{\lambda}^i$ is given by the pair $(\nu_{\lambda}\circ
  \ell_{\lambda},\xi_{\lambda}\circ\omega_{\lambda})$ where $\ell_{\lambda}:{\mathfrak
    N}[\lambda]\to{\mathfrak
    N}$ is the inclusion.

 We clearly have
  $\omega_{\lambda}=\xi_{\lambda}\circ\omega_{\lambda}$. But we also have $\nu_{\lambda}\circ
  (\gamma_{[\lambda]}(.)\gamma_{[\lambda]}^{-1})=\nu_{\lambda}\circ \ell_{\lambda}$.
  To see this recall formula (\ref{nsubquo}) and observe that the conjugation map
  $\gamma_{[\lambda]}(.)\gamma_{[\lambda]}^{-1}$ of formula (\ref{jan}) acts on ${\mathfrak
    N}^{\lambda}$ as the identity and that it induces an isomorphism ${\bf
      N}_{-\lambda}({\mathbb Z}_p)\cap {\mathfrak N}[\lambda]\cong {\bf
      N}_{-\lambda}({\mathbb Z}_p)\cap {\mathfrak N}$.

{\it Step 5:} Let $w\in W$, let $\pi_w:V\to V_{{\mathfrak N}({w})}$ denote the
natural projection. We claim that the
composite $$V\stackrel{\xi_{\lambda}}{\longrightarrow}V\stackrel{\pi_w}{\longrightarrow}
V_{{\mathfrak N}({w})}$$ factors through $\pi_w$. More precisely, we claim that $\xi_{\lambda}({\rm Ker}(\pi_w))\subset {\rm Ker}(\pi_w)$.

It follows from its definition that $\xi_{\lambda}$ respects the action of
$\widetilde{\bf M}_{\lambda}({\mathcal O}_F)=\widetilde{\bf M}_{-\lambda}({\mathcal
  O}_F)$. In particular, $\xi_{\lambda}$ respects the action of the maximal torus $\widetilde{\bf T}_{{\mathbb
    F}_p}$ and hence it respects the weight spaces of $V$
with respect to the action of $\widetilde{\bf T}_{{\mathbb
    F}_p}$. (In fact it acts on $V$ even as a projector.) On the other hand, $\pi_w$ may be viewed as the projection onto the space of
lowest weight with respect to the Borel subgroup corresponding to ${\mathfrak
  N}({w})$, the kernel of $\pi_w$ is the direct sum of the remaining weight
spaces. Together this implies our claim.

{\it Step 6:} Take $x\in H^i({{{{\mathfrak
        N}}}},V)$ for which there is some $w(x)\in W$, $\ell(w(x))=i$, with $\varphi(x)_w=0$ for all $w\ne
w(x)$ and with $\varphi(x)_{w(x)}\ne0$. 

We first claim that for those $\lambda\in\nabla$ with $w(x)\notin
W^{\lambda}$ we have $\rho_{\lambda}^i(x)=0$, or
phrased differently,
$\rho_{\lambda}^i(x)=\theta_{\lambda,x}\iota_{\lambda}^i(x)$ with
$\theta_{\lambda,x}=0\in k$. 

Let $\tau:{\mathfrak N}\to{\mathfrak N}$ denote the composite ${\mathfrak N}\cong{\mathfrak N}[\lambda]\subset {\mathfrak N}$, the (inverse of) the conjugation isomorphism of formula (\ref{jan}) followed by the inclusion of ${\mathfrak N}[\lambda]$ into ${\mathfrak N}$. The map $\rho_{\lambda}^i:H^i({{{{\mathfrak
        N}}}},V)\to H^i({{{{\mathfrak
        N}}}}[\lambda],V)$ can then be identified with the map $H^i({\mathfrak
  N},V)\to H^i({\mathfrak N},V)$ induced by the pair
$(\tau,\xi_{\lambda})$. It follows from Step 5 (together with the obvious fact that
$\tau({\mathfrak N}_w)\subset {\mathfrak N}_w$ for all $w\in W$) that this map sits in a commutative diagram\begin{gather}\xymatrix{H^i({\mathfrak N},V)
   \ar[r]\ar[d]&H^i({{{{\mathfrak N}}}},V)\ar[d]\\ \prod_{w\in W\atop \ell(w)=i}H^i({\mathfrak N}_w,V_{{\mathfrak N}(w)})\ar[r]&
   \prod_{w\in W\atop \ell(w)=i}H^i({\mathfrak N}_w,V_{{\mathfrak
       N}(w)})}\label{jada}\end{gather}where the vertical arrows are the isomorphisms $\varphi$
 (Step 1), and the lower horizontal arrow is the product, over all $w\in W$ with $\ell(w)=i$, of maps \begin{gather}\beta_w^i:H^i({\mathfrak
  N}_{w},V_{{\mathfrak N}({w})})\longrightarrow H^i({\mathfrak
  N}_{w},V_{{\mathfrak N}({w})})\label{wres}\end{gather}induced by the pair $(\tau,\xi_{\lambda})$.  Therefore, by our assumptions on $x$ it is enough to show
that $\varphi(x)_{w(x)}$ vanishes under the map $\beta_{w(x)}^i$. As $V_{{\mathfrak N}({w(x)})}$ is the trivial representation of ${\mathfrak N}_{w(x)}$ this map can be identified with the map \begin{gather}H^i({\mathfrak
  N}_{w(x)},k)\longrightarrow H^i({\mathfrak
  N}_{w(x)},k)\label{wxres}\end{gather}induced by a $k$-endomorphism of $k$ and the injective group homomorphism $\tau_w:{\mathfrak
  N}_{w(x)}\to {\mathfrak
  N}_{w(x)}$ induced by $\tau$. Our assumption on
$\lambda$ and the equivalence (\ref{lambdaw}) imply that ${\mathfrak
  N}_{w(x)}$ is not contained in ${\mathfrak N}[\lambda]$ and hence that $\tau_w$ is not an isomorphism. It therefore follows from Lemma \ref{emerhilf} that the map
(\ref{wxres}) is the zero map. The claim is proven.

Next we claim that, on the other hand, for those $\lambda\in\nabla$ with $w(x)\in
W^{\lambda}$ we have $\rho_{\lambda}^i(x)=\pm\iota_{\lambda}^i(x)$, or phrased
differently, $\rho_{\lambda}^i(x)=\theta_{\lambda,x}\iota_{\lambda}^i(x)$ with
$\theta_{\lambda,x}=\pm1\in k$. For this it is enough, by Step 4, to show that $x$
lies in the image of the map (\ref{bildim}). Consider the commutative
 diagram$$\xymatrix{H^i({\mathfrak N}^{\lambda},V^{{{\bf
     N}_{-\lambda}({\mathbb Z}_p)}})
   \ar[r]\ar[d]&H^i({{{{\mathfrak N}}}},V)\ar[d]\\ \prod_{w\in W^{\lambda}\atop \ell(w)=i}H^i({\mathfrak N}_w,V_{{\mathfrak N}(w)})\ar[r]&
   \prod_{w\in W\atop \ell(w)=i}H^i({\mathfrak N}_w,V_{{\mathfrak
       N}(w)})}.$$Here the upper horizontal is the map\begin{gather}H^i({\mathfrak N}^{\lambda},V^{{{\bf
     N}_{-\lambda}({\mathbb Z}_p)}})\longrightarrow
   H^i(\frac{{{{{\mathfrak N}}}}}{{\bf
     N}_{-\lambda}({\mathbb Z}_p)},V^{{{\bf
     N}_{-\lambda}({\mathbb Z}_p)}})\longrightarrow
   H^i({{{{\mathfrak N}}}},V)\label{bluff}\end{gather}where the first arrow is the
 isomorphism induced by (the inverse of) the natural isomorphism of groups (\ref{nsubquo}). The lower horizontal arrow is the
obvious inclusion. The vertical arrow on the left hand side is induced by the inclusions ${\mathfrak N}_w\to {\mathfrak N}^{\lambda}$ and the projections$$V^{{{\bf
     N}_{-\lambda}({\mathbb Z}_p)}}\longrightarrow (V^{{{\bf
     N}_{-\lambda}({\mathbb Z}_p)}})_{{\mathfrak N}(w)}\cong (V_{{{\bf
     N}_{\lambda}({\mathbb Z}_p)}})_{{\mathfrak N}(w)}\cong V_{{\mathfrak N}(w)};$$here observe that ${{\bf
     N}_{\lambda}({\mathbb Z}_p)}\subset {\mathfrak N}(w)$. Thus, the vertical
 arrow on the left hand side can be identified with the isomorphism
 (\ref{lambdapolotilgal}) (Step 3), while the vertical
 arrow on the right hand side is the isomorphism $\varphi$
 (Step 1). It follows that $x$
lies in the image of the map (\ref{bildim}), as desired.

{\it Step 7:} It remains to show that for $x\in H^i({{{{\mathfrak
        N}}}},V)$ as in Step 6 there is indeed some $\lambda\in\nabla$ with $w(x)\in
W^{\lambda}$. But as we assume $\ell(w(x))=i<|\nabla|$ this follows from the general fact that for any $w\in W$ with $\ell(w)<|\nabla|$ there is some $\lambda\in\nabla$ with $w\in
W^{\lambda}$.\hfill$\Box$\\

{\it Remark:} Once more, the condition $\langle\mu^{(j)}+\rho^{(j)},\check{\beta}\rangle\le
p$ for all $\beta\in\Phi^{(j),+}$ is equivalent with saying that $\mu^{(j)}$
lies 'in the (closed) bottom
alcove' of $X_1^*({\mathcal T}\cap {\mathcal
    G}^{(j)})$. From the strong linkage principle of Andersen and Jantzen
  \cite{jant} it follows that this
  condition (for all $j$) implies that the irreducible rational
  $\widetilde{\bf G}_{F}$-module of highest weight $\mu$ admits a
  $\widetilde{K}$-stable ${\mathcal O}_F$-lattice such that its reduction is an irreducible rational
  $\widetilde{\bf G}_{k_F}$-module (necessarily also of highest weight
  $\mu$), as in case (III).

\section{Integral structures in locally algebraic representations}

Let $E/F/{\mathbb Q}_p$ be finite field exensions; $F$ will maintain its role
it played so far. A locally algebraic representation of ${\widetilde{G}}$ over
$E$ is a
tensor product $W=W^{\rm alg}\otimes_EW^{\rm sm}$ where $W^{\rm alg}$ is an
algebraic (or rational) representation of ${\widetilde{G}}$ over
$E$, and $W^{\rm sm}$ is a smooth representation of ${\widetilde{G}}$ over
$E$. It is of great interest (cf. e.g. \cite{st}) to decide when $W$ admits a
${\widetilde{G}}$-invariant norm, or equivalently, a
$\widetilde{G}$-stable ${\mathcal O}_E$-free submodule containing an $E$-basis
of $W$. 

Suppose that there is a finite dimensional
$E[\widetilde{K}]$-module $V_E$ and a character $\chi_E:{\mathfrak
    H}_{V_E}(\widetilde{G},\widetilde{K})\to E$ such that $$W\cong
  {\rm ind}_{\widetilde{K}}^{\widetilde{G}}V_E\otimes_{{\mathfrak
    H}_{V_E}(\widetilde{G},\widetilde{K}),\chi_E}E\,\,(\,\,=:M_{\chi_E}(V_E)\,\,).$$ 
Suppose that $V$ is a $\widetilde{K}$-stable ${\mathcal O}_E$-lattice in $V_E$. In order to hope for a $\widetilde{G}$-invariant norm on
$W$ we should demand that the restriction of $\chi_E:{\mathfrak
    H}_{V_E}(\widetilde{G},\widetilde{K})\to E$ to the ${\mathcal O}_E$-subalgebra ${\mathfrak
    H}_{V}(\widetilde{G},\widetilde{K})$ takes values in ${\mathcal O}_E$, i.e.
induces a character $\chi:{\mathfrak
    H}_{V}(\widetilde{G},\widetilde{K})\to {\mathcal O}_E$. For a geometric
  interpretation (in terms of the given data $W^{\rm alg}$ and $W^{\rm sm}$)
  of the set of all $\chi_E$ satisfying this condition see \cite{st}. For an
  equivalent characterization of such $\chi_E$ in terms of Jacquet functors
  see \cite{eme}. It seems to be conceivable that this condition is in fact
  also sufficient for the existence of a $\widetilde{G}$-invariant norm on $W$.

Given a $\chi$ as above, the existence of a
$\widetilde{G}$-invariant norm on $W$ immediately follows once we know that
the ${\mathcal O}_E$-module $${\rm ind}_{\widetilde{K}}^{\widetilde{G}}V\otimes_{{\mathfrak
    H}_{V}(\widetilde{G},\widetilde{K}),\chi}{\mathcal O}_E\,\,(\,\,=:M_{\chi}(V)\,\,)$$is
free --- because then $M_{\chi}(V)$ is a
$\widetilde{G}$-stable ${\mathcal O}_E$-free submodule of
$M_{\chi}(V)\otimes_{{\mathcal O}_E}E\cong M_{\chi_E}(V_E)\cong W$. For this
freeness we have given sufficient criteria, and verified them in
important cases, in the present
paper, e.g. Proposition \ref{bettercond}, Corollary \ref{rank2}, Theorem
\ref{modest}.\\

Let us work this out in the case of an unramified smooth principal series representation
$W=P_{\rho}={\rm Ind}_{\widetilde{T}N}^{\widetilde{G}}\rho$; for this discussion we faithfully follow \cite{st}. Let $\rho:\widetilde{T}\to E^{\times}$ be a smooth character and
consider the principal series representation$$P_{\rho}={\rm
  Ind}_{\widetilde{T}N}^{\widetilde{G}}\rho=\{f:\widetilde{G}\to E \mbox{ locally
  constant }\,|\,f(gtn)=\rho^{-1}(t)f(g)\mbox{ for }g\in \widetilde{G},
t\in\widetilde{T}, n\in N\}$$of $\widetilde{G}$ (which acts by left
translations). If $\rho$ in
fact takes values in ${\mathcal O}_E^{\times}$ then the existence of a
$\widetilde{G}$-invariant norm on $P_{\rho}$ is obvious. However, this last condition
is in general not a necessary one.

Assume that $\rho$ is unramified, i.e. that $\rho$ factors through
$\widetilde{T}\to\widetilde{T}/\widetilde{T}\cap  \widetilde{K}\cong
X_*(\widetilde{\bf T})$. Form ${\rm
  ind}_{\widetilde{K}}^{\widetilde{G}}E$ with respect to the trivial
$E[\widetilde{K}]$-module structure on $E$. Then the natural map $${\rm Hom}_{E[\widetilde{G}]}({\rm
  ind}_{\widetilde{K}}^{\widetilde{G}}E,P_{\rho})\cong {\rm
  Hom}_{E[\widetilde{K}]}(E,P_{\rho})\cong P_{\rho}^{\widetilde{K}}$$ provides an
action of ${\mathfrak
    H}_E(\widetilde{G},\widetilde{K})={\rm End}_{\widetilde{G}}({\rm
  ind}_{\widetilde{K}}^{\widetilde{G}}E)$ on the subspace $P_{\rho}^{\widetilde{K}}$ of
$\widetilde{K}$-invariants in $P_{\rho}$. From the Iwasawa decomposition $\widetilde{G}=\widetilde{K}\widetilde{T}N$ it
follows that $P_{\rho}^{\widetilde{K}}$ is one-dimensional, hence ${\mathfrak
    H}_E(\widetilde{G},\widetilde{K})$ acts on $P_{\rho}^{\widetilde{K}}$ through a character $\omega_{\rho}:{\mathfrak
    H}_E(\widetilde{G},\widetilde{K})\to E$. By \cite{kato} Lemma 2.4(i) this is simply the character\begin{gather}\omega_{\rho}=\rho\circ{\mathfrak S}\label{satcha}\end{gather} where now $\rho$ also denotes the
  character $E[X_*(\widetilde{\bf T})]\to E$ induced by $\rho$ and where ${\mathfrak S}:{\mathfrak
    H}_E(\widetilde{G},\widetilde{K})\to E[X_*(\widetilde{\bf T})]$ is the
  Satake map from Proposition \ref{prredsat}. Put$$H_{\rho}=M_{\omega_{\rho}}(E)={\rm
  ind}_{\widetilde{K}}^{\widetilde{G}}E\otimes_{{\mathfrak
    H}_V(\widetilde{G},\widetilde{K}),\omega_{\rho}}E.$$This is known to be an
admissible smooth and finitely generated $\widetilde{G}$-representation
\cite{kato} Theorem 2.7. We have the $\widetilde{G}$-equivariant
map\begin{gather}H_{\rho}\longrightarrow P_{\rho},\quad f\otimes 1\mapsto \sum_{g\in
  \widetilde{G}/\widetilde{K}}f(g)g({\bf 1})\label{cindind}\end{gather}where ${\bf 1}\in
P_{\rho}^{\widetilde{K}}$ denotes the unique $\widetilde{K}$-invariant function with ${\bf
  1}(1_{\widetilde{G}})=1_E$.

In \cite{kato} section 3 it is explained that the above map (\ref{cindind}) is
an isomorphism for sufficiently general $\rho$. (Note that, for example, if $P_{\rho}$ is irreducible --- if $\rho$ is regular, i.e. if $\rho$ is not fixed by any non-zero element in
$W$, a necessary and sufficient criterion for the irreducibility of $P_{\rho}$ is given in \cite{cas} Proposition 3.5 (b) --- then
it is enough to show that (\ref{cindind}) is injective, since it is obviously
non-zero. But to see this injectivity the results of \cite{dat} can be applied.)

If (\ref{cindind}) is an isomorphism, finding a $\widetilde{G}$-invariant norm on
$P_{\rho}$ is equivalent with finding a $\widetilde{G}$-invariant norm on
$H_{\rho}$. Thus, in view of the preceding discussion we ask whether the restriction of $\omega_{\rho}:{\mathfrak
    H}_E(\widetilde{G},\widetilde{K})\to E$ to the subalgebra ${\mathfrak
    H}_{{\mathcal O}_E}(\widetilde{G},\widetilde{K})$ takes values in ${\mathcal O}_E$, hence
induces a character $\omega_{\rho}^0:{\mathfrak
    H}_{{\mathcal O}_E}(\widetilde{G},\widetilde{K})\to {\mathcal O}_E$. By
  Proposition \ref{prredsat} and formula (\ref{satcha}) this is
equivalent with saying that the restriction of $\rho:E[X_*(\widetilde{\bf
  T})]\to E$ to the subalgebra ${\mathcal O}_E[X_*(\widetilde{\bf T})]\cap
E[X_*(\widetilde{\bf T})]^W={\mathcal O}_E[X_*(\widetilde{\bf T})]^W$ of
$E[X_*(\widetilde{\bf T})]$ takes values in ${\mathcal O}_E$.

Given this, the question is whether the ${\mathcal O}_E$-module $$M_{\omega_{\rho}^0}({\mathcal O}_E)={\rm
  ind}_{\widetilde{K}}^{\widetilde{G}}{\mathcal O}_E\otimes_{{\mathfrak
    H}_{{\mathcal O}_E}(\widetilde{G},\widetilde{K}),\omega_{\rho}^0}{\mathcal
  O}_E$$ is free. With Corollary \ref{rank2} and Theorem \ref{modest} we obtain:

\begin{satz}\label{normexis} Suppose that we are in at least one of the
  following two cases (1) or (2): 

(1) $F={\mathbb Q}_p$ and the
  Coxeter numbers of the irreducible components of the root system of
  $\widetilde{G}$ are at most $p+1$.

(2) $|\nabla|\le 2$.

 If the map (\ref{cindind}) is an
  isomorphism and if $\rho|_{{\mathcal O}_E[X_*(\widetilde{\bf T})]^W}$ takes values in ${\mathcal O}_E$ then the principal
  series representation $P_{\rho}$ admits a $\widetilde{G}$-invariant norm.
\end{satz}

The discussion of existence of a $\widetilde{G}$-invariant norm on $H_{\rho}$ is of course the same for {\it any} ${\rho}$, i.e. regardless of whether (\ref{cindind}) is an isomorphism. Then at least the following can be said on the
relationship with $P_{\rho}$. (We literally reproduce this from \cite{st}.) $H_{\rho}$ is generated by a single $\widetilde{K}$-invariant
vector, hence the same is true for all its quotients. Like ${\rm
  ind}_{\widetilde{K}}^{\widetilde{G}}E$ also $H_{\rho}$ has a one dimensional
subspace of $\widetilde{K}$-invariants. Together we
see that $H_{\rho}$ admits a unique $\widetilde{G}$-quotient $\overline{H}_{\rho}$ with a
non-zero (one dimensional)
subspace of $\widetilde{K}$-invariants --- the spherical representation for
$\rho$. Via the map (\ref{cindind}) it may be viewed as the unique irreducible
constituent of $P_{\rho}$ with a non-zero subspace of $\widetilde{K}$-invariants.

\section{Appendix}

\label{appco}

Here we compute the continuous group cohomology with values in ${\mathbb F}_p$
(and hence also with values in any field of characteristic $p$) of the group
${\mathfrak N}$ of upper triangular unipotent matrices in ${\rm
  SL}_3({\mathbb Z}_p)$. We give explicit formulae in terms of normalized
cochains. For example, our conventions are such that$$C_{N}^2({\mathfrak N},{\mathbb F}_p)=\{c:{\mathfrak N}\times {\mathfrak N}\longrightarrow
{\mathbb F}_p\,|\,c((i,1))=c((1,i))=0\mbox{ for all }i\in {\mathfrak N}\}$$is
the group of normalized $2$-cochains,$$Z_N^2({\mathfrak N},{\mathbb F}_p)=\{c\in
C_{N}^2({\mathfrak N},{\mathbb
  F}_p)\,|\,ic(j,k)-c(ij,k)+c(i,jk)-c(i,j)=0\}$$is its subgroup of normalized
$2$-cocycles, and a $2$-cocyle $c\in
Z_N^2({\mathfrak N},{\mathbb F}_p)$ defines the zero class in $H^2({\mathfrak N},{\mathbb F}_p)$ if and only if there is some
$b:{\mathfrak N}\to{\mathbb F}_p$ with $b(1)=0$ (i.e. a normalized $1$-cocycle) satisfying$$c(i,j)=ib(j)-b(ij)+b(i).$$(Of course, as ${\mathfrak N}$ acts trivially on
${\mathbb F}_p$ we could have written $c(j,k)$ resp. $b(j)$ instead of
$ic(j,k)$ resp. $ib(j)$.)\\

 For an element $a\in {\mathfrak N}$ and $1\le i<j\le 3$ we let
$a_{ij}\in{\mathbb Z}_p$ denote the $(i,j)$-entry of $a$. For example, for $a,
b\in {\mathfrak N}$ one computes $(ba)_{12}=a_{12}+b_{12}$, $(ba)_{23}=a_{23}+b_{23}$ and
$(ba)_{13}=a_{13}+b_{13}+b_{12}a_{23}$. For $a\in{\mathbb Z}_p$ we denote by
${a}\in{\mathbb F}_p$ also its residue class in ${\mathbb F}_p$. 

 Define the normal subgroups
\begin{align}{\mathfrak N}^{(12)}&=\{a\in {\mathfrak N} \,|\,a_{12}=0\},\notag\\{\mathfrak N}^{(23)}&=\{a\in {\mathfrak N}\,|\,a_{23}=0\}\notag\end{align} of ${\mathfrak N}$, and the non normal subgroups
\begin{align}{\mathfrak N}^{(12),(13)}&=\{a\in {\mathfrak N} \,|\,a_{12}=a_{13}=0\},\notag\\{\mathfrak N}^{(23),(13)}&=\{a\in {\mathfrak N} \,|\,a_{23}=a_{13}=0\}.\notag\end{align}

\begin{lem}\label{cocoico} $H^1({\mathfrak N},{\mathbb F}_p)$ is generated by the
classes of the $1$-cocycles $\alpha_{12}=[a\mapsto {{a}_{12}}]$ and
$\alpha_{23}=[a\mapsto {{a}_{23}}]$. The classes of the following $2$-cocylces
$\gamma_1$ and $\gamma_2$ form
a basis of $H^2({\mathfrak N},{\mathbb F}_p)$: $$\gamma_1=[(a,b)\mapsto{a_{13}b_{23}+\frac{1}{2}a_{12}(b_{23}^2-b_{23})}],$$$$\gamma_2=[(a,b)\mapsto{a_{12}b_{13}+\frac{1}{2}(a_{12}^2-a_{12})b_{23}}.]$$[Note
that these formulae make sense even if $p=2$ since then $a_{12}^2-a_{12}$ and
$b_{23}^2-b_{23}$ ly in $p{\mathbb Z}_p=2 {\mathbb Z}_2$. For $p\ne2$ the summands
$-\frac{1}{2}a_{12}b_{23}$ could be left out since they are coboundaries.] The cup product $$\cup:H^1({\mathfrak N},{\mathbb F}_p)\times H^1({\mathfrak
  N},{\mathbb F}_p)\to H^2({\mathfrak N},{\mathbb F}_p)$$ is the zero map. The natural maps$$H^1({\mathfrak N},{\mathbb F}_p)\longrightarrow
  H^1({\mathfrak N}^{(12),(13)},{\mathbb F}_p)\oplus H^1({\mathfrak N}^{(23),(13)},{\mathbb
    F}_p)$$$$H^2({\mathfrak N},{\mathbb F}_p)\longrightarrow H^2({\mathfrak N}^{(12)},{\mathbb
    F}_p)\oplus H^2({\mathfrak N}^{(23)},{\mathbb F}_p)$$are isomorphisms. A $3$-cocycle representing a
non-zero element (and hence a basis) in $H^3({\mathfrak N},{\mathbb F}_p)$ is obtained as the
cup product $\zeta=\gamma_1\cup \alpha_{12}$, i.e. the $3$-cocycle $\zeta:{\mathfrak N}\times {\mathfrak N}\times {\mathfrak N}\to{\mathbb F}_p$, $$\zeta=[(a,b,c)\mapsto\gamma_1(a,b)\cdot\alpha_{12}(c)={a_{13}b_{23}c_{12}+\frac{1}{2}a_{12}c_{12}(b_{23}^2-b_{23})}].$$

 \end{lem}

{\sc Proof:} As ${\mathfrak N}$ is a $3$-dimensional
free ${\mathbb Q}_p$-analytic group we have ${\dim}H^0({\mathfrak N},{\mathbb F}_p)=1={\rm
  dim}H^3({\mathfrak N},{\mathbb F}_p)$ and ${\dim}H^1({\mathfrak N},{\mathbb F}_p)={\rm dim}H^2({\mathfrak N},{\mathbb F}_p)$, while all
the other cohomology groups vanish. As ${\mathfrak N}$ is generated by the elements $a$ and $b$ with
$a_{12}=1$ and $a_{13}=a_{23}=0$ and $b_{23}=1$ and $b_{12}=b_{13}=0$ our
claims on $H^1({\mathfrak N},{\mathbb F}_p)={\rm Hom}({\mathfrak N},{\mathbb F}_p)$ follow easily. The restriction of $\gamma_1$
to ${\mathfrak N}^{(12)}$ represents a nonzero class in $H^2({\mathfrak N}^{(12)},{\mathbb F}_p)$, while the restriction of $\gamma_{2}$
to ${\mathfrak N}^{(12)}$ represents the zero class in $H^2({\mathfrak N}^{(12)},{\mathbb F}_p)$, and
dually: the restriction of $\gamma_2$
to ${\mathfrak N}^{(23)}$ represents a nonzero class in $H^2({\mathfrak N}^{(23)},{\mathbb F}_p)$, while the restriction of $\gamma_{1}$
to ${\mathfrak N}^{(23)}$ represents the zero class in $H^2({\mathfrak
  N}^{(23)},{\mathbb F}_p)$. This proves our claims on $H^1({\mathfrak N},{\mathbb F}_p)$. Finally our
claim on $H^3({\mathfrak N},{\mathbb F}_p)$: Since the isomorphism $H^3({\mathfrak N},{\mathbb
  F}_p)\cong H^1({\mathfrak N}/{\mathfrak N}^{(12)},H^2({\mathfrak N}^{(12)},{\mathbb F}_p))$ arising from the Hochschild Serre
spectral sequence sends the class of $\zeta$ to the class of $\gamma_1$ in
$$H^2({\mathfrak N}^{(12)},{\mathbb F}_p)\cong H^1({\mathfrak N}/{\mathfrak N}^{(12)},H^2({\mathfrak N}^{(12)},{\mathbb F}_p)),$$it follows that $\zeta$
indeed represents a non-zero class in $H^3({\mathfrak N},{\mathbb F}_p)$. 
\hfill$\Box$\\

\begin{flushleft} \textsc{Humboldt-Universit\"at zu Berlin\\Institut f\"ur Mathematik\\Rudower Chaussee 25\\12489 Berlin, Germany}\\ \textit{E-mail address}:
gkloenne@math.hu-berlin.de \end{flushleft} 
\begin{thebibliography}{abcdefgh} 
\bibitem{bali}{\it L. Barthel, R. Livn\'{e}}, Irreducible modular representations of ${\rm GL}_2$ of a local field,  Duke Math. J. {\bf 75},  no. 2, 261--292 (1994)

\bibitem{bo}{\it  J. Bella${{\ddot{\i}}}$che, A. Otwinowska}, Platitude du module universel pour ${\rm GL}_3$ en caract\'{e}ristique non banale, Bull. Soc. Math. France  {\bf 131},  no. 4, 507--525 (2003)

\bibitem{berbre} {\it L. Berger, C. Breuil}, Sur quelques repr\'{e}sentations potentiellement cristallines de ${\rm GL}_2(\mathbb Q_p)$, Ast\'{e}risque No. {\bf 330} (2010), 155–211.

\bibitem{bou}{\it Bourbaki, N.} Groupes et alg\`{e}bres de Lie, chap. 4 -- 6, 7 -- 9, Springer, Berlin Heidelberg (2008)

\bibitem{bs}{\it C. Breuil, P. Schneider}, First steps towards $p$-adic
  Langlands functoriality, J. reine angew. Math. {\bf 610}, 149 -- 180 (2007)

\bibitem{bro} {\it K. S. Brown}, Buildings, Berlin-Heidelberg-New York: Springer (1989) 

\bibitem{bt}{\it F. Bruhat, J. Tits}, Groupes r\'{e}ductives sur un corps
  local I, Publ. Math. IHES {\bf 41}, 5--251 (1972)

\bibitem{cas}{\it Casselman}, The unramified principal series of $p$-adic groups. I. The spherical function. Compositio Mathematica {\bf 40}, no 3, p. 387--406 (1980)

\bibitem{dat}{\it J. F. Dat}, Types et inductions pour les repr\'{e}sentations
  modulaires des groupes $p$-adiques, Ann. Sci. \'{E}cole Norm. Sup. (4) {\bf
    32}, no.1, 1--38 (1999)

\bibitem{dat1}{\it J. F. Dat}, Repr\'{e}sentations lisses $p$-temp\'{e}r\'{e}es
  des groupes $p$-adiques, Amer. J. Math. {\bf 131}, 227--255 (2009) 

\bibitem{dwy}{\it W. G. Dwyer}, Homology of Integral Upper-Triangular
  Matrices, Proc. AMS {\bf 94}, no. 3, 523--528 (1985)

\bibitem{eme}{\it M. Emerton}, $p$-adic $L$-functions and unitary completions of representations of $p$-adic reductive groups,  Duke Math. J. {\bf 130},  no. 2, 353--392 (2005)

\bibitem{friedpar}{\it E. M. Friedlander, B. J. Parshall}, Cohomology of
  Infinitesimal and Discrete Groups, Math. Ann. {\bf 273}, 353--374 (1986)

\bibitem{gruen}{\it L. Gr\"unenfelder}, On the homology of filtered and graded
  rings, J. Pure Appl. Algebra {\bf 14}, 21-37 (1979)

\bibitem{hv}{\it G. Henniart, M.F. Vign\'{e}ras}, The Satake isomorphism
  modulo $p$ with weight, preprint

\bibitem{herwe}{\it F. Herzig}, The weight in a Serre-type conjecture for tame $n$-dimensional Galois representations, Duke Math. J. {\bf 149}, 37--116 (2009)

\bibitem{her}{\it F. Herzig}, A Satake isomorphism in characteristic $p$,
  Compositio Math. {\bf 147}, no. 1, 263-283 (2011)

\bibitem{herzclas}{\it F. Herzig}, The classification of irreducible admissible mod $p$ representations of a $p$-adic $GL_n$, Invent. Math. {\bf 186}, no. 2, 373-434 (2011)

\bibitem{hum}{\it J. E. Humphreys}, Modular representations of finite groups of Lie type, London Mathematical Society Lecture Note Series {\bf 326}, Cambridge University Press, Cambridge (2006)

\bibitem{jant}{\it J. C. Jantzen}, Representation of Algebraic Groups,
  2nd. edition, Mathematical surveys and Monographs, vol. {\bf 107}, American
  Mathematical society (2003)

\bibitem{kato}{\it S. Kato}, On eigenspaces of the Hecke algebra with respect to a good maximal compact subgroup of a $p$-adic reductive group,
Math. Ann. {\bf 257}, no. 1, 1--7 (1981)

\bibitem{ln}{\it Z. Lin, D. K. Nakano}, Complexity for modules over finite Chevalley groups and classical Lie algebras, Invent. Math. {\bf 138}, 85--101 (1999) 

\bibitem{pt} {\it P. Polo, J. Tilouine}, Bernstein-Gelfand-Gelfand complexes and cohomology of nilpotent groups over  $\mathbb{Z}\sb {(p)}$ for representations with $p$-small weights, Ast\'{e}risque {\bf 280}, Cohomology of Siegel varieties, 97--135 (2002)

\bibitem{ss}{\it P. Schneider and U. Stuhler}, Representation theory and sheaves on the Bruhat-Tits building, Inst. Hautes \'{E}tudes Sci. Publ. Math. No. {\bf 85}, 97–191 (1997)

\bibitem{st}{\it P. Schneider, J. Teitelbaum}, Banach-Hecke algebras and $p$-adic Galois representations,  Doc. Math.,  Extra Vol., 631--684 (2006)

\bibitem{sergal}{\it J. P. Serre}, Galois Cohomology, Springer Monographs in Mathematics. Springer-Verlag, Berlin (2002)

\bibitem{spri}{\it T. A. Springer}, Linear Algebraic Groups, Progress in Mathematics {\bf 9}, Birkh\"auser, Boston Inc., Boston, MA, second edition (1998)

\bibitem{tits}{\it J. Tits}, Reductive groups over local fields, Proc. AMS
Symp. Pure Math. Vol. XXXIII, Part 1, pp. 29–69 (1979)

\bibitem{vigre}{\it University of Georgia VIGRE Algebra Group}, On Kostant's Theorem for Lie Algebra Cohomology, Contemp. Math. {\bf 478}, p. 39--60 (2009)

\bibitem{vig0}{\it M. F. Vign\'{e}ras}, Cohomology of sheaves on the building
  and $R$-representations, Inventiones Math. {\bf 127}, 349--373 (1997)

\bibitem{vig1}{\it M. F. Vign\'{e}ras}, Alg\'{e}bres de Hecke affines
  g\'{e}n\'{e}riques, Representation Theory {\bf 10}, 1–20 (2006)

\bibitem{vig}{\it M. F. Vign\'{e}ras}, A criterion for integral structures and coefficient systems on the tree of ${\rm PGL}(2,F)$, Pure Appl. Math. Q. {\bf 4} (2008), no. 4, Special Issue: In honor of Jean-Pierre Serre. Part 1, 1291–1316.




\end{thebibliography}
\end{document}